\documentclass{article}
\usepackage[utf8x]{inputenc}
\usepackage{amsmath}
\usepackage{amssymb, amsthm,mathrsfs}
\usepackage{graphicx,subcaption,tikz}
\usetikzlibrary{automata,positioning,calc,intersections,through,backgrounds,patterns,fit}

\begin{document}
\newtheorem{theorem}{Theorem}[section] \newtheorem*{theorem*}{Theorem}
\newtheorem{definition}[theorem]{Definition} \newtheorem*{definition*}{Definition}
\newtheorem{remark}[theorem]{Remark} \newtheorem*{remark*}{Remark}
\newtheorem{lemma}[theorem]{Lemma} \newtheorem{corollary}[theorem]{Corollary}
\newtheorem{property}[theorem]{Property} \newtheorem{proposition}[theorem]{Proposition}
\newenvironment{customdef}[1]
  {\renewcommand\thedefinition{#1}\definition}
  {\enddefinition}

\newcommand{\normone}[1]{\left\lvert #1 \right\rvert}
\newcommand{\norminf}[1]{\left\lVert #1 \right\rVert}

\newcommand{\alphabet}{\mathcal A} \newcommand{\largea}{\mathcal L}
\newcommand{\graph}{G} \newcommand{\subgraph}{F} \newcommand{\subgraphp}{H} \newcommand{\labels}{l} \newcommand{\vertices}{V} \newcommand{\edges}{E}
\newcommand{\simplex}{\Delta} \newcommand{\ppath}{{\gamma^*}} \newcommand{\psimplex}{\simplex^{\! *}}
\newcommand{\gasket}{\mathcal G}
\newcommand{\rauzymat}{M} \newcommand{\rauzymap}{T} \newcommand{\rauzyacc}{T_*}
\newcommand{\roof}{r} \newcommand{\roofacc}{{r_*}} \newcommand{\roofaccn}{r_*^{(n)}}
\newcommand{\parameter}{\simplex^\graph} \newcommand{\suspension}{\widehat{\simplex}^\graph} \newcommand{\subparameter}{\simplex(\subgraph)}
\newcommand{\leb}{\operatorname{Leb}}
\newcommand{\SL}{\operatorname{SL}}
\newcommand{\R}{\mathbb R} \newcommand{\N}{\mathbb N} \newcommand{\Q}{\mathbb Q} \newcommand{\Z}{\mathbb Z}
\newcommand{\Rp}{\R_{+}} \newcommand{\cone}{\Rp^{\alphabet}}
\newcommand{\ccone}{\mathcal C} \newcommand{\conen}{\Rp^{n}}
\newcommand{\supp}{\operatorname{supp}} \newcommand{\diam}{\operatorname{diam}}

\newcommand{\stopping}{\mathcal S} \newcommand{\ccomponent}{\mathscr C}
\newcommand{\jump}{\mathcal J} \newcommand{\pattern}{\mathcal E_\ppath}
\newcommand{\minmax}{\mathcal M} \newcommand{\win}{\mathcal W} \newcommand{\lose}{\mathscr {L}}
\newcommand{\vertex}{v} \newcommand{\vertexout}{\vertex_{\mathrm{out}}} \newcommand{\vertexcount}{\mathcal V}
\newcommand{\coding}{c}
\newcommand{\infpath}{{\overline{\pmb \gamma}}}
\newcommand{\finpath}{\gamma}
\newcommand{\winner}{W}
\newcommand{\loser}{L}

\newcommand{\states}{\mathcal S}

\newcommand{\rauzytype}{non-degenerating}
\newcommand{\Rauzytype}{Non-degenerating}

\newcommand{\unstability}{quick escape property}
\newcommand{\unstable}{quickly escaping}

\title{Dynamical properties of simplicial systems and continued fraction algorithms.}
\author{Charles Fougeron}
\date{}

\maketitle


\abstract{
	We propose a new point of view on multidimensional continued fraction algorithms inspired by Rauzy induction.
	The generic behavior of such an algorithm is described here as a random walk on a graph that we call simplicial system.
	These systems thus provide a family of interesting examples for random walks with memory recorded by a finite dimensional vector.\\

	We introduce a general criterion on these graphs that induces ergodicity together with a bundle of many other dynamical properties for the algorithms.
	In particular, after computing the representation of Brun, Selmer and Arnoux--Rauzy--Poincaré algorithm in this formalism, it provides a unified proof of ergodicity for these classical examples as well as new results such as uniqueness of the ergodic measure and the fact that it induces the measure of maximal entropy on a canonical suspension of the map.\\

	These objects also bring a new perspective to some fractal sets such as Rauzy gaskets.
	We show an explicit upper bound on Hausdorff dimensions of fractals described in this formalism as well as a construction of their measure of maximal entropy.
	This implies in particular that the Rauzy gasket in all dimensions has Hausdorff dimension strictly smaller than its ambient space, as well as sharper bounds on the dimension and an asymptotic result.
}

\newpage
\tableofcontents
\newpage

\section{Introduction}

To compute the best rational approximations of a real number $0 < x < 1$, one classically uses the continued fraction algorithm, also known as the Gauss algorithm.\\
It essentially uses the Gauss map
$$G : x \to \left\{\frac 1 x\right\}$$
which associates to any positive number the integer part of its inverse.
The Gauss algorithm then consists in associating to the number $x$ the sequence of positive integers $a_n := [1/G^{n-1}(x)]$ for $n \geq 1$.
The corresponding sequence of rational numbers
$$ \frac {p_n} {q_n} := \dfrac {1} {a_1 + \dfrac {1} {a_2 + \dfrac {1} {\  \dfrac{\ddots}{a_{n-1} + \dfrac {1} {a_n}}}}}$$
converges to $x$ as $n \to \infty$ and produces the best approximations of $x$ in the sense that for all integer $a, b > 0$, if $|bx-a| \le |q_nx-p_n|$ then $b \ge q_n$.\\

The attempt to generalize this property to simultaneous approximation of vectors by rational numbers --- together with other algebraic motivations on characterization of elements in finite extensions of $\Q$ --- has been the starting point of the theory of Multidimensional Continued Fraction algorithms (MCF).
Jacobi and Poincaré in the 19th century have suggested two different generalizations and a large variety of algorithms have been introduced ever since.
Surprisingly enough, even the question of convergence on each coordinate of a vector for these algorithms does not have a straightforward answer.
This fact is greatly illustrated by Nogueira \cite{Nogueira95} who has showed that the algorithm introduced by Poincaré does not converge for almost every vector.

For more than 30 years, a large community of mathematicians have been working on proving dynamical properties of MCF, such as convergence \cite{Fischer72}, \cite{Nogueira95}, \cite{BertheLabbe13}, as well as further dynamical properties like ergodicity \cite{Schweiger90}, \cite{MessaoudiNogueiraSchweiger09}, \cite{BruinFokkinkKraaikamp15}, construction of invariant measures \cite{ArnouxLabbe18}, \cite{ArnouxSchmidt17} and estimates on the speed of convergence through Lyapunov exponents \cite{Lagarias93}, \cite{Broise-AlamichelGuivarc'h01}, \cite{FougeronSkripchenko19}.\\

The Gauss algorithm is an accelerated version of the projectivization of the Euclidean algorithm on $(x, 1-x)$, by the map defined on the dimension 2 positive cone by
$$F: (x, y)\in\Rp^2 \to
\left\{
\begin{array}{ll}
	(x-y, y) &\text{if } x > y\\
	(x, y-x) &\text{if } x < y
\end{array}
\right. .
$$
MCF algorithms (e.g. every case in the survey \cite{Schweiger00}) can be described as an acceleration of such a map where we subtract some coordinates to others depending on the ordering of the coordinates.\\

Whereas MCF are usually defined (see for instance \cite{Schweiger90} or \cite{Lagarias93}) as iterates of a single map on a $n$-dimensional positive, another natural generalization of the Gauss algorithm is given by Rauzy--Veech induction on interval exchange maps which act on several copies of a positive cone associated to each vertex of a combinatorial graph called \textit{a Rauzy graph}.
This induction is fundamental in the field of Teichmüller dynamics and is a key tool for most of the dynamical results on translation surfaces and Teichmüller flow.
Let us mention some of the results in the field obtained by studying the dynamics of Rauzy--Veech induction: ergodicity of the Teichmüller flow \cite{Veech82} (also proved by \cite{Masur82} with different techniques), introduction of Lyapunov exponents on translation surfaces \cite{Zorich96}, existence and uniqueness of a measure of maximal entropy for the Teichmüller flow \cite{BufetovGurevich11} and exponential mixing \cite{AvilaGouezelYoccoz06}. See \cite{ForniMatheus14} for a nice survey about these results.\\

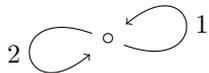
\begin{figure}[!htb]
	\center
	\begin{tikzpicture}[shorten >=1pt,node distance=2cm,on grid,auto]
	   \node (base) {$\circ$};

	   \path[->] (base) edge[loop, out=340, in =40, looseness=13]  node[right] {$1$} (base);
	   \path[->] (base) edge[loop, out=160, in =220, looseness=13]  node[left] {$2$} (base);
	\end{tikzpicture}
	\caption{The Rauzy graph of the Gauss algorithm.}
	\label{gauss:graph}
\end{figure}

The Gauss algorithm is the only example that belongs to both families of continued fraction algorithms and Rauzy--Veech inductions since its Rauzy, graph represented on Figure~\ref{gauss:graph}, has only one vertex.
Notice that our representation of Rauzy diagrams is slightly different from the classical representation where edges are labeled by the words \textit{top} or \textit{bottom} telling which of the top or bottom interval wins whereas we label edges by the corresponding losing letter.

Our initial motivation in the following work is to introduce a notion of Rauzy graph for a general MCF which will be a directed graph labeled by the index of the coordinates of the manipulated vectors.
A path in this graph will encode the combinatorial settings of each step of the MCF algorithm.
We want to take elementary steps in the graph and a given MCF algorithm will in general be an accelerated version of an algorithm defined by a graph.
Let us consider a vector in the positive cone on which the algorithm acts.
From a vertex in the graph we compare all the coordinates of the vector corresponding to the labels of edges going out.
A step of the algorithm is then defined by subtracting the smallest coordinate to all the larger ones and moving to the vertex toward which points the edge labeled by the letter on which the vector is the smallest.
This map can also be described as the inverse of a simple non-negative matrix which products decompose the more common descriptions of multidimensional continued fraction algorithms.
Consistently with Rauzy--Veech induction, we say the label of the smallest \textit{loses} and the other labels of outgoing edges \textit{win}.
This action on the vectors will then be called a \textit{win-lose induction}.
The definition of the graphs is closely related to the idea of \textit{simplicial systems} introduced in \cite{Kerckhoff85} to study unique ergodicity of interval exchange maps, we will thus give the same name to the corresponding labeled directed graphs.
\\

\begin{figure}[t!]
	\center
	\begin{tikzpicture}[shorten >=.5pt,node distance=1cm,on grid,auto, scale=.5]
	   \node (2>3) at (90:2){$\circ$};
	   \node (bot23) at (90:5) {$\bullet$};
	   \node (1>3) at (90:8) {$\circ$};

	   \node(3>1) at (210:2) {$\circ$};
	   \node (bot31) at (210:5) {$\bullet$};
	   \node (2>1) at (210:8) {$\circ$};

	   \node (1>2) at (330:2) {$\circ$};
	   \node (bot12) at (330:5) {$\bullet$};
	   \node (3>2) at (330:8) {$\circ$};

	   \path[->] (2>3) edge[bend right=40]  node[right] {$3$} (bot23);
	   \path[->] (bot23) edge[bend right=40] node[left] {$2$} (2>3);

	   \path[->] (1>3) edge[bend right=40]  node[left] {$3$} (bot23);
	   \path[->] (bot23) edge[bend right=40] node[right] {$1$} (1>3);

	   \path[->] (3>1) edge[bend right=40]  node[above] {$1$} (bot31);
	   \path[->] (bot31) edge[bend right=40] node[below] {$3$} (3>1);

	   \path[->] (2>1) edge[bend right=40]  node[below] {$1$} (bot31);
	   \path[->] (bot31) edge[bend right=40] node[above] {$2$} (2>1);

	   \path[->] (1>2) edge[bend right=40]  node[below] {$2$} (bot12);
	   \path[->] (bot12) edge[bend right=40] node[above] {$1$} (1>2);

	   \path[->] (3>2) edge[bend right=40]  node[above] {$2$} (bot12);
	   \path[->] (bot12) edge[bend right=40] node[below] {$3$} (3>2);

	   \path[->] (2>3) edge[bend right] node[left] {$1$} (3>1);
	   \path[->] (3>1) edge[bend right] node[below] {$2$} (1>2);
	   \path[->] (1>2) edge[bend right] node[right] {$3$} (2>3);

	   \path[->] (1>3) edge[bend left, in=120, out=60] node[right] {$2$} (3>2);
	   \path[->] (2>1) edge[bend left, in=120, out=60] node[left] {$3$} (1>3);
	   \path[->] (3>2) edge[bend left, in=120, out=60] node[below] {$1$} (2>1);
	\end{tikzpicture}
	\caption{Brun algorithm as a simplicial system.}
	\label{graph:brunintro}
\end{figure}
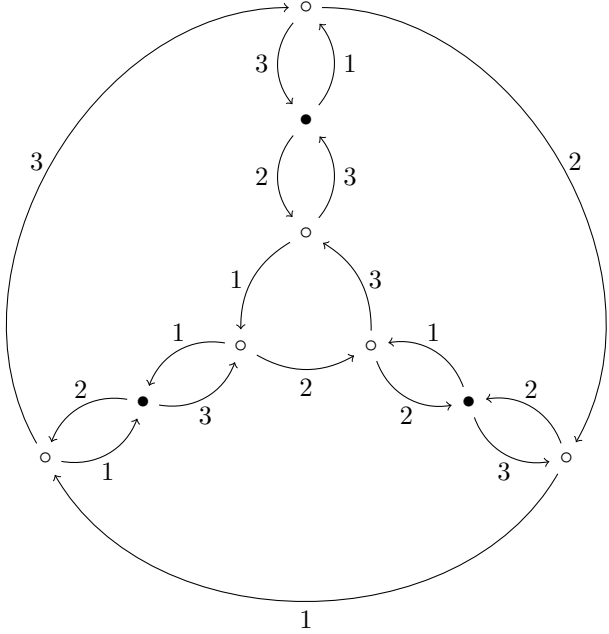

The graphs keep track of comparisons and subtractions on pairs of coordinates that are performed at each step of a given algorithm.
Intuitively, the appearance of several vertices in the graph and thus of several copies of the initial simplex are a consequence of the fact that the domains of definition of a MCF often depend on the relative order of more than two coordinates.
For instance the graph of Brun algorithm in dimension 3 is represented on Figure~\ref{graph:brunintro}.
A computation of this graph can be found in Section~\ref{brun}.

We claim that most of the classical MCF algorithm can be described in this language.
The reader can hopefully convince himself of this claim reading Section~\ref{examples} where several examples are translated in this setting.\\

After giving a general definition of a simplicial system in Section~\ref{definitions}, we introduce a property on the graph that induces ergodicity of the corresponding algorithm.
This property consists in showing that trajectories will go out of degenerate subgraphs in a small time with high probability, we will thus call it the \textit{\unstability\ } of the graph.
The degenerate subgraphs correspond to cases when a subset $\largea$ of the labels can be though as infinitesimally small compared to the others.
In that case, any time there is a comparison between labels in $\largea$ and its complementary, the trajectory will almost surely take the edge labeled in $\largea$.
Thus the degenerate subgraph associated to $\largea$ consists in removing these latter edges not labeled in $\largea$.\\

Our main theorem generalizes the results of \cite{Kerckhoff85}, \cite{BufetovGurevich11} and \cite{AvilaGouezelYoccoz06} to all \unstable\ simplicial systems.
\begin{theorem}
	\label{ergodicity}
	Every \unstable\ simplicial system has a unique ergodic measure equivalent to Lebesgue measure and it induces the unique invariant measure of maximal entropy on its canonical suspension.
\end{theorem}

\begin{figure}[t!]
  \center
  \hspace*{-1cm}
  \begin{subfigure}{0.5\textwidth}
	  \hspace*{-3cm}
	\begin{tikzpicture}[shorten >=.5pt,node distance=1cm,on grid,auto, scale=.5]
	   \node (2>3) at (90:2){$\circ$};
	   \node (bot23) at (90:5) {$\bullet$};
	   \node (1>3) at (90:8) {$\circ$};

	   \node(3>1) at (210:2) {$\circ$};
	   \node (bot31) at (210:5) {$\bullet$};
	   \node (2>1) at (210:8) {$\circ$};

	   \node (1>2) at (330:2) {$\circ$};
	   \node (bot12) at (330:5) {$\bullet$};
	   \node (3>2) at (330:8) {$\circ$};

	   \path[->] (bot23) edge[bend right=40] node[left] {$2$} (2>3);

	   \path[->] (bot23) edge[bend right=40] node[right] {$1$} (1>3);

	   \path[->] (3>1) edge[bend right=40]  node[above] {$1$} (bot31);

	   \path[->] (2>1) edge[bend right=40]  node[below] {$1$} (bot31);
	   \path[->] (bot31) edge[bend right=40] node[above] {$2$} (2>1);

	   \path[->] (1>2) edge[bend right=40]  node[below] {$2$} (bot12);
	   \path[->] (bot12) edge[bend right=40] node[above] {$1$} (1>2);

	   \path[->] (3>2) edge[bend right=40]  node[above] {$2$} (bot12);

	   \path[->] (2>3) edge[bend right] node[left] {$1$} (3>1);
	   \path[->] (3>1) edge[bend right] node[below] {$2$} (1>2);

	   \path[->] (1>3) edge[bend left, in=120, out=60] node[right] {$2$} (3>2);
	   \path[->] (3>2) edge[bend left, in=120, out=60] node[below] {$1$} (2>1);
	\end{tikzpicture}
  \end{subfigure}
  {\LARGE$\rightarrow$}%
  \hspace{1cm}
  \begin{subfigure}{0.2\textwidth}
	  \center
	\begin{tikzpicture}[shorten >=.5pt,node distance=1cm,on grid,auto, scale=.4]
	   \node (bbot12) at (210:5) {$\bullet$};
	   \node (2>1) at (210:8) {$\circ$};

	   \node (1>2) at (330:2) {$\circ$};
	   \node (bot12) at (330:5) {$\bullet$};

	   \path[->] (2>1) edge[bend right=40]  node[below] {$1$} (bbot12);
	   \path[->] (bbot12) edge[bend right=40] node[above] {$2$} (2>1);

	   \path[->] (1>2) edge[bend right=40]  node[below] {$2$} (bot12);
	   \path[->] (bot12) edge[bend right=40] node[above] {$1$} (1>2);
	\end{tikzpicture}
  \end{subfigure}
  \caption{Degenerate subgraph of Brun for $\largea = \{1,2\}$ on the left and its strongly connected components (with multiple vertices) on the right.}
	\label{degenerate:jacobi}
\end{figure}
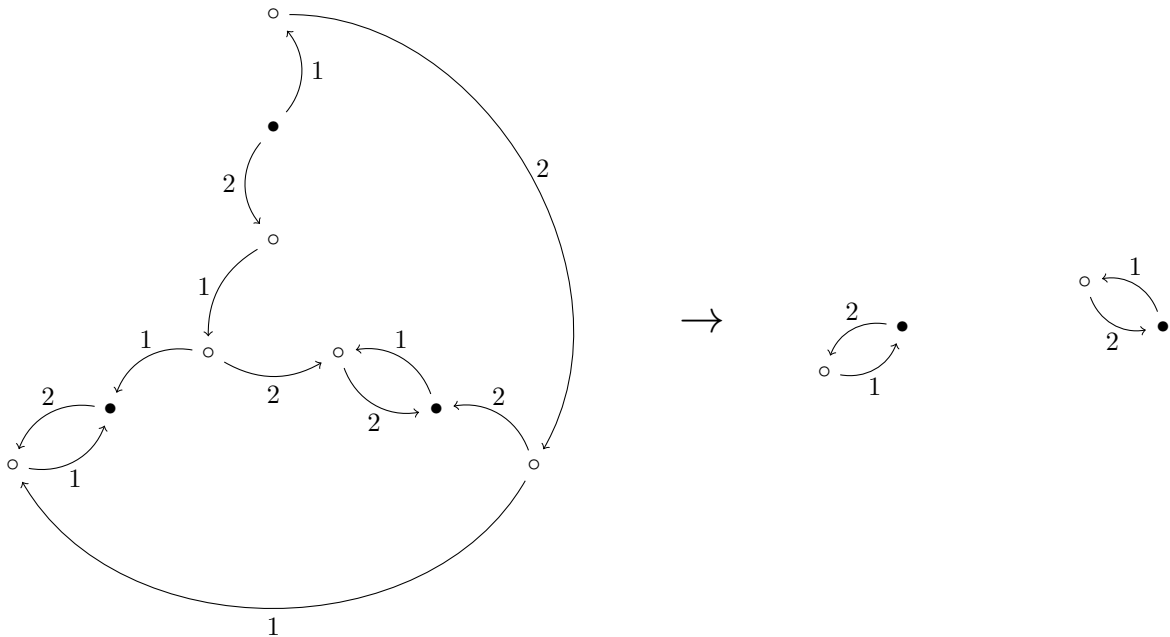

In the case of Brun algorithm the \unstability\ is easy to check.
Strongly connected components of degenerate subgraphs of Brun are always composed of a single loop on a vertex (see Figure~\ref{degenerate:jacobi}).
In other words, the \unstability\ reduces to showing that one letter cannot be the only one losing.
This cannot happen because the coordinates of a given vector are finite.\\

The proof of unique ergodicity of a generic interval exchange transformation of \cite{Kerckhoff85} can be interpreted as a proof of a weak form of the \unstability\ for Rauzy--Veech induction on interval exchanges.
This property has been proved in a different formalism and a stronger sense in \cite{AvilaGouezelYoccoz06} (see especially Appendix~A) which was the main inspiration for this work.

Two key properties on classical Rauzy graphs, noticed in \cite{Kerckhoff85} and later in \cite{ChaikaNogueira13}, are first that every letter has to lose infinitely often and second that in the degenerate case with labels in $\largea$, a letter in $\largea$ will alway lose eventually to a letter in the complementary set.
This latter property can be checked directly by considering the labeling of edges in the strongly connected components of all degenerate subgraphs.
We say a simplicial system satisfying these two properties is \textit{\rauzytype \ }and show the other main theorem of this work :
\begin{theorem}
	\label{unstable}
	\Rauzytype \ simplicial systems are \unstable.
\end{theorem}

Being \rauzytype \ is rather general but the first dynamical property is not always easy to check.
We then give an alternative first property that is purely graph theoretic which we call \rauzytype \ simplicial systems and show that theses two definitions are equivalent.\\

In Section~\ref{examples}, we explain a general strategy to associate a simplicial system to a MCF and show that for a large class of algorithms that these graphs are \rauzytype.
\begin{proposition}
	Brun and Selmer algorithms in all dimension and Arnoux--Rauzy--Poincaré algorithm in dimension 3 are simply connected and \rauzytype.
\end{proposition}
Which provides a unified proof of ergodicity of these algorithms as well as new results on existence and uniqueness of a measure of maximal entropy.
\begin{corollary}
	Brun, Selmer and Arnoux--Rauzy--Poincaré algorithms have a unique ergodic measure equivalent to Lebesgue measure and it induces the unique invariant measure of maximal entropy on its canonical suspension
\end{corollary}

This point of view may also bring a new perspective on Poincaré algorithm in all dimensions, which are the only examples of MCF which are not \rauzytype \ and for which it is not clear that they have stable degenerate subgraphs (except for the case of dimension $3$).
Studying ergodicity of Poincaré algorithm reduces in this formalism to compute fine estimates of the time a path in the graph stays in the degenerate subgraph.
Moreover, this formalism gives a lot of freedom to introduce new examples of ergodic MCF and find algorithms closer to optimality.\\

Another application of these simplicial systems is given by considering subsets of points $\simplex(\subgraph)$ in the simplex of parameter $\simplex$ above a vertex whose path for the win-lose induction remains in a subgraph $\subgraph \subset \graph$.
These sets are \textit{fractal sets} formed as a limit of union of subsimplices.

An important example of such sets is the \textit{Rauzy gasket}.
It has been primarily introduced by Levitt \cite{Levitt93} in connection with the dynamics of partially defined rotations of the circle, and was rediscovered latter by De Leo and Dynnikov \cite{DeLeoDynnikov09} to study particular class of examples for Novikov's conjecture in mathematical physics. It was generalized to all dimensions by \cite{ArnouxStarosta13} in a word combinatorics approach. More recently it was used in Diophantine geometry \cite{GamburdMageeRonan19} to show estimates on the number of integer points on Markoff--Hurwitz varieties.

Using thermodynamical formalism, we give a construction of a measure of maximal entropy on such fractal set and a bound on its Hausdorff dimension.
\begin{theorem}
	If $\subgraph$ is a \unstable\ strict subgraph of $\graph$ and has the same number of distinct labels on its edges then the Hausdorff dimension of $\simplex(\subgraph)$ is strictly smaller than the dimension of its ambient space $\simplex$.
\end{theorem}
As a consequence we generalize the result on the Hausdorff dimension of the Rauzy gasket in \cite{AvilaHubertSkripchenko16b} to Rauzy gaskets of arbitrary dimensions, as introduced in \cite{ArnouxStarosta13}.
\begin{corollary}
	The Rauzy gasket in all dimension has Hausdorff dimension strictly smaller than its ambient space.
\end{corollary}
Moreover, we obtain an explicit bound on the dimension in terms of the solution of an equation on the pressure on a family of potentials.
All the terms concerning thermodynamical formalism will be defined in Section~\ref{thermo}.
\begin{theorem}
	\label{hausdorff}
	Let $\phi$ be the geometric potential and $\kappa$ the unique positive real number satisfying $P(-\kappa \phi) = 0$.
	The Hausdorff dimension of the subset satisfies
	$$\dim_H \left(\simplex(\subgraph)\right) \le d - 1 + \frac {\kappa} {d+1} $$
	where $d$ is the dimension of $\simplex$.
\end{theorem}
It has been proved in \cite{GamburdMageeRonan19} that the solution of the above equation on the pressure for Rauzy gasket, $\gasket^d$, in the simplex of dimension $d$, is equal to the number $\alpha(d-1)$ estimated in \cite{Baragar98} (we detail the correspondence in Section~\ref{gasket}).
Using these estimates we have as a consequence of the previous theorem
\begin{align*}
	\dim_H(\gasket^2) &< 1.825 \\
	\dim_H(\gasket^3) &< 2.7 \\
	\dim_H(\gasket^4) &< 3.612 \\
\end{align*}

Notice that $\gasket^2$ is the classical Rauzy gasket studied in \cite{Levitt93}, \cite{DeLeoDynnikov09} and \cite{AvilaHubertSkripchenko16b}.
The only known bound for its Hausdorff dimension was given by \cite{AvilaHubertSkripchenko16b} where it was proved that $\dim_H(\gasket^2) < 2$.
Moreover, numerical experiments performed in \cite{DeLeoDynnikov09} seem to indicate that $\dim_H(\gasket^2)$ is in the range $[1.7, 1.8]$.\\

And in general,
\begin{equation*}
	\dim_H(\gasket^d) < d-1 + \frac {\log d} {\log 2 \cdot (d+1)} + o(d^{-1.58}).
\end{equation*}
This estimate implies that the difference between the dimension of the ambient space and the Hausdorff dimension of the Rauzy gasket $\gasket^d$  is asymptotically at least $1 - O(\log d / d)$.

\section{Definitions}
\label{definitions}
\subsection{Simplicial systems}
\label{SS}

Let $\graph = (\vertices, \edges)$ be a graph labeled on an alphabet $\alphabet$ by a function $\labels : \edges \to \alphabet$ such that for all $\vertex \in \vertices$  the restriction of $\labels$ to edges starting at $\vertex$ is injective.
We write $e : \vertex \to \vertex'$ if an edge $e$ goes from vertices $\vertex$ to $\vertex'$.\\

Let $\Rp := \{ x \in \R \mid x > 0\}$ and let us consider the norm on $\cone$ defined by $\normone {\lambda} = \sum_{\alpha \in \alphabet} \lambda_\alpha$.
Let $\simplex := \{\lambda \in \cone \mid \normone{\lambda} = 1\}$ be a simplex of dimension $|\alphabet| - 1$.
We associate to graph $\graph$ as above a piecewise projective map
\[ \rauzymap: \parameter \to \parameter, \]
on the parameter space $\parameter := \vertices \times \simplex$.\\

Let $\vertexout$ be the set of all edges going out of $\vertex$.
We defined a partition of $\simplex$ labeled by all $e \in \vertex_{\mathrm{out}}$ with the tiles
$$\simplex^e := \left\{ (\lambda_\alpha)_{\alpha \in \alphabet} \in \simplex \mid \lambda_{\labels(e)} < \lambda_\alpha \text{ for all } \alpha \in \labels(\vertexout) \text{ and } \alpha \ne \labels(e) \right\}.$$
The Rauzy matrix associated to this edge is
$$\rauzymat_e := \operatorname{Id} \ + \sum_{\substack{\alpha \in \labels\left(\vertexout\right)\\ \alpha \neq e}} E_{\alpha, \labels(e)}.$$
Where $E_{a,b}$ is the elementary matrix with coefficient $1$ at row $a$ and column $b$.\\

This implies a partition of $\simplex$.
Let $T: \parameter \to \parameter$, such that for all $\lambda \in \simplex^{e}$ with $e : \vertex \to \vertex'$,
\[\rauzymap(\vertex, \lambda) = \left(\vertex', \rauzymap_e(\lambda)\right),\]
where
\[ \rauzymap_e :
	\left\{
	\begin{array}{lll}
		\simplex^{e} & \to & \simplex \\
		\lambda      & \mapsto & \dfrac {\rauzymat_e^{-1} \lambda}{\normone {\rauzymat_e^{-1} \lambda}}
	\end{array}
	\right. .
\]

\begin{figure}[h!]
  \center
  \begin{subfigure}{0.3\textwidth}
	  \begin{tikzpicture}[scale=.5]
	       \coordinate (A) at (0,0);
	       \coordinate (Ac) at (2,0);
	       \coordinate (B) at (4,0);
	       \draw (A) node[left] {$a$} -- (Ac) -- (B) node[right] {$b$} --++ (120:2)coordinate(Bc) --++ (120:2)coordinate(C) node[above] {$c$} --++ (240:2)coordinate(Cc) -- cycle;
	       \draw[name path=A--Bc, opacity=0] (A) -- (Bc);
	       \draw[name path=B--Cc, opacity=0] (B) -- (Cc);
	       \draw[name path=Ac--C] (Ac) -- (C);
	       \path [name intersections={of=A--Bc and Ac--C,by=E}];
	       \draw (A) -- (B);
	       \draw (B) -- (C);
	       \draw (A) -- (C);
	       \path[pattern=north west lines,pattern color=black] (B)--(Ac)--(C)--cycle;
         \end{tikzpicture}
  \end{subfigure}
  {\LARGE$\xrightarrow{\rauzymap_{a}}$}%
  \hspace{1cm}
  \begin{subfigure}{0.25\textwidth}
	  \begin{tikzpicture}[scale=.5]
	       \coordinate (A) at (0,0);
	       \coordinate (Ac) at (2,0);
	       \coordinate (B) at (4,0);
	       \draw (A) node[left] {$a$} -- (Ac) -- (B) node[right] {$b$} --++ (120:2)coordinate(Bc) --++ (120:2)coordinate(C) node[above] {$c$} --++ (240:2)coordinate(Cc) -- cycle;
	       \draw[name path=A--Bc, opacity=0] (A) -- (Bc);
	       \draw[name path=B--Cc, opacity=0] (B) -- (Cc);
	       \draw[name path=Ac--C, opacity=0] (Ac) -- (C);
	       \path [name intersections={of=A--Bc and Ac--C,by=E}];
	       \draw (A) -- (B);
	       \draw (B) -- (C);
	       \draw (A) -- (C);
	       \path[pattern=north west lines,pattern color=black] (A)--(B)--(C)--cycle;
          \end{tikzpicture}
  \end{subfigure}
  \begin{subfigure}{0.3\textwidth}
	  \begin{tikzpicture}[scale=.5]
	       \coordinate (A) at (0,0);
	       \coordinate (Ac) at (2,0);
	       \coordinate (B) at (4,0);
	       \draw (A) node[left] {$a$} -- (Ac) -- (B) node[right] {$b$} --++ (120:2)coordinate(Bc) --++ (120:2)coordinate(C) node[above] {$c$} --++ (240:2)coordinate(Cc) -- cycle;
	       \draw[name path=A--Bc, opacity=0] (A) -- (Bc);
	       \draw[name path=B--Cc, opacity=0] (B) -- (Cc);
	       \draw[name path=Ac--C, opacity=0] (Ac) -- (C);
	       \path [name intersections={of=A--Bc and Ac--C,by=E}];
	       \draw (A) -- (E);
	       \draw (B) -- (E);
	       \draw (E) -- (C);
	       \draw (A) -- (B);
	       \draw (B) -- (C);
	       \draw (A) -- (C);
	       \path[pattern=north west lines,pattern color=black] (B)--(E)--(C)--cycle;
         \end{tikzpicture}
  \end{subfigure}
  {\LARGE$\xrightarrow{\rauzymap_{a}}$}%
  \hspace{1cm}
  \begin{subfigure}{0.25\textwidth}
	  \begin{tikzpicture}[scale=.5]
	       \coordinate (A) at (0,0);
	       \coordinate (Ac) at (2,0);
	       \coordinate (B) at (4,0);
	       \draw (A) node[left] {$a$} -- (Ac) -- (B) node[right] {$b$} --++ (120:2)coordinate(Bc) --++ (120:2)coordinate(C) node[above] {$c$} --++ (240:2)coordinate(Cc) -- cycle;
	       \draw[name path=A--Bc, opacity=0] (A) -- (Bc);
	       \draw[name path=B--Cc, opacity=0] (B) -- (Cc);
	       \draw[name path=Ac--C, opacity=0] (Ac) -- (C);
	       \path [name intersections={of=A--Bc and Ac--C,by=E}];
	       \draw (A) -- (B);
	       \draw (B) -- (C);
	       \draw (A) -- (C);
	       \path[pattern=north west lines,pattern color=black] (A)--(B)--(C)--cycle;
          \end{tikzpicture}
  \end{subfigure}
  \caption{Action of $\rauzymap_a$ on $\simplex_a$ when $\vertex$ has two or three outgoing edges.}
  \label{action:simplex}
\end{figure}
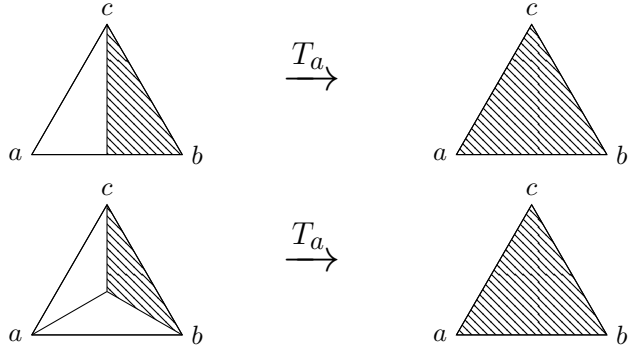
We call the graph $\graph$ a \textit{simplicial system} and the map $T$ its associated \textit{win-lose induction}.\\

\begin{remark}
	\label{homogeneous}
	The map $T_e$ is a projectivized version of linear maps on the cones,
	\[ \widetilde \rauzymap_e :
		\left\{
		\begin{array}{lll}
			\Rp \cdot \simplex^{e} & \to & \cone \\
			\lambda 	     & \mapsto & \rauzymat_e^{-1} \lambda
		\end{array}
		\right. .
	\]
	Similarly, we have a map $\widetilde \rauzymap$ on $\vertices \times \cone$, which will be useful when we will consider suspensions of $T$.\\
\end{remark}

Some dynamical properties of this linear map were studied for a more restrictive generalization of Rauzy--Veech induction in \cite{ChaikaNogueira13} which applies to Selmer and Jacobi--Perron algorithms.
They show an ergodicity property of the linear map with respect to Lebesgue measure which is not an invariant measure.\\

If there is a point in the graph that has no outgoing vertices, the map is not defined and the induction stops.
Such a vertex will be called a \textit{hole}.

The maps we have introduced are not well defined on the boundaries of $\simplex^{e}$.
Our dynamical study will focus in the first place on simplicial systems with no holes and the restriction of the corresponding maps to points for which the induction is defined at all times.
This is the complementary set of countably many codimension one subsets and thus a full Lebesgue measure set.\\

For a vertex $\vertex \in \vertices$, let $\Pi^\infty(\vertex)$ be the set of infinite paths starting at $\vertex$ in $\graph$ and $\Pi^h(\vertex)$ the set of finite path starting at $\vertex$ and ending in a hole.
The win-lose induction induces an injection (excluding a countable union of subset of codimension one)
$$\coding_\vertex: \simplex \to \Pi^h(\vertex) \sqcup \Pi^\infty(\vertex)$$
which associates to $\lambda \in \simplex$ the path followed by the first coordinate of win-lose induction $\rauzymap^n(\vertex, \lambda)$ in the graph until it stops.
If the simplicial system has no hole, then $\Pi^h(\vertex) = \emptyset$.
In the presence of holes, the set of parameters on a given vertex $\vertex \in \vertices$ for which the induction never stops $c_\vertex^{-1}(\Pi^\infty(\vertex))$ will be studied in Section~\ref{fractals} together with other restrictions of paths to a subgraph.\\

In analogy with the standard Rauzy induction on interval exchange transformations (see \cite{Yoccoz10} for an introduction to the subject), we define a loser and winers labels for each edge in the graph.
\begin{definition*}
	At a given vertex $\vertex$ with two or more outgoing edges, we say a letter $\alpha \in \alphabet$ \textbf{loses} along an edge $e$ going out of $\vertex$ if $\labels(e) = \alpha$.
	On the contrary, we say a letter $\beta$ \textbf{wins} along an edge $e$ based at $\vertex$ if there exists another edge $e'$ going out of $\vertex$ such that $\labels(e) \neq \labels(e') = \beta$.
	In both cases we say that $\alpha$ and $\beta$ \textbf{play} along the edge $e$.
\end{definition*}

One can describe the linear Rauzy map as the map which compares the coordinates of all the edges going out of a given vertex $\vertex$ and subtract the smallest to the others, in other terms subtract the losing coordinate to the winning ones.\\

Let $\Pi(\vertex)$ be the set of finite paths in $\graph$ starting at $\vertex$ and not ending in a hole.
For all $\finpath \in \Pi(\vertex)$ we denote the product of matrices
\[\rauzymat_\finpath := \rauzymat_{e_1} \dots \rauzymat_{e_n}\]
and the subsimplex of $\simplex$
\[ \simplex^\finpath := \rauzymat_\finpath \simplex.\]

\begin{proposition}
	For all $\vertex \in \vertices$ and all $\finpath \in \Pi(\vertex)$
	$$\coding_\vertex(\simplex^\finpath) = \Pi(\finpath).$$
	This corresponds to points which associated path starts with $\finpath$.
\end{proposition}

\begin{remark*}
	To clarify the redaction we will use by convention variables of the form $\finpath$ for finite paths in $\Pi(v)$ and $\infpath$ for infinite paths or paths ending in a hole \textit{i.e.} in $\Pi^h(v) \sqcup \Pi^\infty(v)$.
\end{remark*}

\subsection{Projective measures}
\label{proj}
In this subsection we introduce a key concept to study the Lebesgue generic dynamical behavior of a win-lose induction.
The idea will be to study its behavior in a stable family of measures equivalent to Lebesgue.
The feature that enables us to state dynamical results for Lebesgue generic paths is that the action of the induction on this family of measures is tractable through a dual action on a positive vector.

\begin{definition*}
	Let $q \in \Rp^\alphabet$, let $\nu_q$ be the Borel measure on the projective space $P\Rp^\alphabet$, such that for any subset $A \subset P\cone$,
	$$\nu_q (A) := \leb(\Rp A \cap \Lambda_q)$$
	where
	$\Lambda_q = \left\{v \in \cone \, \middle| \, \langle q,v \rangle < 1\right\}$.\\
\end{definition*}
We make the abuse of writing $\nu_q(\simplex)$ for some $\simplex \subset \cone$ while meaning $\nu_q(\Rp \simplex)$.
Moreover, $q$ is a line vector.\\

A fundamental equality is given by
\begin{align}
	\nonumber
	\nu_q (\rauzymat_\finpath \cdot \simplex) = \leb(\rauzymat_\finpath \simplex \cap \Lambda_q)
				      &= \leb(\simplex \cap \Lambda_{q \rauzymat_\finpath}) \\
				      &= \nu_{q \rauzymat_\finpath} (\simplex).
	\label{formula:M}
\end{align}
The vector $q$ keeps track of the way the measure is changed along the induction, we call it the \textit{distortion vector}.\\

An other fundamental equation comes from a computation that can be found in \cite{Veech78} Formula (5.5).
\begin{proposition}[Veech]
	For $\vertex \in \vertices$, $\finpath \in \Pi(\vertex)$ and $q \in \R_{>0}^\alphabet$,
	\begin{align}
		\label{formula:vol}
		&\nu_q (\simplex^{\finpath}) =  \frac 1 {n!}
		\cdot \frac {1}  {(q \rauzymat_\finpath)_1 \dots (q \rauzymat_\finpath)_n}
	\end{align}
	\label{prop:volumes}
\end{proposition}

If $\vertex$ is a vertex and $\finpath$ is a path starting at $\vertex$ we define the probability measure,
\[ \mathbb P_q^\vertex (\finpath) = \frac {\nu_q(\simplex^{\finpath})}{\nu_q(\simplex)}. \]
According to Formula (\ref{formula:vol}),
\[ \mathbb P_q^\vertex (\finpath) = \frac {N(q)}{N(q \rauzymat_\finpath)} \]
where $N(q) = \prod_{a \in \alphabet} q_a$.
\begin{proposition}
	Let $e \in \edges$ such that the label $\labels(e) = \alpha$, then
	\[ \mathbb P_q^\vertex(e) = \frac {q_{\alpha}} {(q \rauzymat_e)_{\alpha}}.\]
\end{proposition}
\begin{proof}
	Just notice that for all $\beta \neq \labels(e)$, $(q \rauzymat_e)_\beta = q_\beta$.
\end{proof}

\begin{definition}
	On a simplicial system we associate a probability spaces to any vertex $\vertex \in \vertices$ and distortion vector $q \in \cone$ formed by
	\begin{itemize}
		\item the sample space $\Pi^h(v) \sqcup \Pi^\infty(v)$,
		\item the set of events given by the $\sigma$-algebra generated by the sets
			$$\Pi^\infty(\finpath_s) := \left\{ \finpath_s \cdot \infpath_e \mid \infpath_e \in \Pi^\infty(\finpath_s \cdot v) \right\}$$
			for $\finpath_s \in \Pi(\vertex)$,
		\item the probability law given by the pushed forward of $\mathbb P_q^v$ by $\coding_\vertex$.
	\end{itemize}
\end{definition}

If a path $\finpath \in \Pi(\vertex)$ can be decomposed $\finpath = \finpath_s \cdot \finpath_e$ where $\finpath_s \in \Pi(\vertex)$ ends at $\vertex'$ and $\finpath_e$ is a path in $\Pi(\vertex')$, one can define conditional probabilities using Formula (\ref{formula:M})
\begin{equation}
	\label{condition}
	 \mathbb P_q^\vertex(\finpath \mid \finpath_s) = \dfrac {\nu_q(\simplex^{\finpath})}{\nu_q(\simplex^{\finpath_s})} =
	 \dfrac {\nu_q(\rauzymat_{\finpath_s} \Delta^{\finpath_e})}{\nu_q(\rauzymat_{\finpath_s} \Delta)} = \mathbb P_{q \rauzymat_{\finpath_s}}^{\vertex'}(\finpath_e).
\end{equation}

If $\Gamma_s$ is disjoint and $\finpath_s \cdot v$ denotes the ending vertex of the path $\finpath_s$ we can decompose the probability
\begin{equation}
	\label{formula:cond}
	\mathbb P_q^\vertex\left(\Gamma \cap \Pi\left(\Gamma_s\right)\right) = \sum_{\finpath_s \in \Gamma_s} \mathbb P_q^{\finpath_s \cdot \vertex} \left(\Gamma \cap \Pi(\finpath_s) \mid \finpath_s\right) \cdot \mathbb P_q^\vertex(\finpath_s)
\end{equation}
where $\Pi(\finpath_s) := \{ \finpath \cdot \finpath_e \mid \finpath_e \in \Pi(\finpath_s \cdot \vertex) \}$ and $\Pi(\Gamma_s) := \bigcup_{\finpath_s \in \Gamma_s} \Pi(\finpath_s)$.

\subsection{Stopping times}
The probability law of the paths strongly depends on the distortion.
Nevertheless some comparison of \textit{stopping times} will have upper or lower bounds independent of the distortion.
This will be a key tool to show Theorem~\ref{unstable} by induction.\\

Let $\mathcal P$ be a property on finite paths, we introduce a random variable
$$ T_{\mathcal P}:
\left\{
\begin{array}{lll}
	\Pi^\infty(\vertex) & \to &\N \cup \{ \infty \}\\
	\infpath      	    & \mapsto  & T_{\mathcal P}(\infpath)
\end{array}
\right.
$$
where $T_{\mathcal P}(\infpath) = \min \{ n \ge 0 \mid \infpath_n \text{ satisfies } \mathcal P \}$ and $\infpath_n$ the prefix of length $n$ of $\infpath$.
We make the abuse of writing the property $\mathcal P$ instead of $T_{\mathcal P}$.\\

The stochastic process formed by the sequence $(\infpath_n)_{n \ge 0}$ is a random path in the graph such that the law for each step only depends on the past and more precisely on the distortion vector.
For this reason, we believe this should be thought as a random walk on $\graph$ which has infinite memory recorded by a finite dimensional vector.
Hence the map $T_{\mathcal P}$ can be referred to as a stopping time for this random walk where the path stops when the property $\mathcal P$ is satisfied.\\

\begin{proposition}
	\label{path_condition}
	For every path $\finpath_s \in \Pi(\vertex)$ and $\infpath_e \in \Pi^\infty(\finpath_s \cdot v)$
	$$T_{\mathcal P}(\finpath_s \cdot \infpath_e) = T_{\finpath_s^{-1} \mathcal P}(\infpath_e)$$
	where $\infpath_e$ satisfies $\finpath_s^{-1} \mathcal P$ iff $\finpath_s \cdot \infpath_e$ satisfies $\mathcal P$.
\end{proposition}

\subsection{Suspension semi-flow}
\label{suspension}
Given a measurable function $f : \parameter \to \Rp$, one can define a suspension of the parameter space,
$$\suspension_f := (\parameter \times \R)/\sim$$
where we use the equivalence relation $(x, t) \sim \left(\rauzymap x, t + f(x)\right)$.
On $\suspension_f$ we define a suspension semi-flow
$$\phi_t : (x, s) \to (x, s + t).$$
Notice that these semi-flows are defined such that the first return map to the section $\parameter \times \{0\}$ is $\rauzymap$ and its return time is $f$.\\

In the case of simplicial systems, there is a canonical suspension coming from the fact that the space
$$\left(\vertices \times \cone\right)/\sim$$
where we identify $(x, s) \sim \widetilde \rauzymap (x, s)$, with $\widetilde T$ the homogeneous win-lose induction defined in Remark~\ref{homogeneous}.
This suspension has a natural semi-flow given for $t \geq 0$ by
$$\psi_t : (x, s) \to (x, e^t \cdot s).$$

The first return map to the section $\parameter \times \{ 0 \}$ for this semi-flow is also equal to $T$ and its first return time is given by a function $r : \parameter \to \Rp$ defined for $e : \vertex \to \vertex'$ by
$$r_e : (\lambda, \vertex) \to - \log \normone{ \rauzymat_e^{-1} \lambda }.$$
Where the norm is the $L^1$ norm.
This semi-flow is the suspension semi-flow on $\suspension_r$ that we call the \textit{canonical suspension semi-flow} associated to a simplicial system.
The function $r$ will be called the \textit{roof function} for the simplicial system.\\

The canonical suspension will play a central role in the study of entropy properties for win-lose inductions.
An important of its features comes from the fact that it is invariant by taking an acceleration of the map induced by a first return to a smaller subsimplex.
A lot of dynamical properties will be deduced from the study of the roof function using thermodynamic formalism in Section~\ref{thermo}.

\section{Generic dynamics in simplicial systems}
\label{properties}

In this section we study generic paths with respect to Lebesgue measure in a simplicial system.
We start by considering a measure associated in the family described above to an arbitrary distortion vector.
Pushing forward the measure after each step of the win-lose induction we obtain another measure in the family for another distortion vector.

All measures with a balanced distortion vector are equivalent with the same constant.
Thus our goal is to show that, starting from any distortion vector, the pushed forward measure will have a balanced one after a reasonable number of steps.
In consequence, the effects of the past on the random walk will be tamed periodically.

These random walk techniques will only be used in this section.
They will be the key to prove that a class of simplicial systems admits a uniformly expanding induced map to some subsimplex and an exponential tail property on its return time.

\subsection{Quick escape property}
In this subsection we introduce a property on random walks induced by the win-lose induction that implies a useful balancing property of the distortion vector and the appearance of any pattern almost surely infinitely many times in a generic path.
This will be the property we will aim to show afterwards starting from a graph theoretic criterion.
In other words the quick escape property is our keystone from the graph to the generic dynamics of the paths of the win-lose induction.\\

We start by introducing some useful properties for which we will compare the stopping times.
Let $\largea \subset \alphabet$, $\tau > 0, K > 0$ and $\ppath$ be a finite path in $\graph$. \\

Let $\jump^\tau$ be the property of a finite path $\finpath$ along which the distortion vector has \textit{jumped} by a factor $\tau$ \textit{i.e.} it satisfies
$$\max q \rauzymat_\finpath \ge \tau \max q$$
where the maximum is taken on all the coordinates of the vector.

\begin{proposition}
	\label{jumpinf}
	Assume that $\graph$ is strongly connected and has not all vertices with only one outgoing edge.
	For all $\tau > 0$ and all $q \in \cone$,
	$$\mathbb P_q(\jump^\tau = \infty) =0.$$
\end{proposition}

\begin{proof}
	For a given finite path $\finpath$ in the graph, let $n$ be the number of times it passes through a vertex of degree strictly larger than $1$.
	Then $$\max q \rauzymat_\finpath  \ge n \min q$$ and $n$ goes to infinity as the length of the path grows.
\end{proof}

An important property to consider on the distortion vector is the balance between its coordinates given by the following definition.
\begin{definition*}
	For $\largea \subset \alphabet$ and $K>1$, we say a distortion vector $q \in \cone$ is $(\largea, K)$-balanced if and only if $$\max_\alphabet q < K \min_\largea q.$$
	We say it is $K$-balanced in the case $\largea = \alphabet$.
\end{definition*}

This definition will be very useful due to the fact that it implies a lower bound on the probability that an edge labeled in $\largea$ is taken.
\begin{proposition}
	\label{balanced-bound}
	Let $\vertex \in \vertices$ be a vertex of degree at least $2$ and $q$ be a $(\largea, K)$-balanced distortion vector, then for all edge $e$ labeled in $\largea$ going out of $\vertex$
	$$\mathbb P_q^\vertex(e) \ge \frac {1} {|\alphabet| \cdot K}$$
	and for all finite path $\finpath$ of length $n$ in $\Pi(\vertex)$ with all labels in $\largea$
	$$\mathbb P_q^\vertex(\finpath) \ge \left(\frac {1} {|\alphabet| \cdot 2^n \cdot K}\right)^{n}.$$
\end{proposition}

Let $\stopping_\largea$ be the property of a finite path $\finpath$ for which
\[\max_{\alphabet \setminus \largea} q \rauzymat_\finpath \ge \min_\largea q \rauzymat_\finpath.\]
The stopping time corresponds to when the distortion on a coordinate outside of the subset $\largea$ reaches the size of the initial distortion on coordinates in $\largea$.

\begin{remark}
	\label{stopping}
	If we start with a $(\largea, K)$-balanced distortion vector, the set $\largea$ form the largest $|\largea|$ coordinates.
	In particular if there is a letter in $\largea$ which wins against a letter in its complementary set, it implies the event $\stopping_\largea$.\\
\end{remark}

\begin{definition*}[Quick escape property]
	We say a simplicial system is \unstable\ if for all non-empty subset $\largea \subsetneq \alphabet$ and all $K>1$ there exist $\tau > 1$ and $\delta > 0$ such that for all vertex $\vertex \in \vertices$ and all $(\largea, K)$-balanced distortion vector $q \in \cone$
\[
	\mathbb P^\vertex_q( \stopping_\largea \le \jump^\tau ) > \delta.
\]\\
\end{definition*}

Let $\minmax_{\largea}$ be the property of a finite path $\finpath$ for which
\[\min_\largea q \rauzymat_\finpath \ge \max_\alphabet q.\]
We denote by $\minmax$ the case $\minmax_\alphabet$.

\begin{lemma}
	\label{minmax}
	If a simplicial system is \unstable\ then there exists $\tau > 1$ and $\delta > 0$ such that for all vertex $\vertex \in \vertices$ and all distortion vector $q \in \cone$
	\begin{equation}
		\mathbb P^\vertex_q ( \minmax \le \jump^\tau) > \delta.
	\end{equation}
\end{lemma}

\begin{proof}
	We show by recurrence on $n$ that there exists $\tau_n > 1$ and $\delta_n > 0$ such that, for all vertex $\vertex$ and all distortion vector $q$, there exists a subset $\largea \subset \alphabet$ of cardinal $n$ which satisfies
	\[
		\mathbb P^\vertex_q( \minmax_\largea \le \jump^{\tau_n}) > \delta_n.
	\]

	\paragraph{Initialization.}
	For $n = 1$, we just have to take $\largea$ to be the singleton of the largest coordinate of $q$.

	\paragraph{Induction.}
	Assume that the property is true for some $n \ge 1$.
	With probability larger than $\delta_n$, the distortion vector will satisfy $\min_\largea(q \rauzymat_\finpath) \ge \max_\alphabet q$ and $\max_\alphabet (q \rauzymat_\finpath) \leq \tau_n \cdot \max_\alphabet q$, in particular it is $(\largea, \tau_n)$-balanced.

	Using the chain rule in Formula (\ref{formula:cond}), we only need to show the induction property with such a distortion vector.
	The \unstability\ tells us that there exists $\tau>1$ and $\delta>0$ such that, with probability larger than $\delta$, there is a letter $\alpha$ in $\largea$ and a letter $\beta$ outside of this set such that $(q \rauzymat_\finpath)_\beta \ge (q \rauzymat_\finpath)_\alpha$ before $q$ jumps by $\tau$.
	But then $(q \rauzymat_\finpath)_\beta \ge (q \rauzymat_\finpath)_\alpha \ge \max_\largea q$ and obviously $(q \rauzymat_\finpath)_\beta \ge q_\beta$ thus $\minmax_{\largea \cup \{\beta\}}$ is satisfied before $\jump^{\tau_n \cdot \tau}$.
\end{proof}

Let $\pattern$ be the property of a finite path which admits $\ppath$ as a suffix.
In other terms, the path $\finpath$ satisfies $\pattern$ if we can factor it into $\finpath = \finpath_0 \cdot \ppath$ for some finite path $\finpath_0$.\\

\begin{corollary}
	If a simplicial system is \unstable\ and strongly connected then for all finite path $\ppath \in \Pi(\vertex_0)$ there exists $K > 1$ and $\delta > 0$ such that for all $q \in \cone$ and all $\vertex \in \vertices$
	\[
		\mathbb P^\vertex_q (\pattern \le \jump^K) > \delta.
	\]
\end{corollary}

\begin{proof}
	Let $\tau$ as in Lemma~\ref{minmax} and let $\infpath \in \Pi(\vertex)$ satisfying $\minmax(\infpath) \le \jump^\tau(\infpath) =: n$.
	Let $\Gamma$ the set of such finite paths $\infpath_n$.
	By the strong connectivity hypothesis there exists a path of minimal length $\finpath_0 \in \Pi(\infpath_n \cdot \vertex)$ from the vertex $\vertex' := \infpath_n \cdot \vertex$ to the vertex $\vertex_0$.
	Notice that
	$$\min q \rauzymat_{\infpath_n} \geq \max q > \frac {1} {\tau} \cdot \max q \rauzymat_{\infpath_n}$$
	thus the distortion vector $q' := \rauzymat_{\infpath_n} q$ is $\tau$-balanced and, using Proposition~\ref{balanced-bound}, there exists $\delta'$ such that
	$$\mathbb P_{q'}^{\vertex'} (\finpath_0 \cdot \ppath) > \delta'.$$
	For $K$ large enough, for every $\vertex'$, the chosen path $\finpath_0 \cdot \ppath$ does not jump of a factor $K$ hence
	$$\mathbb P_{q'}^{\vertex'} \left(\pattern \le \jump^K\right) \ge \mathbb P_{q'}^{\vertex'} (\finpath_0 \cdot \ppath).$$
	Moreover $\finpath^{-1} \pattern = \pattern$ and $\finpath^{-1} \jump^{\tau K} \le \jump^{K}$.
	For the constant $\delta$ induced by Lemma~\ref{minmax} the chain rule implies
	\begin{align*}
		\mathbb P^\vertex_q \left(\pattern \le \jump^{\tau K}\right) & \ge \sum_{\finpath \in \Gamma} \mathbb P_q^\vertex \left(\finpath\right) \cdot \mathbb P_q^\vertex \left(\pattern \le \jump^{\tau K} \mid \finpath\right) \\
							   & \ge \sum_{\finpath \in \Gamma} \mathbb P_q^\vertex \left(\finpath\right) \cdot \mathbb P_{q'}^{\vertex'} \left(\pattern \le \jump^K\right) \\
							   & \ge \delta \cdot \delta'.
	\end{align*}
\end{proof}

\begin{remark}
	\label{recursion}
	By Property~\ref{jumpinf}, the distortion vector $q$ jumps almost surely in finite time for every $\tau > 1$.
	For any choice of finite path $\ppath$ the previous corollary then implies by induction that $\ppath$ appears almost surely in the coding.
	In particular one can define a first return map for the win-lose induction to the subsimplex of parameters whose path starts with $\ppath$.
	This is what is done in subsection~\ref{acceleration} where we show that this acceleration of the algorithm is uniformly expanding which implies ergodicity of the acceleration and of the initial induction.\\

	This sequence of results enable us to derive that some uniformly expanding acceleration of $\rauzymap$ is well defined Lebesgue almost everywhere.
	Such an acceleration will be the starting point for dynamical results and distortion will be of no use in the dynamical study of the algorithm.
\end{remark}

The following estimate is a discrete version of an exponential tail property that will be essential to apply thermodynamic formalism in the next section.
It implies that return times for the aforementioned acceleration of the induction can be though of as bounded.

Although we state the result for arbitrary distortion, we will only use it for one specific vector, namely for $\mathbf 1 := (1, \dots, 1)$ at each coordinates.
Similarly to the previous remark distortion was a key tool to derive such a structural result but will be of no further use in the thermodynamical study.
All results will then be stated for one Lebesgue measure $\nu := \nu_{\mathbf 1}$.
\def \exptail {\eta}
\begin{corollary}
	\label{proba:tail}
	If a simplicial system is \unstable\ and strongly connected then for any path $\ppath$ there exists $C>1, \exptail > 0$ such that for all $\vertex \in \vertices$, $\tau > 1$ and all $q \in \cone$
	$$\mathbb P^\vertex_q (\jump^\tau \leq \pattern) < C \cdot \tau^{-\exptail}.$$
\end{corollary}

\begin{proof}
	Consider $K$ and $\delta$ as in the previous corollary.
	For all $\tau > 1$, if $\tau > K^n$ for some integer $n$, we have, using the chain rule and Proposition~\ref{jumpinf},
	$$\mathbb P^\vertex_q (\jump^\tau \leq \pattern) \leq
	  \mathbb P^\vertex_q (\jump^{K^n} \leq \pattern) < (1-\delta)^n.$$
	  Thus, taking $n = \left\lceil \frac {\log \tau} {\log K} \right\rceil$,
	  $$(1-\delta)^n \le (1-\delta)^{\frac {\log \tau} {\log K} -1} = \frac {1} {1-\delta} \cdot \exp\left(\log(1-\delta) \cdot \frac{\log \tau}{\log K}\right).$$
	  Hence, for $C = \frac{1}{1-\delta}$ and $\exptail = - \frac{\log(1-\delta)}{\log K}$, we have
	  $$\mathbb P^\vertex_q (\jump^\tau \leq \pattern) < C \cdot \tau^{-\exptail}.$$
\end{proof}

\subsection{Criterion}
\label{criterion}

In this section we develop a graph theoretic criterion for the \unstable \ property.

\subsubsection{Counter example}
We give a simple example of a subgraph which prevents a simplicial system to be \unstable.
We will call such a subgraph \textit{stable}.
It will be a motivation for a criterion on subgraphs inducing the \unstability\ introduced in the next paragraph.\\

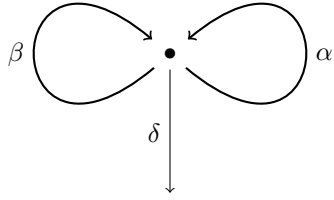
\begin{figure}[!htb]
	\center
	\begin{tikzpicture}[shorten >=1pt,node distance=2cm,on grid,auto]
	   \node (base) {$\bullet$};
	   \node (a) [below of =base] {};
	   \path[->] (base) edge[loop, thick, out=320, in =40, looseness=20]  node[right] {$\alpha$} (base);
	   \path[->] (base) edge[loop, thick, out=220, in =140, looseness=20]  node[left] {$\beta$} (base);
	   \path[->] (base) edge node[left] {$\delta$} (a);
        \end{tikzpicture}
	\caption{A stable subgraph.}
	\label{counter}
\end{figure}

Assume in the graph $\graph$ there exists a vertex $\vertex$ with three outgoing edges as in Figure~\ref{counter}.
Where all labels are distinct and the edge labeled $\delta$ points to any vertex in $\graph$.\\

To show stability, we will use the stopping time $\lose_\delta$ which corresponds to the first time $\delta$ loses.

The following lemma states that if the ratio between $q_\delta$ and $q_\alpha$ or $q_\beta$ is large enough then the probability that a path leaves the subgraph in finite time is small.
\begin{proposition}
	For all $q \in \cone$
	$$\mathbb P^\vertex_q ( \lose_\delta < \infty) \le \frac {\phi + \phi^{-1}} {1 - \phi^{-1}} \cdot \frac {q_\delta} {\min (q_\alpha, q_\beta)}$$
	where $\phi$ is the golden ratio.
\end{proposition}

\begin{proof}
	Let us assume that $q_\beta \le q_\alpha$.
	Notice that
	$$ \mathbb P^\vertex_q ( \delta \mid \beta \text{ or } \delta ) = \frac {q_\delta} {q_\delta + q_\beta} \le \frac {q_\delta} {q_\beta}.$$
	As $q_\delta$ is unchanged before it loses and $q_\beta$ is non decreasing, using the chain rule of Formula~(\ref{formula:cond}) we have
	\begin{align*}
		\mathbb P^\vertex_q ( \lose_\delta < \lose_\beta ) &= \sum_{n=0}^\infty \mathbb P^\vertex_q (\win_\alpha = n+1) \cdot \mathbb P^\vertex_{q \rauzymat_{\alpha}^n} ( \delta \mid \beta \text{ or } \delta )\\
								   &\le \mathbb P^\vertex_q (\win_\alpha < \infty) \cdot \frac {q_\delta} {q_\beta} .
	\end{align*}

	When $\beta$ loses, in the new distortion vector we have $q'_\beta \ge q_\alpha + q_\beta \ge q'_\alpha$.
	We induce this argument, switching $q'_\alpha$ and $q'_\beta$ to always have $q'_\beta \le q'_\alpha$.
	Thus we obtain a sequence of distortion vectors $q^{(n)}$ such that $q^{(0)} = q$, $q_\alpha ^{(n+1)} = q^{(n)}_\alpha + q^{(n)}_\beta$, $q_\beta ^{(n+1)} = q_\alpha^{(n)}$, $q^{(n)}_\delta = q_\delta$ and
	\begin{align*}
		\mathbb P^\vertex_q ( \lose_\delta < \infty) &\le \sum_{n=0}^\infty \frac {q^{(n)}_\delta} {q^{(n)}_\beta} \\
							     &\le \frac {q_\delta} {q_\beta} \cdot \sum_{n=1}^\infty F_n^{-1}
	\end{align*}
	where $F_n$ is the Fibonacci sequence.\\

	To compute this sum, notice that for $\phi$ the golden ratio $v_1 = (1, \phi)$ and $v_2 = (1, -\phi^{-1})$ are eigenvectors of eigenvalues $\phi$ and $-\phi^{-1}$ for the matrix associated to Fibonacci sequence.
	And $\phi v_1 + \phi^{-1} v_2 = (\phi + \phi^{-1}, \phi + \phi^{-1})$.
	Hence
	$$F_n = \frac {\phi^n - (-\phi^{-1})^n} {\phi + \phi^{-1} }$$
	and
	$$\sum_{n=1}^\infty F_n^{-1} \le (\phi + \phi^{-1}) \cdot \sum_{n=1}^\infty \phi^{-n} = \frac{\phi + \phi^{-1}}{1- \phi^{-1}}.$$

\end{proof}

This proposition implies that if a generic path visits $\vertex$ with $q_\delta$ at least three times as large as $q_\alpha$ and $q_\beta$ infinitely often then almost surely it will stay eventually in the subgraph of Figure~\ref{counter}.\\

This stability phenomenon is due to the fact that we have a subgraph with two edges of a subset $\largea$ of labels, assumed to be large, that play with each other.
Moreover the only edges leaving the subgraph are labeled outside of $\largea$ and play against edges labeled in $\largea$.
The distortion of the labels in this subset then increases exponentially fast, leaving few chances to lose for labels leaving the subgraph.\\

\subsubsection{First criterion}

We introduce a property on simplicial systems that prevents the phenomenon described above and will imply the \unstability.
This property is satisfied by a very large class of examples such as the Rauzy--Veech induction and most of multidimensional continued fractions algorithms as showed in Section~\ref{examples}.\\

The main idea here will be to consider degenerations of the induction where for some subsets of labels $\largea \subset \alphabet$ the distortion vector at these corresponding coordinates is infinitely larger than for others.
In particular, when we are on a vertex that has an outgoing edge labeled in $\largea$, any edges with a label outside of $\largea$ will almost surely not be chosen.\\

Let us denote by $\graph_\largea$ the subgraph of $\graph$, with the same set of vertices $\vertices$ and a set of edges defined as follows. For any $\vertex \in \vertices$,
\begin{itemize}
	\item if there is at least one edge in the outgoing edges $\vertexout$ labeled in $\largea$, the set of outgoing edges in $\graph_\largea$ is
	$$\vertexout^\largea = \{ e \in \vertexout \mid \labels(e) \in \largea \},$$
	\item otherwise,
	$$\vertexout^\largea = \vertexout.$$
\end{itemize}

\begin{definition}
	\label{rauzy}
	We say a simplicial system is strongly \rauzytype \ if
	\begin{enumerate}
		\item for all vertices every letter wins and loses in almost every path with respect to Lebesgue measure,
			\label{play}
		\item for all $\emptyset \subsetneq \largea \subsetneq \alphabet$ and all $\ccomponent$ strongly connected component of $\graph_\largea$ one of the following is true :
			\begin{enumerate}
				\item for all vertex $\vertex$ in $\ccomponent$ the cardinality $|\labels(\vertexout) \cap \largea| \le 1$,
					\label{card}
				\item for all vertex in $\ccomponent$ there is a path in $\graph$ labeled in $\largea$ leaving $\ccomponent$.
					\label{leave}
			\end{enumerate}
			\label{subset}
	\end{enumerate}
\end{definition}
The last property can be reformulated as: no letter in $\largea$ can win against another letter in $\largea$ in any strongly connected component of $G_\largea$ except if there is an edge labeled in $\largea$ leaving the component.
The first condition is a dynamical property, which makes it more difficult to check.
Using a result on subset of parameter in Section~\ref{fractals} we give after that section an equivalent definition that is purely graph theoretic.

\begin{proposition}
	The simplicial system associated to Rauzy--Veech induction on an irreducible interval exchange is \rauzytype \ and strongly connected.
\end{proposition}

\begin{proof}
	For Rauzy--Veech induction, each vertex has exactly two edges going in and two going out.
	Hence the connected component corresponding to the algorithm must be strongly connected.\\

	Property \ref{play} of \rauzytype \ comes from the observation that after a finite number of steps the subset of letters that never lose or win must form an invariant subinterval.
	This would contradict irreducibility (see \cite{Yoccoz10} for details).\\

	Assume that, as in the degenerate subgraph $\graph_\largea$, a subset of labels $\largea$ always loses against labels in its complementary set $\overline \largea$.
	If an interval labeled in $\overline \largea$ is at the right-hand side extremity of the interval exchange it can only lose to a letter in $\overline \largea$ thus there will always remain an interval labeled in $\overline \largea$ at the extremity after an arbitrary number of inductions.
	In such a configuration we cannot have two letters in $\largea$ playing with each other.

	If the two extremal intervals are labeled in $\largea$ then there is a path labeled in $\largea$ to an interval exchange with an extremal interval labeled in $\overline \largea$ since all labels in this complementary set win almost surely in finite time.
	This implies Property \ref{subset}.
\end{proof}

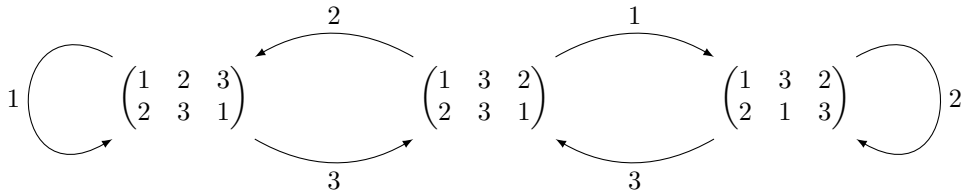
\begin{figure}[!h]
	\center
  \begin{tikzpicture}[>=stealth]
    \node (middle)
    { $ \begin{pmatrix}
          1 & 3 & 2\\
          2 & 3 & 1
        \end{pmatrix} $ };
    \node (left) [left=3cm]
    { $ \begin{pmatrix}
          1 & 2 & 3\\
          2 & 3 & 1
        \end{pmatrix} $ };
    \node (right) [right=3cm]
    { $ \begin{pmatrix}
          1 & 3 & 2\\
          2 & 1 & 3
        \end{pmatrix} $ };

    \draw[->,>=latex] (middle) to[bend left] node[above] {$1$} (right);
    \draw[->,>=latex] (right) to[bend left] node[below] {$3$} (middle);
    \draw[->,>=latex] (middle) to[bend right] node[above] {$2$} (left);
    \draw[->,>=latex] (left) to[bend right] node[below] {$3$} (middle);
    \draw[->,>=latex] (right) to [out=30, in=330, looseness=4] node[right] {$2$} (right);
    \draw[->,>=latex] (left) to [out=150, in=210, looseness=4] node[left] {$1$} (left);
  \end{tikzpicture}
  \caption{Rauzy graph for 3-IET.}
  \label{iet}
\end{figure}

	As an illustration, the reader can check directly these properties on the Rauzy graph for $3$-interval exchange transformations represented on Figure \ref{iet}.\\

\begin{remark*}
	The fully subtractive algorithm in dimension 3 or larger (see Section~\ref{examples}) provides a simple example of a simplicial system that is neither \rauzytype \ nor \unstable\ nor ergodic.
	The Poincaré algorithm in dimension 4 is a case that is not of \rauzytype \ but which is conjecturally ergodic.\\
\end{remark*}

Let $\alpha \in \alphabet$, $q \in \cone$ and $\tau > 1$.
We introduce some other useful properties on path for the following.\\

Let $\win_\alpha$ be the property of a path for which the letter $\alpha$ wins at its last step.
The set $\Gamma_\vertex(P_\alpha)$ can be thought of as the set of all paths stopping whenever $\alpha$ wins.\\

Similarly $\win^k_\alpha$ corresponds to whether the letter $\alpha$ wins at the last step of $\alpha$ or the length of the path is larger or equal to $k$.
It stops whenever $\alpha$ has won or the length of the path is $k$.\\

Let $\jump_{\alpha,q}^\tau$ be the property of a finite path $\finpath$ along which the distortion vector $q$ has \textit{jumped} by a factor $\tau$ \textit{i.e.} it satisfies
$$\left(q \rauzymat_\finpath\right)_\alpha \ge \tau \cdot q_\alpha.$$
We write $\jump_\alpha^\tau $ when $q$ corresponds to the distortion vector defining the measure.\\

\begin{lemma}
	\label{first Kerckhoff}
	In a simplicial system, for all $\alpha \in \alphabet$, $v \in V$, $\tau>1$ and $q \in \cone$,
	the probability that the vector $q$ jumps by $\tau$ on coordinate $\alpha$ before the letter $\alpha$ wins satisfies
	\[ \mathbb P^v_q\left( \jump_\alpha^ \tau < \win_\alpha < \infty \right) \le \frac 1 \tau.\]
\end{lemma}

\begin{proof}
	We prove by induction on $k$ that for all $k \in \mathbb N$
	\[ \mathbb P^v_q\left( \jump_\alpha^ \tau < \win_\alpha^ k \right) \le \frac 1 \tau.\]

	\noindent For $k=1$
	\[ \mathbb P^v_q\left( \jump_\alpha^ \tau < \win_\alpha^ 1 \right) = 0,\]
	since $\jump_\alpha^\tau$ is never satisfied by the empty path
	which keeps the vector $q$ unchanged. \\

	\noindent Assume now that the inequality is true for some $k$.
	As above, if $e \in \win_\alpha$ then $\win_\alpha^{k+1} = 1$ and $\mathbb P^v_q\left(\jump_\alpha^\tau < 1\right) = 0$.
	Hence
	\begin{align*}
		\mathbb P^v_q\left( \jump_\alpha^\tau < \win_\alpha^ {k+1} \right) &=
		\sum_{e \in \vertex_{\mathrm{out}}}
			\mathbb P^v_q\left(\jump_\alpha^\tau < \win_\alpha^ {k+1} \mid e \right)
			\cdot \mathbb P^v_q(e) \\
		&= \sum_{e \not\in \win_\alpha}
			\mathbb P^v_q\left(\jump_\alpha^\tau < \win_\alpha^ {k+1} \mid e \right)
			\cdot \mathbb P^v_q(e)
	\end{align*}
	If $e = (\vertex, \vertex') \not \in \win_\alpha$ and $q' := q \rauzymat_e$,
	observe that
	\begin{align*}
		e^{-1} \jump_{\alpha, q}^ \tau = \jump_{\alpha, q'}^{\tau'}, \\
		e^{-1} \win_\alpha^{k+1} = \win_\alpha^{k},
	\end{align*}
	where $\tau' = \tau \cdot \frac {q_\alpha} {q'_\alpha}$.\\

	Thus according to Proposition \ref{path_condition}
	\begin{align*}
		\mathbb P^v_q\left(\jump_{\alpha,q}^\tau < \win_\alpha^ {k+1} \mid e \right) = \mathbb P_{q'}^{\vertex'} \left(\jump_{\alpha,q'}^{\tau'} < \win^k_{\alpha} \right).
	\end{align*}

	If $\tau' \le 1$, $\mathbb P_{q'}^{\vertex'} \left(\jump_{\alpha,q'}^{\tau'} < \win^k_{\alpha} \right) \le 1 \le \frac {1} {\tau'}$.
	Otherwise, the recurrence hypothesis applied to the constant $\tau'$ implies
	\begin{align*}
		\mathbb P^v_q\left(\jump_{\alpha,q}^\tau < \win_\alpha^ {k+1} \mid e \right) \leq \frac 1 \tau \cdot \frac {q'_\alpha} {q_\alpha}.
	\end{align*}

	First assume that the label $\alpha$ appears in the vertices going out of $\vertex$, then there is only one edge that does not satisfies $\win_\alpha$ :
	it is the unique edge such that $\labels(e_\alpha)=\alpha$.
	Thus
	\[
		\mathbb P^v_q\left( \jump_\alpha^ \tau < \win_\alpha^{k+1} \right)
		< \frac 1 \tau \cdot \frac {q'_\alpha} {q_\alpha} \cdot \mathbb P^v_q(e_\alpha)
		= \frac 1 \tau.
	\]

	On the contrary, if the label $\alpha$ does not appear in the vertices leaving $\vertex$, we always have $q'_\alpha = q_\alpha$ and
	\[
		\mathbb P^v_q\left( \jump_\alpha^ \tau < \win_\alpha^{k+1} \right)
		<  \sum_{e \not \in \win_\alpha} \frac 1 \tau \cdot \mathbb P^v_q(e) \leq \frac 1 \tau.
	\]
	Hence
	\[
		\lim_{k \to \infty} \mathbb P^v_q\left( \jump_\alpha^ \tau < \win_\alpha^{k} \right) = \mathbb P^v_q \left( \jump_\alpha^ \tau < \win_\alpha < \infty \right) \le \frac {1} {\tau}.
	\]
\end{proof}

\begin{corollary}
	\label{rauzy:lower}
	Let $\vertex \in \vertices$, if a letter $\alpha$ wins in almost every path in $\Pi(\vertex)$, for all $\tau>1$ and $q \in \cone$
	\[ \mathbb P^\vertex_q\left( \win_\alpha \leq \jump_\alpha^ \tau \right) > 1 - \frac {1} {\tau}.\]
\end{corollary}
This lemma is the key ingredient to show one of our main theorems.
\begin{theorem}
	\label{rauzy-unstable}
	A strongly \rauzytype \ simplicial system is \unstable.
\end{theorem}
To show this theorem we will use an intermediate proposition on strongly connected components.
Let $\ccomponent$ be a strongly connected component of $\graph_\largea$.
The property $\stopping_\largea^\ccomponent$ is true if the path satisfies $\stopping_\largea$ or if it goes through an edge in $\graph_\largea$ outside of $\ccomponent$.
In other words, the stopping time corresponds to the state when the value of the distortion on a coordinate in $\alphabet \setminus \largea$ reaches the level of the initial distortion on $\largea$ or leaves the strongly connected component $\ccomponent$.

\begin{proposition}
	\label{induction-rauzy}
	Let $\vertex \in \ccomponent$, for all $\tau>1$, $K>1$ and $q$ $(\largea, K)$-balanced distortion vector,  if $\ccomponent$ satisfies, in Definition~\ref{rauzy}, property \ref{card} then
	\[
		\mathbb P^\vertex_q( \stopping_\largea^\ccomponent \le \jump^{\tau} ) > \frac {1} {|\alphabet| \cdot \tau K} \cdot \left(1 - \frac {1} {\tau}\right).
	\]
	If it satisfies \ref{leave} then
	\[
		\mathbb P^\vertex_q( \stopping_\largea^\ccomponent \le \jump^{\tau} ) > \left(\frac {1} {2^{|\vertices|} \cdot K}\right)^{|\vertices|} .
	\]
\end{proposition}
\begin{proof}
	We start with the case when $\ccomponent$ satisfies \ref{card}.
	Corollary~\ref{rauzy:lower} implies that for all vertex $\vertex$ and $\tau > 1$
	$$\mathbb P^\vertex_q \left(\min_{\alpha \in \largea} \win_\alpha \le \min_{\alpha \in \largea} \jump_\alpha^\tau \right) > 1 - \frac {1} {\tau}.$$
	Let us introduce two sets
	$$A := \left\{ \infpath \in \Pi^\infty(\vertex) \mid \min_{\alpha \in \largea} \win_\alpha(\infpath) \le \min_{\alpha \in \largea} \jump_\alpha^\tau(\infpath) \right\}$$
	and
	$$S := \left\{ \infpath \in \Pi^\infty(\vertex) \mid \stopping_\largea^\ccomponent(\infpath) \le \jump^{\tau}(\infpath) \right\}.$$
	To prove the proposition, we will show the inequality
	$$\mathbb P^\vertex_q (S) \ge \frac {1} {|\alphabet| \cdot \tau K} \cdot \mathbb P^\vertex_q (A).$$\\

	We first separate cases when the jump happens on a label in or out of $\largea$ \textit{i.e.} we show inequalities for the intersection with the set
	$$B := \left\{ \infpath \in \Pi^\infty(\vertex) \mid \max_\alphabet q \rauzymat_{\infpath_n} = \max_{\largea} q \rauzymat_{\infpath_n} \text{ where } n = T_{\jump^{\tau}}(\infpath) < \infty \right\}$$
	and its complementary $\overline B$.
	Notice that when $\infpath \in B$, $\min_{\alpha \in \largea} \jump_\alpha^\tau \le \jump^\tau.$\\

	Let $\infpath \in \overline B$, as $n = \jump^\tau(\infpath) < \infty$ almost surely by property \ref{play} of the strongly \rauzytype \ simplicial systems we can assume that
	$$\max_\alphabet q \rauzymat_{\infpath_n} = \max_{\alphabet \setminus \largea} q \rauzymat_{\infpath_n},$$
	thus there exists $\beta$ in $\alphabet \setminus \largea$ such that $(q \rauzymat_{\infpath_n})_\beta \ge \tau \max q$ and for every other letter $\alpha$ in $\alphabet$, $(q \rauzymat_{\infpath_n})_\alpha < \tau \max q$.
	Hence $\max_{\alphabet \setminus \largea} q \rauzymat_{\infpath_n} \ge \min_\largea q \rauzymat_{\infpath_n}$ and $\stopping_\largea(\infpath) \le \jump^{\tau}(\infpath)$ in other words $\infpath \in S$.
	This implies $\overline B \subseteq S$.\\

	Another important separation will be made with the set $C$ of paths $\infpath$ in $\Pi^\infty(\vertex)$ such that $m = \min_{\alpha \in \largea} \win_\alpha(\infpath) - 1 < \infty$ and the last vertex of the path before a letter in $\largea$ wins $\vertex' := \infpath_{m} \cdot \vertex$ has at least two edges labeled in $\largea$.
	In notations
	$$C := \left\{ \infpath \in \Pi^\infty (\vertex) \mid m = \min_{\alpha \in \largea} \win_\alpha(\infpath)-1 < \infty \text{ and } |l\left((\infpath_m \cdot \vertex)_{\mathrm{out}}\right) \cap \largea| \ge 2 \right\}.$$

	By inclusion $\overline B \subseteq S$ we have
	$$
		\mathbb P^\vertex_q(A \cap \overline B \cap \overline C) \le \mathbb P^\vertex_q(S \cap \overline B \cap \overline C).
	$$
	Notice that for $\infpath \in \overline C$, the first time a letter in $\largea$ wins it wins against a letter outside of $\largea$.
	Assume moreover that $\infpath \in \overline{B} \cap \overline{C}$ then it satisfies $\stopping_\largea(\infpath)$ before $\min_{\alpha \in \largea} \win_\alpha(\infpath) \le \min_{\alpha \in \largea} \jump^\tau_\alpha(\infpath) \le \jump^\tau(\infpath)$.
	Hence
	$$
		\mathbb P^\vertex_q(A \cap B \cap \overline C) \le \mathbb P^\vertex_q(S \cap B \cap \overline C).
	$$
	These two inequalities imply
	\begin{equation}
		\label{Ccomp}
		\mathbb P^\vertex_q(A \cap \overline C) \le \mathbb P^\vertex_q(S \cap \overline C).
	\end{equation}
	\\

	Let $\infpath \in A \cap C$, notice that for all edges $e \in (\infpath_m \cdot \vertex)_\mathrm{out}$ and all $\infpath_e \in \Pi^\infty (\infpath_m \cdot e \cdot \vertex)$ we have
	$$\infpath_m \cdot e \cdot \infpath_e \in A \cap C$$
	since by definition of $C$ for any edge going out of $\infpath_m \cdot \vertex$ at least one letter in $\largea$ wins.
	Let $\Gamma$ be the set of finite paths $\infpath_m$ for all $\infpath \in A \cap C$.
	Then
	\begin{equation}
		\label{sumgamma}
		\mathbb P_q^\vertex (A \cap C) = \sum_{\finpath \in \Gamma} \mathbb P_q^\vertex (\finpath).
	\end{equation}

	As $q$ is by assumption $(\largea, K)$-balanced, the property of $A$ implies that the distortion at $\infpath_m \cdot \vertex$ satisfies $\max_\largea q \rauzymat_{\infpath_m} \le \tau K \min_\largea q$.
	If one of the label outside of $\largea$ has the largest distortion coefficient then its probability is bounded below by $\frac {1} {|\alphabet|}$.
	Otherwise the previous remark implies that $\rauzymat_{\infpath_m} q$ is $(\largea, \tau K)$-balanced.
	Thus there exists an edge $e_{\infpath_m} \in (\infpath_m \cdot \vertex)_\mathrm{out}$ labeled outside of $\largea$ or going out of $\ccomponent$ such that
	\begin{equation}
		\label{Pe}
		\mathbb P_q^\vertex (\infpath_m \cdot e_{\infpath_m} \mid \infpath_m) \ge \frac {1} {|\alphabet| \cdot \tau K}.
	\end{equation}

	Let us show that for all $\infpath_e' \in \Pi^\infty (\infpath_m \cdot e_{\infpath_m} \cdot \vertex)$ the induced infinite path $\infpath^* := \infpath_m \cdot e_{\infpath_m} \cdot \infpath_e'$  belongs to $S \cap C$.
	If $\infpath^* \in \overline B$ we have already seen above that $\infpath^* \in S$.
	If $\infpath^* \in B$ then by the choice of $e_{\infpath_m}$ and property of B $\infpath^* \in S$.
	Hence, using this fact, (\ref{sumgamma}) and (\ref{Pe}),
	\begin{align}
	\nonumber
	\mathbb P_q^\vertex (S \cap C) & \ge \sum_{\finpath \in \Gamma} \mathbb P_q^\vertex (\finpath \cdot e_{\finpath}) \ge \frac {1} {|\alphabet| \cdot \tau K} \sum_{\finpath \in \Gamma} \mathbb P_q^\vertex (\finpath_m) \\
				       & \ge \frac {1} {|\alphabet| \cdot \tau K} \cdot \mathbb P_q^\vertex (A \cap C).
				       \label{AC}
	\end{align}

	Formula (\ref{Ccomp}) and (\ref{AC}) finally imply
	\[
		\mathbb P_q^\vertex (S) \ge \frac {1} {|\alphabet| \cdot \tau K} \cdot \mathbb P_q^\vertex (A).
	\]

	In the case $\ccomponent$ satisfies \ref{leave}, there is a finite path $\finpath$ from $\vertex$ labeled in $\largea$ leaving $\ccomponent$.
	Up to removing loops one can assume that $\largea$ as at most $|\ccomponent|$ steps.
	Where $|\ccomponent|$ is the number of vertices in $\ccomponent$.
	As $q$ is assumed to be $(\largea, K)$-balanced it remains at least $(\largea, 2^{|\ccomponent|} K)$-balanced at each step of the path.
	By Proposition~\ref{balanced-bound},
	\[
		\mathbb P^\vertex_q( \stopping_\largea^\ccomponent \le \jump^{\tau} ) > \left(\frac {1} {2^{|\ccomponent|} \cdot K}\right)^{|\ccomponent|} .
	\]

\end{proof}

Recall that the strongly connected components of a directed graph induce an acyclic graph called condensation of a directed graph (see for instance section~3.4 of \cite{BondyMurty08}).
In the condensation graph there exists minimal vertices \textit{i.e.} vertices from which there are no outgoing edges.
The corresponding strongly connected components are called minimal components.
We end the proof of the theorem by inducing on the distance, called height, of a strongly connected component to minimal ones in the condensation graph.
\begin{proof}[Theorem~\ref{rauzy-unstable}]
	Let $\emptyset \subsetneq \largea \subsetneq \alphabet$ and $K>1$.
	Assume that the vertex $\vertex$ is in a minimal strongly connected component $\ccomponent$.
	In this case there are no edges in $\graph_\largea$ going out of the strongly connected component thus $\stopping_\largea^\ccomponent = \stopping_\largea$ and Lemma~\ref{induction-rauzy} implies for every $\tau > 1$ and $\delta_0 := \frac {1} {|\alphabet|\cdot \tau K} \cdot \left(1 - \frac {1} {\tau}\right)$ we have the lower bound of the \unstability\
	\[
		\mathbb P_q^\vertex \left(\stopping_\largea \le \jump^\tau\right) > \delta_0.
	\]

	We prove by induction that if $\vertex$ is in a strongly connected component of height $h$ in the tree formed by the condensation graph then
	\[
		\mathbb P_q^\vertex \left(\stopping_\largea \le \jump^{\tau^h}\right) > \delta_0^h.
	\]
	Let $\vertex$ be in a component of height $h$ and let $\infpath$ satisfying $\stopping_\largea (\infpath) \le \jump^\tau (\infpath) =: n$.
	Recall that by property $\ref{rauzy}$ of strongly \rauzytype \ simplicial systems $n<\infty$ almost surely.
	For $\Gamma$ the set of such finite paths $\infpath_n$, by Proposition~\ref{induction-rauzy},
	$$\sum_{\finpath \in \Gamma} \mathbb P_q^v (\finpath) > \delta_0.$$
	In the case $\stopping_\largea = \stopping_\largea^\ccomponent$ the lower bound is satisfied.
	Let us assume $\stopping_\largea > \stopping_\largea^\ccomponent$ then the path in time $\stopping_\largea^\ccomponent$ goes through an edge in $\graph_\largea$ that leaves the component $\ccomponent$ to go to another strongly connected component $\ccomponent'$ of height $h-1$.
	Notice that in this case, for $q' := q \rauzymat_{\infpath_n}$, $\infpath_n^{-1} \stopping_{\largea, q} = \stopping_{\largea, q'}$ and $\infpath_n^{-1} \jump^{\tau^h}_q \ge \jump^{\tau^{h-1}}_{q'}$ thus
	\[
		\mathbb P_q^\vertex \left(\stopping_\largea \le \jump^{\tau^h} \mid \infpath_n\right) \ge \mathbb P_q^{\vertex'} \left(\stopping_\largea \le \jump^{\tau^{h-1}}\right) > \delta_0^{h-1}
	\]
	where $\vertex' := \finpath_n \cdot \vertex \in \ccomponent'$.
	Hence using the chain rule of Formula~(\ref{formula:cond})
	\begin{align*}
		\mathbb P_q^\vertex \left(\stopping_\largea \le \jump^{\tau^h}\right) =
		\sum_{\finpath \in \Gamma} \mathbb P_q^\vertex \left(\gamma\right) \cdot \mathbb P_q^\vertex \left(\stopping_\largea \le \jump^{\tau^h} \mid \finpath \right)
		& \ge \sum_{\finpath \in \Gamma} \mathbb P_q^{\vertex} (\finpath) \cdot \delta_0^{h-1} \\
		& > \delta_0^h.
	\end{align*}
\end{proof}

\subsubsection{Graph criterion}

Theorem~\ref{fractal:measure} enables us to relax Property~\ref{play} on winning and losing letters in Definition~\ref{rauzy} of strongly \rauzytype \ simplicial systems to a purely graph theoretic one.
\begin{customdef}{1'}
	\label{rauzyprime}
	We say a simplicial system is \rauzytype \ if
	\begin{enumerate}
		\item from every vertex there exists a path along which each label in $\alphabet$ appears,
			\label{everyletter}
		\item for all $\emptyset \subsetneq \largea \subsetneq \alphabet$ and all $\ccomponent$ strongly connected component of $\graph_\largea$ one of the following is true :
			\begin{enumerate}
				\item for all vertex $\vertex$ in $\ccomponent$ the cardinality $|\labels(\vertexout) \cap \largea| \le 1$,
				\item from every vertex in $\ccomponent$ there is a path in $\graph$ labeled in $\largea$ leaving $\ccomponent$.
			\end{enumerate}
	\end{enumerate}
\end{customdef}

This section will be dedicated to proving the following corollary of Theorem~\ref{fractal:measure}.
\begin{theorem}
	\label{strongweak}
	If a simplicial system is \rauzytype \ then it is strongly \rauzytype.
\end{theorem}

A first step in the proof if the following lemma.
It enables us to relax Property~\ref{play} that each label wins and loses almost surely to just each label plays almost surely.
\begin{lemma}
	\label{losewin}
	If in a simplicial system $\graph$ there is some label $\alpha \in \alphabet$ such that, for all vertex and all distortion vector, $\alpha$ plays almost surely then for all vertex the label $\alpha$ wins and loses infinitely many times almost surely.
\end{lemma}

\begin{proof}
If a label plays almost surely for all vertex and all distortion vector then it plays almost surely infinitely many times.
This is shown recursively by conditioning by finite paths at the time when the given label plays.
Hence if a label has the property of the lemma it either wins of loses infinitely many times almost surely.\\

If $\alpha$ wins infinitely many times then by Lemma~\ref{first Kerckhoff}, for all $\tau > 1$, the distortion vector at coordinate $\alpha$ jumps almost surely by $\tau$.
In particular the distortion at $\alpha$ goes almost surely to infinity and thus $\alpha$ loses infinitely many times.\\

Let us now assume that $\alpha$ loses infinitely many times.
For all $n \ge 1$, we introduce a stopping time $\lose_\alpha^n$ which corresponds to the $n$-th time $\alpha$ loses.
The lemma is a consequence of the following proposition.
\begin{proposition}
	For $\graph$ and $\alpha$ as above, for all $n \ge 1$ and $q \in \cone$
	$$\mathbb P^\vertex_q \left(\lose_\alpha^{(n)} < \win_\alpha\right) \le \frac {q_\alpha} {q_\alpha + n q_\delta}$$
	where $q_\delta = \min_{\beta \in \alphabet} q_\beta$.
\end{proposition}

\begin{proof}
	For all vertex $\vertex$ let $\Gamma(\vertex)$ be the set of finite paths starting at $\vertex$ and ending at a vertex with two outgoing edges one of which is labeled by $\alpha$ and such that $\alpha$ neither wins or loses in the path.
	Let $\Gamma^n := \left\{ (\finpath_1, \dots, \finpath_n) \mid \finpath_1 \in \Gamma(\vertex), \ \finpath_2 \in \Gamma(\finpath_1 \cdot \alpha \cdot \vertex), \ \dots, \ \finpath_n \in \Gamma(\finpath_1 \cdot \alpha \cdot \finpath_2 \cdot \alpha \dots \finpath_{n-1} \cdot \alpha \cdot \vertex) \right\}$.
	For all $(\finpath_1, \dots, \finpath_n) \in \Gamma^n$ and all $1 \le i \le n$ we introduce the notation
	$$v^{i} = \finpath_1 \cdot \alpha \dots \finpath_i \cdot \vertex$$
	and
	$$q^{i} = q \rauzymat_{\finpath_1 \cdot \alpha \dots \finpath_i}.$$

	We can decompose
	\begin{align*}
		\mathbb P^\vertex_q \left(\lose_\alpha^{(n)} < \win_\alpha \right) &= \sum_{(\finpath_1, \dots, \finpath_n) \in \Gamma^n} \mathbb P^\vertex_q \left(\finpath_1 \cdot \alpha \dots \finpath_n \cdot \alpha\right)\\
										   &= \sum_{(\finpath_1, \dots, \finpath_n) \in \Gamma^n} \mathbb P^\vertex_q \left(\finpath_1\right) \cdot \mathbb P^{\vertex^1}_{q^1} (\alpha) \dots \mathbb P^{\alpha \cdot \vertex^{n-1}}_{q^{n-1} \rauzymat_\alpha} (\finpath_n) \cdot \mathbb P^{\vertex^{n}}_{q^{n}} (\alpha) \\
\end{align*}
For a given element of $\Gamma^n$
$$ \mathbb P^{\vertex^1}_{q^1} (\alpha) \dots \mathbb P^{\vertex^{n}}_{q^{n}} (\alpha) = \frac {q^{1}_{\alpha}} {(q^{1} \rauzymat_\alpha)_\alpha} \dots \frac {q^{n}_{\alpha}} {(q^{n} \rauzymat_\alpha)_\alpha}.$$
As $\alpha$ does not play in paths $\finpath_i$ the distortion vector is unchanged on coordinate $\alpha$ thus $q^{i+1}_\alpha = \left(q^i \rauzymat_\alpha\right)_\alpha$ for all $i \ge 1$ and $q^1_\alpha = q_\alpha$.
Moreover each time $\alpha$ loses the distortion vector increases by at least $q_\delta$ on coordinate $\alpha$.
Hence
$$ \mathbb P^{\vertex^1}_{q^1} (\alpha) \dots \mathbb P^{\vertex^{n}}_{q^{n}} (\alpha) \le \frac {q_{\alpha}} {q_\alpha + n q_\delta}$$
which implies
$$\mathbb P^\vertex_q \left(\lose_\alpha^{(n)} < \win_\alpha \right) \le \frac {q_{\alpha}} {q_\alpha + n q_\delta}.$$
\end{proof}
\end{proof}

The theorem is now a consequence of the following lemma.
\begin{lemma}
	\label{zeromes}
	Let $\graph$ be a \rauzytype \ simplicial system.
	Assume that for some $\vertex \in \vertices$ all labels in $\loser \subset \alphabet$ win and lose almost surely infinitely many times in paths of $\Pi^\infty(\vertex)$.
	Let $\subgraph$ be a \ subgraph of $\graph$ containing only vertices of $\graph$ such that its outgoing edges are all labeled in $\loser$ and there exists a path of $\graph$ in $\Pi^\infty(\vertex)$ labeled in $\loser$ that leaves $\subgraph$.
	Then the subset of parameters in the simplex associated to $\Pi^\infty(\vertex)$ which corresponds to paths remaining in $\subgraph$ has zero Lebesgue measure.
\end{lemma}

\begin{proof}
	Let $\stopping^\subgraph$ be the time when a path visits an edge out of $\subgraph$.
	By assumption every letter $\alpha$ in $\loser$ wins in almost every path thus by Lemma~\ref{first Kerckhoff} for all $\tau > 1$ and $q \in \cone$
	$$\mathbb P^\vertex_q (\win_\alpha \le \jump_\alpha^\tau) > 1 - \frac {1} {\tau}.$$
	Now we can reproduce the proof of Proposition~\ref{rauzy-unstable} by taking for all $\largea \subsetneq \loser$ the stopping time $\min(\stopping^\ccomponent_\largea, \stopping^\subgraph)$ and $\jump^\tau_{\loser}$ instead of $\stopping^\ccomponent_\largea$ and $\jump^\tau$.
	We then get a modified \unstability \textit{i.e.} for all non-empty subset $\largea \subsetneq \loser$ and all $K>1$ there exists $\tau > 1$ and $\delta > 0$ such that for all vertex $\vertex \in \vertices$ and all $(\largea, K)$-balanced distortion vector $q \in \cone$
	$$\mathbb P^\vertex_q \left( \min( \stopping_\largea, \stopping^\subgraph) \le \jump^\tau_\loser\right) > \delta.$$
	Lemma~\ref{minmax} can also be reproduced to get some $\tau > 1$ and $\delta > 0$ such that for all vertices $\vertex \in \vertices$ and all distortion vector $q \in \cone$
	$$\mathbb P^\vertex_q \left( \min( \minmax_{\loser}, \stopping^\subgraph) \le \jump^\tau_\loser\right) > \delta.$$
	By assumption there is a path in $\graph$ leaving $\subgraph$ such that only labels in $\loser$ play.
	Notice that for some $\infpath \in \Pi^\infty(\vertex)$ such that $n = \minmax_\loser(\infpath) \le \jump^\tau_{\loser}(\infpath)$ then $q' = q \rauzymat_{\infpath_n}$ is $(\loser,\tau)$-balanced and by Proposition~\ref{balanced-bound}
	$$\mathbb P^{\infpath_n \cdot \vertex}_{q'} \left( \stopping^\subgraph \le \jump^\tau_\loser \right) > \left(\frac {1} {2^{|\vertices|} \cdot \tau }\right)^{|\vertices|}.$$
	Hence there exists $\delta' > 0$ such that for all $q \in \cone$
	$$\mathbb P^\vertex_q \left( \stopping^\subgraph \le \jump^{\tau^2}_\loser \right) > \delta \cdot \delta'$$
	which implies that almost surely a path will go out of $\subgraph$.
\end{proof}

\begin{proof}[Theorem~\ref{strongweak}]
	Let $\graph$ be a \rauzytype \ simplicial system.
	There exists $\vertex$ a vertex of $\graph$ such that paths starting at $\vertex$ come back in finite time with positive probability.
	Let $\loser$ be the maximal subset of labels in $\alphabet$ such that for a path in $\Pi^\infty(\vertex)$ a label $\alpha$ in $\loser$ plays infinitely often almost surely.
	By definition, a vertex which has at least one of its outgoing edges labeled in $\overline {\loser}$ is visited almost surely finitely many times.
	Thus up to taking several iteration of the first return map to $\vertex$, there exists a subset of paths in $\Pi^\infty(\vertex)$ of positive Lebesgue measure which do not visit these latter vertices.
	Let $\subgraph$ be the subgraph of $\graph$ to which we remove all the edges labeled in $\overline {\loser}$.
	By Property~\ref{everyletter} of Definition~\ref{rauzyprime} for \rauzytype \ simplicial systems there exists a path from $\vertex$ in which all label in $\alphabet$ appear.
	This path leaves $\subgraph$ as soon as it goes through an edge not labeled in $\loser$.
	But the set of paths remaining in $\subgraph$ has positive Lebesgue measure which contradicts Lemma~\ref{zeromes}.
\end{proof}

\section{Ergodic measure and fractal dimension}

Building on the results of the previous section we show ergodicity results on quickly escaping simplicial systems.
Our main tool will be thermodynamical formalism applied to an induced map of the win-lose induction on a sub-simplex of the parameter space.
The previous section provides us with integrability conditions to use this formalism.

In the last part of this section we study subsets of parameters for which the win-lose induction remains in a subgraph of a simplicial system.
These sets generalize some classical fractal sets such as the set of real numbers with bounded continued fractions expansion or the Rauzy gasket.
As a consequence of all the work done before we can give an upper bound to the Hausdorff dimension of such fractal sets.
This improves in particular the previous bound obtained by Avila--Hubert--Skripchenko \cite{AvilaHubertSkripchenko16b} for the Rauzy gasket and gives new bounds for its generalizations to higher dimensions.

\subsection{A uniformly expanding acceleration}
\label{acceleration}
The win-lose induction associated to a simplicial system is not uniformly hyperbolic.
As for Markov shifts on a finite number of states, a useful idea is to wait until a given path $\ppath$ appears in the coding.
Remark~\ref{recursion} states that such an acceleration of the win-lose induction can be defined for \unstable\ simplicial systems.
For a good choice of path, this acceleration will be uniformly hyperbolic.
If we assume that the chosen path starts and ends at the same vertex, this is exactly considering the first return map on the given vertex to the subsimplex $\simplex_\ppath := \rauzymat_\ppath \simplex$, which we denote by
$$T_\ppath : \simplex_\ppath \to \simplex_\ppath.$$

\begin{proposition}
	\label{uniform}
	If $\rauzymat_{\ppath}$ is a positive matrix, $\rauzymap_\ppath$ is uniformly expanding. We say in this case that $\ppath$ is a positive path.
\end{proposition}

Before giving a proof of this proposition, we need to introduce some definitions.
For any two vectors $v, w \in \cone$, let
\[
	\alpha (v,w) := \max_{a \in \alphabet} \frac {v_a} {w_a} \text{,} \hspace{1cm} \beta (v,w) := \min_{a \in \alphabet} \frac {v_a} {w_a}
\]
and
\[
	d(v,w) := \log \frac {\alpha(v,w)} {\beta(v,w)}.
\]
One can check that $d$ is a complete metric on the projectivization of $\cone$ called the \textit{Hilbert metric}.
This metric has the very useful feature that any linear map induced by a positive matrix is contracting with respect to it.
\begin{proposition}
	\label{Hilbert}
	For any non-negative matrix $M$, we have
	\[
	d(Mv,Mw) \leq d(v,w),
	\]
	moreover if $M$ is positive, there exists $\theta < 1$ such that
	\[
	d(Mv,Mw) \leq \theta d(v,w).
	\]
\end{proposition}

\begin{proof}
	This is a well known property of Hilbert metrics, the proof can be found \textit{e.g.} in \cite{Viana97}.
\end{proof}

\begin{proof}[Proposition~\ref{uniform}]
	Let $\ppath \cdot \finpath_e$ be the Rauzy path for a given point in $\simplex_\ppath$ until its first return.
	The inverse of the Rauzy map is a projectivization of the linear map $\rauzymat_\ppath \rauzymat_{\finpath_e}$, which is, according to Lemma~\ref{Hilbert}, the composition of a weakly contracting map and a contraction a contracting map with coefficient $\theta < 1$ for the Hilbert metric on $\simplex$.
	Hence the inverse of the Rauzy map is contracting for the coefficient $\theta$ depending only on $\ppath$.
	Moreover, by positivity, $\simplex_\ppath$ is precompact in $\simplex$, thus the Hilbert metric is equivalent to all finite metric on this space.
\end{proof}

\def \cylinder {\underline w}
\begin{remark*}
	We will use the notations $\psimplex$ for $\simplex_\ppath$, $\rauzyacc$ for $T_\ppath$.
\end{remark*}

In the case of \unstable\ simplicial system on can always find such an acceleration.
\begin{proposition}
	\label{positive}
	If $\graph$ is a \unstable\ strongly connected simplicial system then it admits a positive path.
\end{proposition}

\begin{proof}
	Notice that the \unstability\ implies in particular that $\graph_\largea$ is a strict subgraph of $\graph$.
	Thus for any subset of letter $\largea$ there is a vertex such that there is an outgoing edge labeled in $\largea$ and an other labeled in the complementary set.
	We can then construct a positive path by recurrence.
\end{proof}

\begin{corollary}
	If a simplicial system is \unstable\ and strongly connected then its win-lose induction is ergodic with respect to the any invariant measure absolutely continuous with respect to Lebesgue measure.
\end{corollary}
We will show the existence of such a measure in Corollary~\ref{invariant_measure}.

\begin{proof}
	Let $\ppath$ be a positive path as constructed in Proposition~\ref{positive}, the acceleration $\rauzyacc$ is uniformly expanding thus any invariant subset has Lebesgue measure either $0$ or $1$.
	Thus $\rauzyacc$ and by extension $T$ are ergodic with respect respect to the invariant measure absolutely continuous with respect to Lebesgue measure.
\end{proof}

We will assume in the following that we are given a positive path $\ppath$ that starts and ends at the same vertex of the graph.

\subsection{Thermodynamic formalism}
\label{thermo}

Considering a well chosen acceleration associated to a positive path $\ppath$ that starts at some vertex $\vertex$ we have reduced the dynamical study of a quickly escaping win-lose induction to the study of a map on a simplex $\rauzyacc : \psimplex \to \psimplex$.

The domains of definition for this map tile $\psimplex$ by subsimplices and are labeled by $S$ the set of loops in $\graph$ starting at $v$ and that do not contain the path $\ppath$.
We denote by $\psimplex_{w}$ the domain corresponding to $w \in \states$ in $\psimplex$ which stands for the subsimplex for which the coding of $\rauzyacc$ starts by $w$.\\

The map $\rauzyacc$ is then conjugated to the full shift on $\Sigma = \states^{\N}$.
To every cylinder of the shift $\underline w = [w_1, \dots, w_n]$ we can associate the subsimplex $\psimplex_{\underline w}$ of points with the corresponding coding.
We will make the abuse of also calling this subsimplex a cylinder.\\

It will turn out fruitful to study the canonical suspension of the win-lose induction, introduced in Section~\ref{suspension}, by considering its suspension time as a potential in the thermodynamical study of this full shift.

With this approach we will be able to show that there exists a unique invariant measure of maximal entropy for the canonical suspension of the win-lose induction.
The general strategy consists in showing existence of Gibbs measures for a family of potentials parametrized in $\R$.
To these potentials is associated a pressure for which, when zero, the associated Gibbs measure induces a measure of maximal entropy on the suspension.

\subsubsection{Properties on the norm}
We list some easy but nonetheless useful properties for the $L^1$ norm on the simplex.
\begin{proposition}
	\label{proposition:norm}
	Let $\finpath$ a finite path in $\graph$, $\lambda \in \simplex$ and $\lambda' = \rauzymat_\finpath^{-1} \lambda$.\\

	If $\normone \lambda = 1$ then
	\[
		\lambda = \frac {\rauzymat_\finpath \lambda'} {\normone {\rauzymat_\finpath \lambda'}}.
	\]

	If $\normone {\lambda'} = 1$ then
	\[
		\normone{\rauzymat_\finpath^{-1} \lambda} = \frac {1} {\normone {\rauzymat_\finpath \lambda'}}.
	\]

	Moreover if we can decompose $\finpath = \finpath_1 \cdot \finpath_2$ and if $\lambda_1 = \frac {\rauzymat_{\finpath_1}^{-1} \lambda} {\normone{\rauzymat_{\finpath_1}^{-1} \lambda}} $,
	\[
		\normone{\rauzymat_\finpath^{-1} \lambda} = \normone{\rauzymat_{\finpath_2}^{-1} \lambda_1} \cdot \normone {\rauzymat_{\finpath_1}^{-1} \lambda}.
	\]
\end{proposition}

\begin{proposition}
	\label{L1}
	Let $v, w \in \cone$,
	\[
		\frac {\normone{v}} {\normone{w}} \leq \max_{\alpha \in \alphabet} \frac {v_\alpha} {w_\alpha}.
	\]
\end{proposition}

\subsubsection{Roof function}
The canonical suspension flow on $\parameter$ can be defined on the base $\psimplex$ with an \textit{accelerated} roof function defined, for $x \in \psimplex$, as
\[
	\roofacc(x) = \roof(x) + \roof(Tx) + \dots + \roof(T^{n-1} x) = - \log \normone {\rauzymat_\finpath^{-1} x}
\]
where $n \ge 1$ is the smallest integer such that $T^n x \in \psimplex$ and $\finpath$ is the finite path in the graph which is the coding of $x$ until it returns to $\psimplex$.
The second equality uses Proposition \ref{proposition:norm}.
\begin{remark*}
	This potential is similar to the geometric potential in the context of thermodynamic formalism.
\end{remark*}

Let $0<\theta<0$ be the constant associated to the matrix $\rauzymat_{\ppath}$ by Lemma~\ref{Hilbert}.
We show that the accelerated roof function is Hölder of order $\beta := \log(1/\theta)$.
\begin{proposition}
	\label{holder}
	For all $x, y \in \psimplex$ in the same $n$-cylinder $\psimplex_{\cylinder}$, where $n \ge 1$,
	\[
		| \roofacc(x) - \roofacc(y) | \leq \theta^{n+1} \cdot \diam(\psimplex).
	\]
\end{proposition}

\begin{proof}
	Let $x, y \in \psimplex$ be in the same cylinder $\psimplex_{\cylinder}$ which corresponds to the path
	$$\finpath = \ppath \cdot w_1 \cdot \ppath \cdot w_2 \dots \ppath \cdot w_n.$$
	Then according to Proposition~\ref{proposition:norm},
	\[
		| \roofacc(x) - \roofacc(y) | = \left| \log \frac {\normone{\rauzymat_{\ppath w_1}^{-1} y}} {\normone{\rauzymat_{\ppath w_1}^{-1} x}} \right|
		= \left| \log \frac {\normone{\rauzymat_{\ppath w_1} x}} {\normone{\rauzymat_{\ppath w_1} y}} \right|.
	\]
	By Proposition~\ref{L1}, up to switching $x$ and $y$, we can bound the distance by the Hilbert metric.
	\begin{equation}
		|\roofacc(x) - \roofacc(y)| \leq d(\rauzymat_{\ppath w_1} x, \rauzymat_{\ppath w_1} y).
		\label{roof inequality}
	\end{equation}
	Let $x',y' \in \psimplex$ be such that $x = \rauzymat_\finpath x'$ and $y = \rauzymat_\finpath y'$.
	With the notation $\finpath = \ppath \cdot w_1 \cdot \finpath'$, using Lemma~\ref{Hilbert},
	\begin{align*}
		d(x, y) & = d(\rauzymat_{\ppath} \rauzymat_{w_1} \rauzymat_{\finpath'} x', \rauzymat_{\ppath} \rauzymat_{w_1} \rauzymat_{\finpath'} y') \\
			& \leq \theta \cdot d(\rauzymat_{w_1} \rauzymat_{\finpath'} x', \rauzymat_{w_1} \rauzymat_{\finpath'} y') \\
			& \leq \theta \cdot d(\rauzymat_{\finpath'}x', \rauzymat_{\finpath'}y').
	\end{align*}
	By induction on $n$, we obtain
	$$d(\rauzymat_{\ppath w_1} x, \rauzymat_{\ppath w_1} y) \leq \theta \cdot d(x,y) \le \theta^{n+1} \cdot \diam(\psimplex).$$
\end{proof}
\begin{remark}
	\label{holder2}
	We only have used in the proof of the previous proposition the fact that matrices $M_\finpath$ are non-negative and $M_\ppath$ is positive.
	Hence if we consider a matrix representation of the set of path $\Pi(\vertex)$ with the same properties we still get a Hölder roof function.
\end{remark}

Recall that $\nu$ is the Lebesgue measure defined as in Section~\ref{proj} for a vector $q = (1, \dots, 1)$.
The following lemma is a key property to apply many thermodynamic formalism theorems.
\begin{lemma}
	\label{exponential_tail}
	The accelerated roof function $\roofacc$ has exponential tail, \textit{i.e.} there exists $0 < \sigma$ such that
	$$\int_{\psimplex} e^{\sigma \roofacc} \mathrm d \nu < \infty.$$
\end{lemma}

\begin{proof}
	Notice that there exists $C, \exptail > 0$ such that, for any $q \in \Rp$ and all $\tau > 1$,
	\begin{equation}
		\label{slice}
		\nu_q \{ x \in \psimplex \mid \roofacc(x) \geq \log \tau \} \leq C \tau^{-\exptail}.
	\end{equation}
	It follows from Corollary~\ref{proba:tail} since the above set is included in the subset of $\psimplex$ that satisfy the property $\jump^\tau \leq \pattern$.\\

	Now Formula (\ref{slice}) implies that, cutting into pieces where $\log \tau^n < \roofacc(x) \le \log \tau^{n+1}$, for all $\sigma < \exptail$,
	\begin{align*}
		\int_{\psimplex} e^{\sigma \roofacc} \mathrm d \nu_q & \leq \sum_{n=0}^\infty (\tau^{n+1})^\sigma \cdot C \cdot (\tau^n)^{-\exptail}\\
								     & = C \cdot \tau \cdot \sum_{n=0}^\infty (\tau^{\sigma-\exptail})^n = C \cdot \frac {\tau} {1 - \tau^{\sigma-\exptail}}
	\end{align*}
\end{proof}

We denote by $\sigma_0$ the supremum of such $\sigma$.
As $e^{\sigma \roofacc}$ is positive and increasing in $\sigma$, the integral for $\sigma = \sigma_0$ is infinite.

\subsubsection{Estimates on the Jacobian}
We give some useful properties on the Jacobian of the win-lose induction which follow from a computation that can be found \textit{e.g.} in \cite{Veech78}.
\begin{proposition}
	\label{jacobian:roof}
	For all $x \in \psimplex$ the Jacobian of the win-lose induction satisfies
	\[
		|D\rauzyacc (x)| = e^{|\alphabet|\roofacc(x)}.
	\]
\end{proposition}

\begin{corollary}
	\label{measure}
	There exists $Q>0$ such that for all $1$-cylinder $\psimplex_w$ and all $x \in \psimplex_w$
	$$\frac {1} {Q} \cdot |D\rauzyacc (x)|^{-1} \leq {\nu(\psimplex_{w})} \leq Q \cdot |D\rauzyacc (x)|^{-1}.$$
	In particular, there exists $Q >0$ such that for all $\kappa>0$ and $x \in \psimplex_w$
	$$\frac {1} {Q} \cdot {\nu(\psimplex_{w})}^{\kappa/|\alphabet|} \leq e^{-\kappa \roofacc(x)} \leq Q \cdot {\nu(\psimplex_{w})}^{\kappa/|\alphabet|}.$$
\end{corollary}

\begin{proof}
	If $x, y$ are in the same $1$-cylinder $\psimplex_{w}$, using (\ref{roof inequality}),
	\begin{equation*}
		|\roofacc(x) - \roofacc(y)| \leq \diam(\psimplex) < \infty.
	\end{equation*}
	Thus, using Proposition \ref{jacobian:roof}, there exists $Q'>0$ such that
	\begin{equation}
		\label{distortion_ratio}
		\frac {1} {Q'} \cdot |D\rauzyacc(x)| \leq |D\rauzyacc(y)| \leq Q' \cdot |D\rauzyacc(x)|.
	\end{equation}
	The restriction $\rauzymap_{*w} := {\rauzyacc}_{|\psimplex_w}$ is invertible and for all $v \in \simplex$ we have $y := \rauzymap_{*w} v \in \simplex$.
	Thus, integrating $v$ over $\simplex$,
	$$\frac {1} {Q'} \cdot |D\rauzyacc (x)|^{-1} \cdot \nu(\psimplex) \leq \int_{\psimplex} |D\rauzyacc(\rauzymap_{*w}^{-1} v)|^{-1} d \nu(v) \leq Q' \cdot |D\rauzyacc (x)|^{-1} \cdot \nu(\psimplex)$$
	and
	$$\int_{\psimplex} |D\rauzyacc(\rauzymap_{*w}^{-1} v)|^{-1} d \nu(v) = \int_{\psimplex} |D\rauzymap_{*w}^{-1}| d \nu = \nu(\psimplex_w).$$
\end{proof}

\begin{remark}
	\label{mixing}
	The argument also works for a $n$-cylinders indexed by $\mathbf w = [w_1, \dots, w_n]$ and $\rauzyacc^n$.
	When integrating with respect to $v \in \psimplex$ changing variable $x = (\rauzymap_{* \mathbf w}^n)^{-1} v$ this implies that for all finite cylinders $\mathbf {w_1}, \mathbf {w_2}$,
	$$ \frac {1} {Q} \le \frac {\nu(\psimplex_{\mathbf {w_1} \mathbf {w_2}})} {\nu(\psimplex_{\mathbf {w_1}}) \nu(\psimplex_{\mathbf {w_2}})} \le Q.$$
	Which is exactly the bounded distortion property of \cite{AvilaForni07}.
\end{remark}

\begin{corollary}
	\label{invariant_measure}
	There exists a unique ergodic $\rauzyacc$-invariant Borel probability measure $\mu$ absolutely continuous with respect to Lebesgue measure $\nu$.
	Moreover the logarithm of its density $|\log \frac {d\mu} {d\nu}|$ is bounded by a constant at almost every point.
\end{corollary}

\begin{proof}
	This is a direct application of Lemma 4.4.1 in \cite{Aaronson97} which shows existence of an invariant measure absolutely continuous with respect to Lebesgue measure with an almost everywhere bounded logarithm of density.
	Remark~\ref{mixing} implies ergodicity of $\rauzyacc$ for Lebesgue measure and thus for the invariant measure.
	Remark~\ref{recursion} implies conservativity with respect to Lebesgue thus there exists a unique absolutely continuous invariant measure (see for instance Theorem~1.5.6 in \cite{Aaronson97}).
\end{proof}

\subsubsection{Gibbs measures and Gurevic--Sarig pressure.}

\begin{definition*}
	Let $\mu$ be a $\sigma$-invariant (shift) Borel probability measure on a countable Markov chain $\Sigma$.
	For any continuous function $\phi : \Sigma \to \R$, $\mu$ will be called a \textit{Gibbs measure for the potential $\phi$} if there exists $Q>0$ and $P$ such that for every path $\finpath = \finpath_0 \cdot \finpath_n$ and every $x$ in the cylinder $[x_1, \dots, x_n]$
	\begin{equation}
		\label{Gibbs}
		\frac {1} {Q} \leq \frac {\mu\left([x_1, \dots, x_n]\right)} {\exp \left(\sum_{k=0}^{n-1} \phi(\sigma^k (x)) - Pn\right)} \leq Q.
	\end{equation}
	$P$ is unique and is called the topological pressure of $\phi$.
\end{definition*}

In the following we consider the potential functions $\phi_\kappa = - \kappa r$ for  $\kappa \ge 0$.
When there is no ambiguity we will denote one of these functions simply by $\phi$.
We will show that they satisfie good properties to induce existence and uniqueness of Gibbs measures.\\

We need some further definitions in order to introduce existence theorems.
As we deal with a coding on a countable alphabet, we first need to check a technical property on this coding called "Big Image and big Preimage".
\begin{definition*}
	The \textit{BIP property} is the existence of $w_1, \dots, w_m \in \states$ tiles of the Markov partition, such that for all $v \in \states$, there exists $1\leq k,l \leq m$ such that $\rauzyacc w_k \cap v$ and $\rauzyacc v \cap w_l$ are not empty.
\end{definition*}
This property is obviously true in our case since each Markov tile is sent to the whole domain by $\rauzyacc$.\\

The Ruelle operator $L_\phi$ associated to a potential function $\phi$ is an operator acting on the space of continuous functions.
For a function $f$ on $\psimplex$ and $x \in \psimplex$ it is defined by
$$(L_\phi f)(x) = \sum_{\rauzyacc(y)=x} e^{\phi(y)} f(y).$$
As explained in \cite{Sarig15}: "the analysis of thermodynamic limits reduces to the study of the asymptotic behavior of $L^n_\phi f$ as $n \to \infty$ for \textit{sufficiently many} functions $f$ ".
One of the key to understand this behavior is to first understand the limit of $\frac {1} {n} \log L^n_\phi f$.
In particular, it can be compared to the following quantities.\\

For $w \in \states$, let
	$$Z_n(\phi, w) = \sum_{\rauzyacc^n(x) = x, \, x_0=w} e^{\phi_n(x)},$$
with $\phi_n = \phi + \phi \circ \rauzyacc + \dots + \phi \circ \rauzyacc^{n-1}$ and $x_0$ is the tile in $S$ to which $x$ belongs.
According to Theorem 4.3 in \cite{Sarig15} the limit
\begin{equation}
	\label{eq_limit}
	P_G(\phi) := \lim_{n \to \infty} \frac {1} {n} \log Z_n(\phi, w)
\end{equation}
exists for all $w \in \states$ and is independent of $w$.
Moreover, if $\|{L_\phi 1}\|_\infty < \infty$, then $P_G(\phi) < \infty$.

\begin{definition*}
	$P_G(\phi)$ is called the Gurevic--Sarig pressure of $\phi$.
\end{definition*}
This is a relevant quantity to consider according to Theorem 4.4 of \cite{Sarig15} since, when $P_G(\phi)$ is finite, it is equal to the limit of $\frac {1} {n} \log L^n_\phi f$ for a large class of functions.
It is not always the case for characteristic functions $1_{[w]}$ but we can still make the following remark.
\begin{remark}
	\label{pressure_lower_bound}
	As a consequence of Sarig's generalized Ruelle--Perron--Frobenius (Theorem~4.9 in \cite{Sarig15}) for all $w \in \states$ and $x \in \psimplex$
	$$P_G(\phi) \ge \lim_{n \to \infty}\frac {1} {n} \log \left(L^n_\phi 1_{[w]}\right)(x).$$
\end{remark}

\begin{definition*}
	The potential function $\phi$ has bounded variations if and only if $\sum_{n=2}^\infty var_n(\phi) < \infty$ where
	$$var_n(\phi) = \sup \{ |\phi(x) - \phi(y)| : x_i = y_i, i = 1, \dots, n \}.$$
\end{definition*}
Notice that the Hölder property proved in Lemma~\ref{holder} implies that, for all $\kappa$, $\phi$ has bounded variations and $var_1(\phi) < \infty$.\\

These definitions enable us to state the key theorem in this section.
It gives a criterion for uniqueness of a Gibbs measure for a given potential function.
The following formulation of Sarig theorem is taken from Theorem~4.6 in \cite{Pesin14}.
\begin{theorem*}[Sarig \cite{Sarig03}]
Assume that the potential $\phi$ has summable variations.
Then $\phi$ admits a unique $\rauzyacc$-invariant Gibbs measure $\mu_\phi$ if and only if
\begin{itemize}
	\item  X satisfies the BIP property;
	\item the Gurevic--Sarig pressure $P_G(\phi) < \infty$ and $var_1 \phi < \infty$.
\end{itemize}
In this case, the topological pressure and Gurevic--Sarig pressure coincide.
\end{theorem*}
The hypothesis of this theorem are satisfied for large $\kappa$.
Indeed the following lemma which implies finiteness of the pressure in this case.
\begin{lemma}
	\label{bounded}
	The pressure $P_G(\phi)$ is finite if and only if $\kappa > |\alphabet| - \sigma_0$.
\end{lemma}

\begin{proof}
	As noticed above, in order to show that $P_G(\phi)$ is finite we only need to prove that $L_\phi 1$ is finite for all $\kappa > |\alphabet|-\sigma_0$.
	By definition,
	\[
		 L_\phi 1 = \sum_{\rauzyacc(y) = x} e^{\phi(y)} \le (Q')^{\kappa/|\alphabet|} \cdot \sum_{w \in \states} e^{-\kappa \roofacc(w)}.
	\]
	Where the last inequality is a consequence of Formula~(\ref{distortion_ratio}) and $S$ is a choice of representative of every $1$-cylinders.

	\paragraph{Necessary condition.}
	Finiteness of the pressure thus follows from the following proposition.
	\begin{proposition}
		\label{prop_bounded_pressure}
		There exists a constant $K$ such that for all $0< \sigma < \sigma_0$ and all $\kappa > |\alphabet| - \sigma$
		\begin{equation}
			\label{bounded_pressure}
			\sum_{w \in \states} e^{-\kappa \roofacc(w)} \le K \cdot \frac {I_\sigma} {1 - e^{|\alphabet| - \sigma - \kappa}}.
		\end{equation}
	\end{proposition}

	\begin{proof}
	Let $Y(N)$ be the set of representatives $w \in \states$ for which $N \leq \roofacc(w) < N+1$ then
	\begin{equation}
		\label{geometric}
		\sum_{w \in \states} e^{-\kappa \roofacc(w)} = \sum_{N=0}^\infty \ \sum_{w \in Y(N)} e^{-\kappa \roofacc(w)} \leq \sum_{N=0}^\infty |Y(N)|  e^{-\kappa N}.
	\end{equation}
	Using Formula~(\ref{distortion_ratio}) and Lemma~\ref{exponential_tail}, for all $0 < \sigma < \sigma_0$,
	\[
	\sum_{w \in Y(N)} \int_{\psimplex_w} e^{\sigma \roofacc(w)} \mathrm d \nu \le Q' \cdot I_\sigma < \infty.
	\]
	where $I_\sigma$ stands for $\int_{\psimplex} e^{\sigma \roofacc} \mathrm d \nu$.
	And
	\[\sum_{w \in Y(N)} \int_{\psimplex_w} e^{\sigma \roofacc(w)} \mathrm d \nu \geq e^{\sigma N} \sum_{w \in Y(N)} \nu(\psimplex_w).\]
	Moreover, according to Corollary~\ref{measure}, for all $w \in Y(N)$,
	$$\nu(\psimplex_w) \geq Q^{-1} \cdot e^ {-|\alphabet|(N+1)}.$$
	Hence
	\[
		I_\sigma \ge (Q\cdot Q')^{-1} \cdot |Y(N)| \cdot e^{\sigma N - |\alphabet|(N+1)},
	\]
	and
	\[
		|Y(N)| \le Q \cdot Q' \cdot I_\sigma \cdot e^{|\alphabet|} \cdot e^{(|\alphabet|-\sigma) N}.
	\]
	Thus the geometric sum in (\ref{geometric}) is bounded for $\kappa > |\alphabet| - \sigma$ by
	\begin{equation*}
		Q \cdot Q' \cdot {I_\sigma} \cdot e^{|\alphabet|} \cdot \sum_{N=0}^\infty e^{(|\alphabet| - \sigma - \kappa)N} = Q \cdot Q' \cdot {I_\sigma} \cdot \frac{e^{|\alphabet|}} {1 - e^{|\alphabet|-\sigma-\kappa}}.
	\end{equation*}
	\end{proof}

	\paragraph{Sufficient condition}
	Using Remark~\ref{pressure_lower_bound} it will be enough to show that whenever the tail integral is infinite, for some $x \in \psimplex$ and $w \in \states$, we have
	$$\frac {1} {n} \log\left(L_{\phi}^n 1_{[w]}\right)(x) \to \infty.$$
	Let $w_k$ be the representative of the $1$-cylinder containing $\rauzyacc^k y$ in $S$.
	By Formula~(\ref{distortion_ratio}) for $w$ in $S$,
	\begin{align*}
		L_\phi^n 1_{[w]} (x) &= \sum_{\rauzyacc^n (y) = x} e^{-\kappa \cdot \roofaccn(y)} 1_{[w]}(y) \\
				     &\ge (Q')^{-n \cdot {\kappa}/{|\alphabet|}} \sum_{w_1, \dots, w_{n-1} \in \states} e^{-\kappa \cdot \roofacc(w)} e^{-\kappa \cdot \roofacc(w_1)} \dots e^{-\kappa \cdot \roofacc(w_{n-1})}.
	\end{align*}
	Hence, we only need to show that for $\kappa \le |\alphabet| - \sigma_0$
	$$\sum_{w \in \states} e^{-\kappa \roofacc(w)} = \infty.$$
	As previously we split the sum
	\begin{equation*}
		\sum_{w \in \states} e^{-\kappa \roofacc(w)} = \sum_{N=0}^\infty \ \sum_{w \in Y(N)} e^{-\kappa \roofacc(w)} \ge \sum_{N=0}^\infty |Y(N)|  e^{-\kappa (N+1)}.
	\end{equation*}
	Now
	\begin{align*}
		I_\sigma & < Q \cdot \sum_{N=0}^\infty |Y(N)| \cdot e^{\sigma (N+1) - |\alphabet|N}\\
			 & = Q \cdot e^\sigma \cdot \sum_{N=0}^\infty |Y(N)| \cdot e^{-(|\alphabet| - \sigma)N}.
	\end{align*}
	In particular
	$$
		\sum_{N=0}^\infty |Y(N)| \cdot e^{-(|\alphabet| - \sigma_0)N} = \infty.
	$$

	\begin{proposition}
		\label{existence}
		For all $\kappa> |\alphabet| - \sigma_0$ there exists a unique Gibbs measure $\mu_\phi$ of potential $\phi = -\kappa \cdot \roofacc$.
		Moreover this measure has finite entropy.
	\end{proposition}

	\begin{proof}
		The existence of the Gibbs measure is a direct consequence of the Hölder property in Lemma~\ref{holder}, Lemma~\ref{bounded} and Theorem~1 in \cite{Sarig03} quoted above.

		Using inequality $x+1 \le e^x$, for $\sigma$ as in Lemma~\ref{exponential_tail},
		$$ \int_{\psimplex} \sigma \roofacc d\mu_\phi \le \int_{\psimplex} (e^{\sigma \roofacc} - 1) d \mu_\phi \le \left\|\frac {d\mu_\phi}{d\nu}\right\|_\infty \int_{\psimplex} e^{\sigma \roofacc} d \nu < \infty.$$
		Hence $\int_{\psimplex} - \phi d\mu_\phi < \infty$ and by Lemma~\ref{bounded} $P_G(\phi) < \infty$.
		By Theorem~2 in \cite{Sarig01a} this implies that the entropy of $\mu_\phi$ is finite.
	\end{proof}
\end{proof}

\subsection{Suspension flow}
	\label{naturalext}
	There is an easy way to construct a natural extension for a full shift on a countable alphabet by extending it to bi-infinite words.
	The canonical suspension then extends to a flow on the suspension of the natural extension.\\

	As in Section~\ref{suspension}, we consider for $f : \parameter \to \Rp$ the suspensions flow $\Phi$ on $\suspension_f$.
	Denote by $\mathcal M_{T,f}$ the set of $\rauzymap$-invariant Borel probability measures with $\mu(f) := \int_\parameter f d\mu < \infty$.
	Every $\Phi$-invariant Borel probability measure $\widetilde \mu$ on $\suspension_f$ can be decomposed as a product of a measure $\mu \in \mathcal M_{T,f}$ and the Lebesgue measure on fibers.
	Namely,
	$$\widetilde \mu_f = \left(\mu(f)\right)^{-1} (\mu \times \lambda)_{|\suspension_f}.$$

	The Kolmogorov--Sinai entropy of the flow for this measure is written $h(\Phi, \widetilde \mu)$ and satisfies Abramov's formula
	$$h(\Phi, \widetilde \mu) = \frac {h(T, \mu)} {\mu(f)}$$
	where $h(T,\mu)$ is the Kolmogorov--Sinai entropy for $T$.
	In this setting the topological entropy can be defined as
	$$h_{top}(\Phi) = \sup_{\mu \in \mathcal M_{T,f}} h(\Phi, \widetilde \mu_f).$$
	The induced measure $\widetilde \mu_f$ for $\mu \in \mathcal M_{T,f}$ at which this supremum is achieved (and by extension $\mu$ itself) is referred to as a \textit{measure of maximal entropy}.\\

	We will first consider the canonical suspension of $T$ by $\roof$ and in the next subsection by a perturbed roof function.
	As noticed previously, the suspension flow on $\parameter$ for the roof function and the one on $\psimplex$ for the accelerated roof function are conjugate.
	Moreover, the exponential tail integrals are equal for these two suspension flows with the same measure on the base restricted to $\psimplex$.

	We use from now on the representation of the suspension on the base $\psimplex$ with the map $\rauzyacc$ conjugated to a full shift.

\begin{proposition}
	\label{proposition:kappa}
	The pressure $P_G(-\kappa \cdot \roofacc)$ vanishes at a unique value $\kappa = |\alphabet|$.
	The associated Gibbs measure $\mu_\phi$, as in the previous proposition, is the unique measure of maximal entropy for the canonical suspension of the win-lose induction.
\end{proposition}

\begin{proof}
	Let $\mu$ be a measure as in Corollary~\ref{invariant_measure}, we show that it is the unique Gibbs measure for potential $- |\alphabet| \cdot \roofacc$.
	According to Corollary~\ref{measure}, there exists a constant $Q > 0$ such that
	$$\frac {1} {Q} \le \frac{\nu([x_1, \dots, x_n])}{\exp \left(\sum_{k=0}^{n-1} - |\alphabet| \cdot r\left(\sigma^k(x)\right)\right)} \le Q.$$
	As $\mu$ is such that $|\log \frac {d\mu} {d\nu}|$ is bounded at almost every point then it also satisfies the same property for another constant $Q$.
		Thus it is by definition a Gibbs measure for the potential $- |\alphabet| \cdot \roof$ and this potential has zero topological pressure.
		The function $P_G(-\kappa \cdot \roofacc)$ hence vanishes at $|\alphabet|$.

	Moreover, $P_G\left(-\left(|\alphabet| - \sigma_0\right) \cdot \roofacc\right) = \infty$ by Lemma~\ref{bounded} and $P_G$ is a convex decreasing continuous function of $\kappa$ (see Theorem 4.6 of \cite{Sarig15}).
	Thus $|\alphabet|$ is the unique value such that $P_G(-\kappa \cdot \roofacc) = 0$.

	\begin{figure}[h!]
		\center
		\includegraphics[scale=.7]{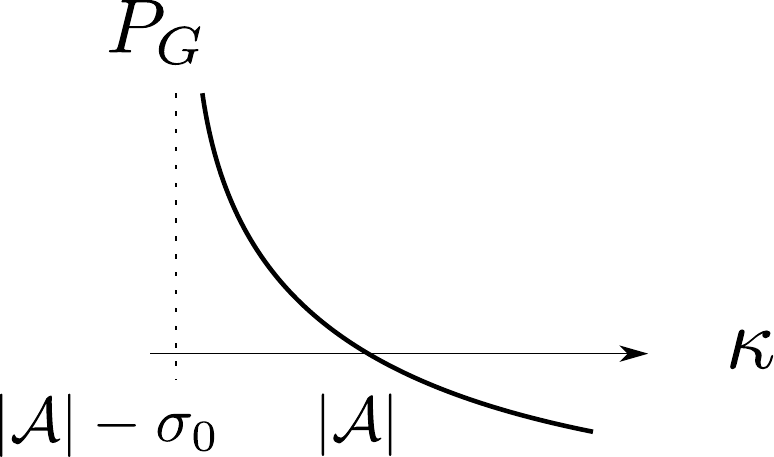}
	\end{figure}

	As $\mu$ has finite entropy it is an equilibrium measure and it satisfies the \textit{variational principle} for the topological pressure (see Section~5.3 in \cite{Sarig15}).
	Thus $\mu$ maximizes, over all Borel measures such that $\mu(\phi) > -\infty$ \textit{i.e.} $\mu(\roofacc) < \infty$, the quantity
	$$h(\rauzyacc, \mu) - \int_{\parameter} \kappa_0 \roofacc \mathrm d \mu.$$
	The maximum is equal to the pressure, here 0.
	Hence $h(\Phi, \widetilde \mu) = \frac {h(\rauzyacc, \mu)} {\mu_0(r)} = |\alphabet|$ is maximal.
	Theorem 1.1 in \cite{BuzziSarig03} tells us that there is at most one such maximizing measure.
\end{proof}

We summarize these results in the following theorem.
\begin{theorem}
	\label{entropy}
	The measure of maximal entropy for the suspension of a \unstable\ win-lose induction is the suspension of the unique $\rauzyacc$-invariant Borel measure $\mu$ absolutely continuous with respect to Lebesgue measure introduced in Corollary~\ref{invariant_measure}.\\
	Moreover, the entropy of the measure on the suspension is equal to $|\alphabet|$.
\end{theorem}

	\subsection{Subgraph parameter space}
	\label{fractals}

	Let $\subgraph$ be a \unstable\ subgraph of $\graph$.
	Similarly to $\parameter$ we denote by $\simplex^\graph(\subgraph) \subset \parameter$ or simply $\subparameter$ the subset of points whose path belongs to $\subgraph$.
	Using the coding functions $c_\vertex$ introduced in Section~\ref{SS}, there is a natural bijection, defined outside of a countable union of codimension one subsets,
	$$\iota : \simplex^\subgraph \to \simplex^\graph(\subgraph).$$

	Let us assume that $\subgraph$ is strongly connected and all letters in $\alphabet$ appear as a label in the subgraph.
	Then there exists a positive path $\ppath \in \subgraph$ which we will use to define an acceleration of the win-lose induction.
	Notice that
	$$\roofacc_\graph \circ \iota = \roofacc_\subgraph + \delta$$
	where $\roofacc_\graph$ and $\roofacc_\subgraph$ are the $\ppath$-accelerated roof functions on $\graph$ and $\subgraph$ respectively and $\delta$ is a non-negative function.
	For convenience in the writing of the proofs we will simply denote $\roofacc_\graph$ by $\roofacc$.\\

	The following theorem bounds the Hausdorff dimension of $\subparameter$ in terms of the value at which the pressure of a family of potentials vanishes.
	This may recall of Bowen's computation of Hausdorff dimensions of quasicircles in \cite{Bowen79}.\\

	\begin{theorem}
		\label{fractal:measure}
		Let $\subgraph$ be a strongly connected \unstable\ subgraph of $\graph$ which has the same set of labels, then the Hausdorff dimension of $\subparameter$ satisfies
		$$ \dim_H \subparameter \le |\alphabet| - 2 + \frac{\kappa_\subgraph}{|\alphabet|}$$
		where $\kappa_\subgraph$ is the unique positive real number satisfying
		$$P_G\left(-\kappa_\subgraph \cdot (\roofacc_\subgraph + \delta)\right) = 0.$$
	\end{theorem}

	\begin{proof}
		The map $\roofacc_\subgraph + \delta$ can be seen as a roof function associated to an alternative non-negative matrix representation $M'$ of $\Pi(\vertex)$ such that if $M_\ppath > 0$ then $M'_{\ppath} > 0$.
		By Remark~\ref{holder2} it is also a Hölder roof function as in Lemma~\ref{existence}.
		Thus there exists a Gibbs measure $\widetilde \mu$ associated to the potential $- \kappa \cdot (\roofacc_\subgraph + \delta)$ for all $\kappa$ such that the pressure is finite.\\

		As $\delta \ge 0$, the pressure satisfies for all $\kappa \ge 0$
		$$P_G \left(- \kappa \cdot (\roofacc_\subgraph + \delta)\right) \le P_G\left( - \kappa \cdot \roofacc_\subgraph\right).$$
		Thus $P_G \left(- |\alphabet| \cdot (\roofacc_\subgraph + \delta)\right) \le 0$.
		As $P_G(0) = \infty$, by continuity and strong convexity of the pressure, there exists a unique $\kappa_\subgraph$ such that $P_G\left(- \kappa_\subgraph \cdot (\roofacc_\subgraph + \delta)\right) = 0$.\\

		We consider the pushed forward measure $\mu := \iota_* \widetilde \mu$ which is also a Gibbs measure for the potential equal to $- \kappa_\subgraph \cdot \roofacc$ on $\subparameter$ and $0$ elsewhere.
		According to Formula (\ref{Gibbs}) in the definition of Gibbs measures, there exists $Q > 0$ such that for all $x$ in the intersection of the cylinder ${\cylinder} = [w_1, \dots, w_n]$ and $\subparameter$
		\[
			\frac {1} {Q} \cdot \exp \left(- \sum_{k=0}^{m-1} \kappa_\subgraph \cdot \roofacc(\rauzyacc^k (x))\right) \leq {\mu(\psimplex_{\cylinder})}  \leq Q \cdot \exp \left(- \sum_{k=0}^{m-1} \kappa_\subgraph \cdot \roofacc(\rauzyacc^k (x))\right).
		\]

		We use for convenience the notations $a \lesssim b$ and $a \simeq b$ to state that there exists a constant $Q$ depending only on the choice of graph and subgraph such that $a \le Q \cdot b$ and $\frac {1} {Q} \cdot b \le a \le Q \cdot b$ respectively.\\

		Corollary~\ref{measure} implies that,
		\[
			\exp \left(|\alphabet| \sum_{k=0}^{m-1} \roofacc(\rauzyacc^k (x))\right) = |D\rauzyacc^m (x)| \simeq \frac {1} {\nu(\psimplex_{\cylinder})}.
		\]
		Thus, for all cylinders intersecting $\subparameter$,
		\begin{equation}
			\label{mugibbs}
			\mu(\psimplex_{\cylinder}) \simeq {\nu(\psimplex_{\cylinder})}^{\kappa_\subgraph / |\alphabet|}.
		\end{equation}
		Let us introduce the notation $\alpha := \kappa_\subgraph / |\alphabet|$.\\

		\def \fam {\mathcal F}
		\def \famc {\fam_n^{\le C}}
		\def \famcc {\fam_n^{> C}}
		For $\epsilon > 0$, let $\fam$ be a family of cylinders for $\rauzyacc$ that intersect $\subparameter$ and such that $\nu(\psimplex_{\cylinder}) < \epsilon$.
		Such a family exists since $\rauzyacc$ is uniformly expanding by Proposition~\ref{uniform}.\\

		As noticed in \cite{AvilaDelecroix16}, simplices with bounded diameter satisfy a property which relates their Lebesgue measure with the number of balls necessary to cover them.
		This will be very useful to bound Hausdorff dimensions.
		\begin{proposition}
			\label{minimal}
			There exists $K>0$ such that for all simplex $\simplex$ of dimension $d$, measure $m$ and diameter less than $1$, the minimal number of balls of radius $0 < \rho \le m$ required to cover $\simplex$ satisfies
			$$N_\rho \le K \cdot \frac {m} {\rho^d}.$$
		\end{proposition}
		In the case of $m = \rho$ this implies $N_\rho \le K \cdot \rho^{1-d}$.
		Thus for all $\cylinder \in \fam$ one can find a covering $\{B_i\}$ by less than $K \cdot \nu\left(\psimplex_{\cylinder}\right)^{1-d}$ balls of radius $\nu \left( \psimplex_{\cylinder} \right) < \epsilon$.
		For this covering we then have
		$$\sum_i (\diam B_i) ^\delta \le \sum_{\cylinder \in \fam} K \cdot \nu \left(\psimplex_{\cylinder}\right)^{1-d} \cdot \nu \left(\psimplex_{\cylinder}\right)^\delta.$$
		By Formula~\ref{mugibbs},
		$$\sum_{\cylinder \in \fam} \nu \left(\psimplex_{\cylinder}\right)^{\kappa_\subgraph / |\alphabet|} \simeq \sum_{\cylinder \in \fam} \mu \left(\psimplex_{\cylinder}\right) = \mu \left(\bigcup_{\cylinder \in \fam} \psimplex_{\cylinder}\right).$$
		Thus if $1-d+\delta \ge \kappa_\subgraph / |\alphabet|$ then $\sum_i (\diam B_i) ^\delta$ is bounded uniformly for all $\epsilon$.
		Then $\dim_H \simplex(\subgraph) \le d - 1 + \kappa_\subgraph / |\alphabet|$.
	\end{proof}

	To continue our study we need to have more control on the induced measure on the set $\subparameter$.
	\begin{proposition}
		\label{finite_entropy_tilde}
		The measure $\widetilde \mu$ introduced in the previous theorem has finite entropy.
	\end{proposition}

	\begin{proof}
		As in Proposition~\ref{existence} we only need to prove that $\widetilde \mu(\roofacc \circ \iota) < \infty$.
		Similarly to the proof of Proposition~\ref{prop_bounded_pressure} let us consider $Y(N)$ the set of representatives $w \in \states$ (for the acceleration on $\subgraph$) for which $N \leq \left(\roofacc \circ \iota\right) (w) < N+1$, then
	\begin{align*}
		\int_{(\simplex^\subgraph)^*} \roofacc \circ \iota \, d\mu & \le \sum_{N = 0}^\infty \sum_{w \in Y(N)} (N+1) \cdot \widetilde \mu(\psimplex_w)\\
									   & \lesssim \sum_{N = 0}^\infty (N+1) \cdot |Y(N)|  \cdot e^{-\kappa_\subgraph \cdot N}.
	\end{align*}

	Recall that there exists $\sigma > 0$ such that
	\[
		I_\sigma := \int_{(\simplex^\subgraph)^*} e^{-\sigma (\roofacc_\subgraph + \delta)} d\nu \le \int_{(\simplex^\subgraph)^*} e^{-\sigma \roofacc_\subgraph} d \nu < \infty.
	\]
	Let $\sigma_0$ be the supremum of such $\sigma$.
	By the same argument as in Lemma~\ref{bounded}, the pressure of $- \kappa \cdot \roofacc \circ \iota$ is infinite for $\kappa = |\alphabet| - \sigma_0$ thus $\kappa_\subgraph > |\alphabet| - \sigma_0$.
	And
	\[
		|Y(N)| \lesssim I_\sigma \cdot e^{(|\alphabet| - \sigma) N}.
	\]
	Hence
	\[
		\int_{(\simplex^\subgraph)^*} \roofacc \circ \iota \, d\mu \lesssim I_\sigma \cdot \sum_{N = 0}^\infty N \cdot e^{(|\alphabet|-\sigma)N} \cdot e^{-\kappa_\subgraph \cdot N}.
	\]
	Choosing $\sigma < \sigma_0$ such that $\kappa_\subgraph > |\alphabet|-\sigma$ we have a converging sum.
	\end{proof}

	\begin{corollary}
	\label{fractal:dimension}
	Let $\subgraph$ be a subgraph of $\graph$ as in the previous theorem such that there is a vertex in $\subgraph$ for which one outgoing edge is in $\graph$ and not in $\subgraph$.
	The Hausdorff dimension of the parameter subset $\subparameter$ is strictly smaller than the dimension of $\parameter$.
	\end{corollary}

\begin{proof}
	As we saw in the previous theorem, the pulled back roof function  $\roofacc \circ \iota$ is also Hölder and thus has bounded variations.
	Moreover, as $\delta \ge 0$ the pressure satisfies
	$$P_G\left(-|\alphabet| \cdot (\roofacc_\subgraph + \delta)\right) \le P_G\left(-|\alphabet| \cdot \roofacc_\subgraph\right) = 0.$$
	Thus the pressure for the potential $-\kappa \cdot (\roofacc \circ \iota)$ is zero for a unique $\kappa_\subgraph \le |\alphabet|$.
	We prove in the following that the former inequality is strict.\\

	Let $\widetilde \mu_0$ and $\widetilde \mu_1$ be the Gibbs measures associated to potentials $-\kappa_\subgraph \cdot (\roofacc \circ i)$ and $-|\alphabet| \cdot \roofacc_\subgraph$ respectively.
	We have proved that these two measures have finite entropy, thus by the variational principle and Buzzi--Sarig theorem they are the unique measures satisfying respectively
	$$P_G\left(-\kappa_\subgraph \cdot (\roofacc_\subgraph + \delta)\right) = h(\rauzyacc^\subgraph, \widetilde \mu_0) - \int_{\simplex_\subgraph} \kappa_\subgraph \cdot (\roofacc_\subgraph + \delta) \ d \widetilde \mu_0$$
	and
	$$P_G\left(-|\alphabet| \cdot \roofacc_\subgraph\right) = h(\rauzyacc^\subgraph, \widetilde \mu_1) - \int_{\simplex_\subgraph} |\alphabet| \cdot \roofacc_\subgraph \ d \widetilde \mu_1.$$
	Each measure maximizes the right-hand side quantity and the maximum value is zero in both cases.\\

	Assume $\kappa_\subgraph = |\alphabet|$ then
	$$ h(\rauzyacc^\subgraph, \widetilde \mu_0) - \int_{\simplex_\subgraph} \kappa_\subgraph \cdot (\roofacc_\subgraph + \delta) \ d \widetilde \mu_0
	\le  \left(h(\rauzyacc^\subgraph, \widetilde \mu_1) - \int_{\simplex_\subgraph} |\alphabet| \cdot \roofacc_\subgraph \ d \widetilde \mu_1\right) - \kappa_\subgraph \cdot \widetilde \mu_0 (\delta) .$$
	Thus $\widetilde \mu_0 (\delta) \le 0$ and as $\delta \ge 0$ this implies that $\delta = 0$ almost everywhere.\\

	Let us pick a cylinder associated to a path that passes through a vertex to which we have removed an edge.
	The map $\delta$ is positive on this cylinder which has positive measure.
	This is a contradiction, hence $\kappa_\subgraph < |\alphabet|$.
\end{proof}

The proof of Theorem~\ref{fractal:measure} together with Proposition~\ref{finite_entropy_tilde} also imply a result on the existence of a measure of maximal entropy similar to Theorem~\ref{entropy}.
\begin{theorem}
	\label{entropy_fractal}
	The measure of maximal entropy for the restriction of the suspension of the win-lose induction defined on $\parameter$ to a parameter subset $\subparameter$, where $F$ has \unstability, is the product of the unique Gibbs measure associated to $\kappa_\subgraph$ and the Lebesgue measure on fibers.
	Moreover, its entropy is equal to $\kappa_\subgraph$.
\end{theorem}

	Notice in particular that this gives an intrinsic definition of $\kappa_\subgraph$ as the topological entropy of the suspension flow restricted to $\subparameter$ which hence does not depend on the choice of acceleration.

\newpage
\section{Continued fraction algorithms}
\label{examples}

In this section, we describe how to associate to a large set of examples of linear simplex-splitting MCF algorithm (in the sense of Lagarias \cite{Lagarias93}) a conjugate simplicial system.
This hopefully will make the general algorithm clear.
We are able to check the \unstability\ using the criterion introduced in Section~\ref{criterion} for all known ergodic algorithms we are considering.
The only limit case in which our criterion does not apply is given by the Poincaré algorithm in dimensions larger or equal to $4$.  \\

We will start with two simple examples, the fully subtractive and Poincaré algorithms, for which it is easy to derive from their classical description an associated simplicial system.
One of the reason that make these examples easier to describe in terms of simplicial systems is the fact that their domains of definition are all sent to the whole simplex by the corresponding map.\\

We then present a general strategy to compute these simplicial systems and apply it to Brun and Selmer algorithms.
We finish by computing a simplicial system which induces the Rauzy gasket in every dimension.
This will induce in particular a simplicial system description of Arnoux--Rauzy--Poincaré algorithm.\\

As a consequence we have a unified proof that Brun and Selmer and Arnoux--Rauzy--Poincaré algorithms are ergodic for their unique invariant measure equivalent to Lebesgue. Moreover this measure induces the unique measure of maximal entropy on their canonical suspension.\\

Ergodicity for Brun and Selmer algorithms in all dimension is due to Schweiger \cite{Schweiger00}, for Arnoux--Rauzy--Poincaré it has been proved in \cite{BertheLabbe13}.
The result on Hausdorff dimension has been proved in dimension 2 in \cite{AvilaHubertSkripchenko16b}.

\subsection{Two full-image examples}

\subsubsection{Fully subtractive algorithms}
The fully subtractive algorithm in dimension $3$ can be described by the map, defined at almost every point,
$F: (x_1, x_2, x_3)\in\Rp^3 \to (x'_1, x'_2, x'_3),$
where if $\{i,j,k\} = \{1,2,3\}$ and $x_i > x_j > x_k$,
\begin{equation*}
	x'_i = x_i - x_k,\
	x'_j = x_j - x_k,\
	x'_k = x_k.
\end{equation*}

This map corresponds to a step for the win-lose induction in the graph with one vertex and three edges of distinct labels, represented below.

\begin{figure}[!htb]
	\center
	\begin{tikzpicture}[shorten >=1pt,node distance=2cm,on grid,auto]
	   \node (base) {$\circ$};
	   \path[->] (base) edge[loop, out=60, in =150, looseness=13]  node[above=.2cm] {$1$} (base);
	   \path[->] (base) edge[loop, out=180, in =270, looseness=13]  node[left=.2cm] {$2$} (base);
	   \path[->] (base) edge[loop, out=300, in =20, looseness=13]  node[right=.2cm] {$3$} (base);
        \end{tikzpicture}
\vspace*{-.7cm}
\end{figure}

This first example has stable subgraphs, in which the orbits will eventually be trapped.
This corresponds to the behavior proved in \cite{Nogueira95} for the 3-dimensional Poincaré algorithm where one coordinate remains much bigger than the two others which decrease very fast by applying a continued fraction algorithm to them.\\

This construction generalizes to fully subtractive algorithms in dimension $n > 3$ by taking a single edge with $n$ loops labeled by $n$ different letters.\\

\subsubsection{Poincaré algorithms}
Poincaré algorithm has been introduced by Poincaré as a generalization of the continued fraction algorithm and was later studied and generalized in \cite{Nogueira95}.
It can be described by the map $$F: (x_1, x_2, x_3)\in\Rp^3 \to (x'_1, x'_2, x'_3),$$
where if $\{i,j,k\} = \{1,2,3\}$ and $x_i > x_j > x_k$,
\begin{equation*}
	x'_i = x_i - x_j,\
	x'_j = x_j - x_k,\
	x'_k = x_k.
\end{equation*}

\begin{figure}[!htb]
	\center
	\begin{tikzpicture}[shorten >=1pt,node distance=2cm,on grid,auto]
	   \node (base) {$\circ$};
	   \node (2) [below=of base] {$\bullet$};
	   \node (1) [left=4cm of 2] {$\bullet$};
	   \node (3) [right=4cm of 2] {$\bullet$};
	   \node (12) [below right= of 1] {$\circ$};
	   \node (13) [below left= of 1] {$\circ$};
	   \node (21) [below right= of 2] {$\circ$};
	   \node (23) [below left= of 2] {$\circ$};
	   \node (31) [below right= of 3] {$\circ$};
	   \node (32) [below left= of 3] {$\circ$};

	   \path[->] (base) edge node[above left] {$1$} (1);
	   \path[->] (base) edge node[left] {$2$} (2);
	   \path[->] (base) edge node[above right] {$3$} (3);

	   \path[->] (1) edge node[above right] {$2$} (12);
	   \path[->] (1) edge node[above left] {$3$} (13);
	   \path[->] (2) edge node[above right] {$1$} (21);
	   \path[->] (2) edge node[above left] {$3$} (23);
	   \path[->] (3) edge node[above right] {$1$} (31);
	   \path[->] (3) edge node[above left] {$2$} (32);

	\end{tikzpicture}
	\caption{Poincaré algorithm as a simplicial system.}
	\label{poincare:sisy}
\end{figure}
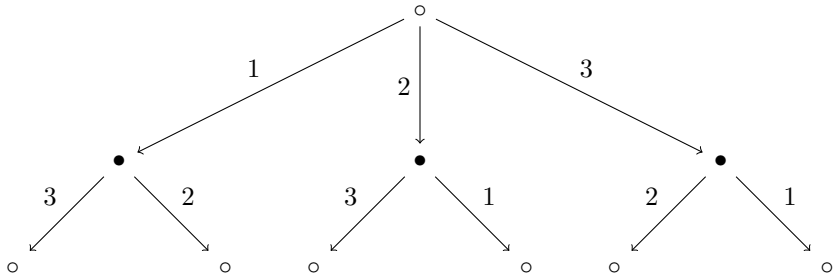

This map corresponds to the first return map of the simplicial system represented on Figure~\ref{poincare:sisy} to the white node (where all white nodes are identified).
The first step is determining which coordinate is the smallest of the three and subtracting it to the other two.
The second step is comparing the two initially largest coordinates and subtracting the smallest to the largest.
This is precisely describing Poincaré algorithm.

Notice that the induction associated to the subgraph $\graph_{\{1,2\}}$, where $G$ is the graph represented on Figure~\ref{poincare:sisy}, is equivalent to Rauzy induction on two intervals.
As for fully subtractive algorithms, this subgraph is stable by a result of Nogueira \cite{Nogueira95}.
For a higher dimension $n$ Poincaré algorithm this tree graph construction generalizes by starting with a vertex of degree $n$ with $n$ edges labeled by distinct letters and removing the outgoing edge of the ingoing label for each new vertex; when there is only one label left, we identify the vertex to the root.

Dimensions $n \ge 4$ are the only classical examples to our knowledge for which the criterion developed in Section~\ref{criterion} does not apply and which does not has obvious stable subgraphs.

\subsection{Other examples}
We first deal with examples that do not have full image.
Let $$I_1, \dots, I_n \subset \simplex$$ be all the different image sets of the domain on which the given algorithm is a linear map.
In the following examples these domains of definition correspond to the different cases depending on the order of the coordinates and will thus be indexed by the corresponding permutations.
Moreover, the image sets $I_k$ will form a finite cover of the set $\simplex$.\\

Let $\pi$ be the finite-to-one projection from the disjoint union of the sets $I_k$ to $\simplex$.
We will construct a simplicial system for which a first return to a given set of vertices of the win-lose induction map $\rauzyacc$ satisfies $\pi \circ \rauzyacc = F \circ \pi$ and thus has the same dynamical properties as $F$.\\

If $d$ is the dimension of the simplex $\simplex$, by assumption on simplex-splitting MCF, for all $k$, there exists a matrix in $\SL(d+1, \Z)$ that sends projectively $I_k$ to $\simplex$.
We make the further assumption that the inverse of these matrices are non-negative.
In the examples we consider, the image sets are a union of domain sets up to higher codimension subsets.
Consider the graph whose vertices are all the image sets and draw an edge between $I_k$ and $I_l$ if there is a domain contained in $I_k$ which is sent to $I_l$ by a matrix in $\SL(d+1, \Z)$.

\begin{remark*}
	If for some given MCF this condition is not met, one can try to divide the domains of definition in smaller piece.
\end{remark*}

\begin{proposition}
	If two non-negative matrices in $\SL(d+1, \Z)$ have the same projective action on the extremal points of $\simplex$ then they are equal.
\end{proposition}

\begin{proof}
	Let $v_1, \dots, v_{d+1}$ be the vectors defining the extremal points of $\simplex$.
	Assume the images of these vectors by the first matrix are $w_1, \dots, w_{d+1}$.
	For the second matrix they must be by assumption $\alpha_1 w_1, \dots, \alpha_{d+1} w_{d_1}$.
	Moreover, as the matrices are both non-negative of determinant $1$, we have $\prod_{k=1}^{d+1} \alpha_k = 1$, hence $\alpha_1 = \dots = \alpha_{d+1} = 1$.
\end{proof}

In particular, it is enough to describe the action of a simplicial system on the extremal points of its linear domains to fully characterize it.
As we are reduced to the full image case, it is enough to find a graph that splits each simplex $I_k$ into the domain subsimplices it contains and to connect the endpoints of this graph with the corresponding image sets.
This will define the right simplicial system up to permutation of the extremal points of the simplex.
Checking the action on the extremal points of the simplex will be dealt with in the following by discussing labeling of the length vectors coordinates.
This is in general straightforward but will have to be discussed further in the case of Selmer algorithm.

\begin{remark*}
	This issue can always been dealt with up to taking a finite number of copies of the image set with different labellings.
\end{remark*}

\subsubsection{Brun algorithms}
\label{brun}
The Brun algorithm, introduced by Brun in 1957, is described in dimension 3 by the map
$F: (x_1, x_2, x_3)\in\mathbb{R}^3_+ \to (x'_1, x'_2, x'_3),$
where if $\{i,j,k\} = \{1,2,3\}$ and $x_i > x_j > x_k$,
\begin{equation*}
	x'_i = x_i - x_j,\
	x'_j = x_j,\
	x'_k = x_k.
\end{equation*}

The definition domains of this map are given by the order of the coordinates and the action of the map on these domains is described by Figure~\ref{action:brun}.
The Figure gives the action on the extremal points up to permutation, to specify it let us remark that each small triangle is sent to the large one which has a common side with the small one and contains the central point of the simplex.\\

The image sets as introduced above are all the $6$ halves of the simplex which we will denote by the relation on two coordinates that define them.
They are represented on Figure~\ref{image:brun}.

\begin{figure}[h!]
	\begin{subfigure}{.5\textwidth}
		\hspace{-.5cm}
	  \begin{subfigure}{0.4\textwidth}
	     \begin{tikzpicture}[scale=.5]
	       \coordinate (a) at (0,0);
	       \coordinate (ac) at (2,0);
	       \coordinate (b) at (4,0);
	       \draw (a)node[below]{$(1,0,0)$} -- (ac) -- (b)node[below]{$(0,1,0)$} --++ (120:2)coordinate(bc) --++ (120:2)coordinate(c)node[above]{$(0,0,1)$} --++ (240:2)coordinate(cc) -- cycle;
	       \draw[name path=a--bc] (a) -- (bc);
	       \draw[name path=b--cc] (b) -- (cc);
	       \draw[name path=ac--c] (ac) -- (c);
	       \path [name intersections={of=a--bc and ac--c,by=e}];
	       \fill (a) circle (2pt);\fill (b) circle (2pt);\fill (c) circle (2pt);
	       \fill (ac) circle (2pt);\fill (bc) circle (2pt);\fill (cc) circle (2pt);\fill (cc) circle (2pt);
	       \path[pattern=north west lines,pattern color=black] (a)--(ac)--(e)--cycle;
	     \end{tikzpicture}
	  \end{subfigure}
	  {\large$\xrightarrow{F}$}%
	  \begin{subfigure}{0.1\textwidth}
	  \end{subfigure}
	  \begin{subfigure}{0.3\textwidth}
	     \begin{tikzpicture}[scale=.5]
	       \coordinate (a) at (0,0);
	       \coordinate (ac) at (2,0);
	       \coordinate (b) at (4,0);
	       \draw (a) -- (ac) -- (b) --++ (120:2)coordinate(bc) --++ (120:2)coordinate(c) --++ (240:2)coordinate(cc) -- cycle;
	       \draw[name path=a--bc] (a) -- (bc);
	       \draw[name path=b--cc] (b) -- (cc);
	       \draw[name path=ac--c] (ac) -- (c);
	       \path [name intersections={of=a--bc and ac--c,by=e}];
	       \fill (a) circle (2pt);\fill (b) circle (2pt);\fill (c) circle (2pt);
	       \fill (ac) circle (2pt);\fill (bc) circle (2pt);\fill (cc) circle (2pt);\fill (cc) circle (2pt);
	       \path[pattern=north east lines,pattern color=black] (a)--(b)--(bc)--cycle;
	     \end{tikzpicture}
	  \end{subfigure}
	  \caption{action on simplicial domains.}
	  \label{action:brun}
  \end{subfigure}
  \hfill
	\begin{subfigure}{.5\textwidth}
		\vspace*{.6cm}
	  \begin{subfigure}{0.2\textwidth}
	     \begin{tikzpicture}[scale=.5]
	       \coordinate (a) at (0,0);
	       \coordinate (ac) at (2,0);
	       \coordinate (b) at (4,0);
	       \draw (a)-- (ac) -- (b)--++ (120:2)coordinate(bc) --++ (120:2)coordinate(c)--++ (240:2)coordinate(cc) -- cycle;
	       \draw (a) -- (b) node[pos=.6, above]{$2\!>\!3$};
	       \draw (b) -- (c) node[pos=.55, left]{$3\!>\!2$};
	       \draw (a) -- (bc);
	     \end{tikzpicture}
	  \end{subfigure}
	  \hfill
	  \begin{subfigure}{0.2\textwidth}
		\begin{tikzpicture}[scale=.5]
	       \coordinate (a) at (0,0);
	       \coordinate (ac) at (2,0);
	       \coordinate (b) at (4,0);
	       \draw (a)-- (ac) -- (b)--++ (120:2)coordinate(bc) --++ (120:2)coordinate(c)--++ (240:2)coordinate(cc) -- cycle;
	       \draw (a) -- (b) node[pos=.4, above]{$1\!>\!3$};
	       \draw (a) -- (c) node[pos=.55, right]{$3\!>\!1$};
	       \draw (b) -- (cc);
	     \end{tikzpicture}
	  \end{subfigure}
	  \hfill
	  \begin{subfigure}{0.2\textwidth}
	     \begin{tikzpicture}[scale=.5]
	       \coordinate (a) at (0,0);
	       \coordinate (ac) at (2,0);
	       \coordinate (b) at (4,0);
	       \draw (a)-- (ac) -- (b)--++ (120:2)coordinate(bc) --++ (120:2)coordinate(c)--++ (240:2)coordinate(cc) -- cycle;
	       \draw (c) -- (ac) node[pos=.5, below=.7cm,right]{$ 2\!>\!1$};
	       \draw (c) -- (ac) node[pos=.5, below=.7cm,left]{$1\!>\!2 $};
	     \end{tikzpicture}
	  \end{subfigure}
		\vspace*{.5cm}
	  \caption{image domains.}
	  \label{image:brun}
     \end{subfigure}
\end{figure}

\begin{figure}[!h]
	\center
	\begin{tikzpicture}[shorten >=1pt,node distance=2cm,on grid,auto]
	   \node[state] (2>3) {$2>3$};

	   \node[state] (3>1) [right=of 2>3] {$3>1$};

	   \node[state] (1>2) [right=of 3>1] {$1>2$};

	   \node[state] (2>1) [below=of 1>2] {$2>1$};

	   \node[state] (1>3) [left=of 2>1] {$1>3$};

	   \node[state] (3>2) [left=of 1>3] {$3>2$};

	   \node[state, dashed] (2>3g) [right=of 1>2] {$2>3$};
	   \node[state, dashed] (3>2g) [right=of 2>1] {$3>2$};

	   \path[->] (2>3) edge  (3>1);
	   \path[->] (2>3) edge[loop above,looseness=6]  (2>3);
	   \path[->] (2>3) edge[bend right] (1>3);

	   \path[->] (3>1) edge  (1>2);
	   \path[->] (3>1) edge[loop above,looseness=6]  (3>1);
	   \path[->] (3>1) edge[bend right] (2>1);

	   \path[->,dashed] (1>2) edge  (2>3g);
	   \path[->] (1>2) edge[loop above,looseness=6]  (1>2);
	   \path[->,dashed] (1>2) edge[bend right] (3>2g);

	   \path[->] (2>1) edge  (1>3);
	   \path[->] (2>1) edge[loop below,looseness=6]  (2>1);
	   \path[->] (2>1) edge[bend right] (3>1);

	   \path[->] (1>3) edge  (3>2);
	   \path[->] (1>3) edge[loop below,looseness=6]  (1>3);
	   \path[->] (1>3) edge[bend right] (2>3);

	   \path[->,dashed] (3>2g) edge  (2>1);
	   \path[->] (3>2) edge[loop below,looseness=6]  (3>2);
	   \path[->,dashed] (3>2g) edge[bend right] (1>2);
	\end{tikzpicture}
	\caption{Combinatoric of Brun algorithm image domains.}
	\label{domains:brun}
\end{figure}
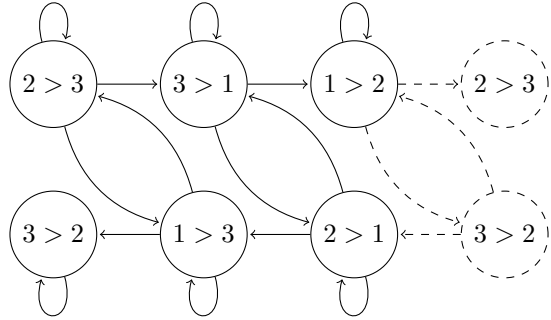

Each of these halves of the simplex is itself cut into three parts that are sent by Brun algorithm to three different halves.
The combinatoric of theses domains are represented in Figure~\ref{domains:brun}.
Where the dashed arrows and states are identified with the states of same label.\\

\begin{figure}[h!]
	\center
	\begin{tikzpicture}[shorten >=.5pt,node distance=1cm,on grid,auto, scale=.5]
	   \hspace{-2cm}
	   \node[state] (2>3) {$2>3$};
	   \node (bot23) [below=of 2>3] {$\bullet$};
	   \node (bbot23) [below=of bot23] {$\bullet$};

	   \node[state] (3>1) [right=4cm of 2>3] {$3>1$};
	   \node (bot31) [below=of 3>1] {$\bullet$};
	   \node (bbot31) [below=of bot31] {$\bullet$};

	   \node[state] (1>2) [right=4cm of 3>1] {$1>2$};
	   \node (bot12) [below=of 1>2] {$\bullet$};
	   \node (bbot12) [below=of bot12] {$\bullet$};

	   \node[state] (2>1) [below=of bbot12] {$2>1$};
	   \node[state] (1>3) [left=4cm of 2>1] {$1>3$};
	   \node[state] (3>2) [left=4cm of 1>3] {$3>2$};

	   \node (1>2g) [left=4cm of 2>3] {};
	   \node (bot12g) [below=of 1>2g] {};
	   \node (bbot12g) [below=of bot12g] {};
	   \node (2>1g) [below=of bbot12g] {};

	   \node (2>3g) [right=4cm of 1>2] {};
	   \node (bot23g) [below=of 2>3g] {};
	   \node (bbot23g) [below=of bot23g] {};
	   \node (3>2g) [right=4cm of 2>1] {};

	   \path[->] (2>3) edge[bend right=40]  node[left] {$3$} (bot23);
	   \path[->] (bot23) edge[bend right=40] node[right] {$2$} (2>3);
	   \path[->] (bot23) edge[bend right=20] node {$1$} (1>3);

	   \path[->] (3>1) edge[bend right=40]  node[left] {$1$} (bot31);
	   \path[->] (bot31) edge[bend right=40] node[right] {$3$} (3>1);
	   \path[->] (bot31) edge[bend right=20] node {$2$} (2>1);

	   \path[->] (1>2) edge[bend right=40]  node[left] {$2$} (bot12);
	   \path[->] (bot12) edge[bend right=40] node[right] {$1$} (1>2);
	   \path[->,dashed] (bot12) edge[bend right=20] node {$3$} (3>2g);
	   \path[->,dashed] (bot12g) edge[bend right=20] node {$3$} (3>2);

	   \path[->] (3>2) edge[bend right=40]  node[right] {$2$} (bbot23);
	   \path[->] (bbot23) edge[bend right=40] node[left] {$3$} (3>2);
	   \path[->,dashed] (bbot23) edge[bend right=20] node {$1$} (1>2g);
	   \path[->,dashed] (bbot23g) edge[bend right=20] node {$1$} (1>2);

	   \path[->] (1>3) edge[bend right=40]  node[right] {$3$} (bbot31);
	   \path[->] (bbot31) edge[bend right=40] node[left] {$1$} (1>3);
	   \path[->] (bbot31) edge[bend right=20] node {$2$} (2>3);

	   \path[->] (2>1) edge[bend right=40]  node[right] {$1$} (bbot12);
	   \path[->] (bbot12) edge[bend right=40] node[left] {$2$} (2>1);
	   \path[->] (bbot12) edge[bend right=20] node {$3$} (3>1);

	   \path[->] (2>3) edge node[above] {$1$} (3>1);
	   \path[->] (3>1) edge node[above] {$2$} (1>2);
	   \path[->,dashed] (1>2) edge node[above] {$3$} (2>3g);
	   \path[->,dashed] (1>2g) edge node[above] {$3$} (2>3);

	   \path[->] (1>3) edge node[below] {$2$} (3>2);
	   \path[->] (2>1) edge node[below] {$3$} (1>3);
	   \path[->,dashed] (3>2) edge node[below] {$1$} (2>1g);
	   \path[->,dashed] (3>2g) edge node[below] {$1$} (2>1);
	\end{tikzpicture}
	\caption{Brun algorithm as a simplicial system.}
	\label{graph:brun}
\end{figure}
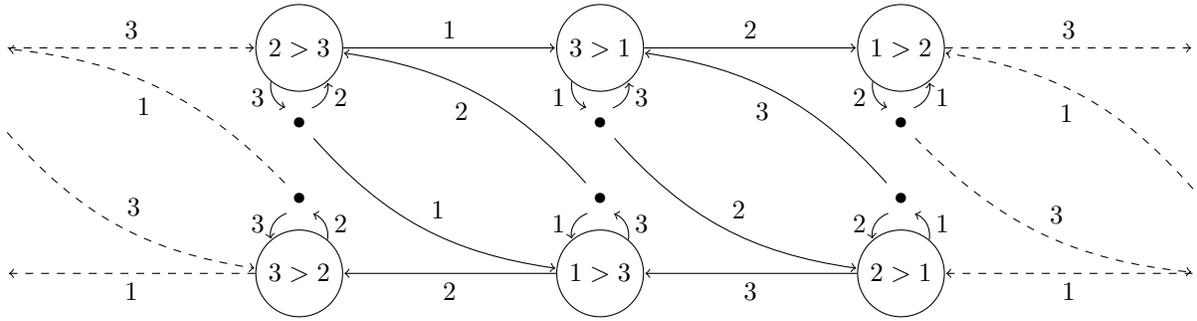

We can convert the three cuts in the simplex to a sequence of comparison between the three coordinates, as in Figure~\ref{graph:brun}.
Where the dashed arrow on left and right are identified with one another.

The actions on the three subsimplices in the image domains can be described by the graph in Figure~\ref{graph:brun}.\\

Following the arguments developed in the beginning of the section, we obtain the following proposition, which will generalize to higher dimensions.
\begin{proposition}
	Let $\rauzyacc$ be the first return map of the win-lose induction induced by the simplicial system defined on Figure~\ref{graph:brun} to the white circle vertices, then we have $\pi \circ \rauzyacc = F \circ \pi$.
\end{proposition}

In dimension $3$ there is an extra symmetry in the graph that enables us to define a conjugate algorithm on a simplicial system on $9$ vertices.
Indeed each top black vertex has its vertex labeled and pointing to a vertex exactly as for the black vertex at its bottom right.
Thus we can identify these pairs of black vertices to obtain a smaller graph (the one represented in Figure~\ref{graph:brunintro} defining the same algorithm.

Moreover we only need to check the \unstability\ for $2$ letters subgraphs.
In $G_{1,2}$, the strongly connected components are two loops around $1>2$ and $2>1$ which are clearly \unstable.
The same is true for any two letters and implies the following proposition.

\begin{proposition}
	Brun algorithm in dimension $3$ is simply connected and \rauzytype.\\
\end{proposition}

In particular Theorem~\ref{ergodicity} and Theorem~\ref{unstable} give an alternative proof of the following result as well as unicity of the ergodic measure and the fact that it induces the measure of maximal entropy on the canonical suspension.
\begin{theorem*}[\cite{Schweiger79}, \cite{Schweiger91}]
	The sorted Brun algorithm in dimension $3$ admits an invariant ergodic measure equivalent to Lebesgue measure.
\end{theorem*}
The sorted algoritm corresponds to the same maps on the subset of vectors with ordered coordinates composed with a permutation that order the coordinates afterwards.

This construction can be generalized to all dimensions.
For any $n \geq 2$, the Brun algorithm is defined by the map,
$$F: (x_1, \dots, x_n)\in\conen \to (x'_1, \dots, x'_n),$$
where for $\sigma \in \mathfrak S_n$ defined such that $x_{\sigma_1} > \dots > x_{\sigma_n}$,
\begin{align*}
	x'_{\sigma_1} &= x_{\sigma_1} - x_{\sigma_2} \\
	x'_{\sigma_i} &= x_{\sigma_i} \text{ for all } i \geq 2.
\end{align*}
The domains of definition depend again on the order of the coordinates.
They can be labeled by permutations in $\mathfrak S_n$ and will be denoted by $D_\sigma$ for any $\sigma \in \mathfrak S_n$.
For any $\sigma \in \mathfrak S_n$ the corresponding domain is sent bijectively by $F$ to the subsimplex defined by the equation $x'_{\sigma_2} > \dots > x'_{\sigma_n}$ which will be denoted by $I_\sigma$.
We change basis to have a simplex corresponding to a whole positive cone and for which the labels are compatible:
$$y_{\sigma_n} = x'_{\sigma_n}, \ y_{\sigma_{n-1}} = x'_{\sigma_{n-1}} - x'_{\sigma_n}, \ \dots \ , \ y_{\sigma_2} = x'_{\sigma_2} - x'_{\sigma_3} \ \text{ and } y_{\sigma_1} = x'_{\sigma_1}.$$

In $I_\sigma$, the coordinate $x'_{\sigma_1}$ can be in any position, in other words, $I_\sigma = \bigcup_{k = 1}^{n} D_{(1 \dots k) \sigma}$.
Thus the corresponding combinatoric graph has vertices from $I_\sigma$ to all $I_{(1 \dots k) \sigma}$ with $1 \le k \le n$.
Now the algorithm can be decomposed into first checking if $x'_{\sigma_1}$ is smaller than $x'_{\sigma_n}$, if so, $F$ sends the domain in $I_{(1\dots n) \sigma}$, otherwise, we check if $x'_{\sigma_1}$ is smaller than $x'_{\sigma_{n-1}}$, if so, $F$ sends the domain to $I_{(1\dots(n-1)) \sigma}$ and so on and so forth\dots.\\

One can check that this corresponds for a simplicial systems on coordinates $y$, to compare $y_{\sigma_n}$ and $y_{\sigma_1}$, then if $y_{\sigma_1}$ wins, compare $y_{\sigma_{n-1}}$ and $y_{\sigma_1}$ (since $y_{\sigma_1}$ will be equal to $x'_{\sigma_1}-x'_{\sigma_n}$), \dots.\\

\begin{figure}[!h]
	\center
	\begin{tikzpicture}[shorten >=1pt,node distance=2cm,on grid,auto]
	   \node[state] (sigma) {$I_\sigma$};
	   \node (2) [right=of sigma] {$\bullet$};
	   \node (3) [right=of 2] {$\dots$};
	   \node (4) [right=of 3] {$\bullet$};
	   \node (5) [right=of 4] {$\bullet$};
	   \node[state] (nsigma) [right=of 5] {$I_\sigma$};

	   \node (1sigma) [below=of sigma] {$I_{(1\dots n)\sigma} $};
	   \node (2sigma) [below=of 2] {$I_{\left(1 \dots (n-1)\right)\sigma}$};
	   \node (4sigma) [below=of 4] {$I_{\left(123\right)\sigma}$};
	   \node (5sigma) [below=of 5] {$I_{\left(12\right)\sigma}$};

	   \path[->] (sigma) edge node[left] {$\sigma_1$} (1sigma);
	   \path[->] (2) edge node[left] {$\sigma_1$} (2sigma);
	   \path[->] (4) edge node[left] {$\sigma_{1}$} (4sigma);
	   \path[->] (5) edge node[left] {$\sigma_{1}$} (5sigma);

	   \path[->] (sigma) edge node[above] {$\sigma_{n}$} (2);
	   \path[->] (2) edge node[above] {$\sigma_{n-1}$} (3);
	   \path[->] (3) edge node[above] {$\sigma_4$} (4);
	   \path[->] (4) edge node[above] {$\sigma_3$} (5);
	   \path[->] (5) edge node[above] {$\sigma_2$} (nsigma);
	\end{tikzpicture}
	\caption{Brun algorithm as a simplicial system in dimension $n$.}
	\label{graph:brun_multi}
\end{figure}
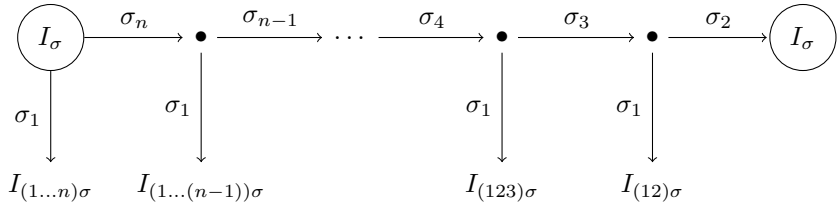

This description is giving us the corresponding vertices and labels between the image domains, it is represented on Figure~\ref{graph:brun_multi}.

To each image domain $I_\sigma$ corresponds a white circle vertex in the simplicial system.
As described in the introduction to the section and in the case of dimension 3, we relate this win-lose induction to Brun algorithm.
\begin{proposition}
	Let $\rauzyacc$ be the first return map of the win-lose induction induced by the simplicial system on Figure~\ref{graph:brun_multi} to the white circle vertices, then we have $\pi \circ \rauzyacc = F \circ \pi$.
\end{proposition}

\begin{proposition}
	Brun algorithm in all dimensions is simply connected and \rauzytype.
\end{proposition}

This also imply an alternative proof of another of Schweiger's theorems.
\begin{theorem*}[\cite{Schweiger00}]
	Sorted Brun algorithms in all dimensions admit an invariant ergodic measure equivalent to Lebesgue measure.
\end{theorem*}
Our proof actually shows something stronger for the unsorted algorithms.
\begin{theorem}
	Brun algorithms in all dimensions admit a unique invariant ergodic measure equivalent to Lebesgue measure which induces the unique measure of maximal entropy on its canonical suspension.
\end{theorem}

\begin{proof}
	The graph is clearly strongly connected and all labels in $\{1, \dots, n\}$ appear at least once.\\

	Let us denote by $\graph$ the graph in Figure~\ref{graph:brun_multi} and let $\largea$ be a non trivial subset of $\alphabet$.
	In the subgraph $\graph_\largea$, for a white circle vertex labeled by $I_\sigma$, the accelerated win-lose induction to the next white circle vertex acts on the permutation $\sigma$ by inserting $\sigma_1$ in front of the last losing label.
	Let us define $m := \max \{ i \geq 0 \mid \sigma([n-i+1, n]) \subset \largea \}$ then the set $M := \sigma([n-m+1,n]) \cap \largea$ is non-decreasing for the accelerated induction.
	Hence this set is invariant in a strongly connected component.

	Now consider a white circle vertex labeled $I_\sigma$ with an associated invariant set $M \subset \largea$.
	The intermediate black vertices all have two edges one of which is labeled by $\sigma_1$.
	Thus if $\sigma_1$ is not in $\largea$ then there are no vertices with more than one edge labeled in $\largea$.

	If $\sigma_1$ is in $\largea$, the connected component contains only the sequence of edges labeled by $\sigma_n, \sigma_{n-1}, \dots, \sigma_{n-m+1}$, since the last vertex in this sequence points either to a black vertex with a label in $\largea$ or to a white circle vertex labeled by $\sigma_1$ in $\largea$.
	For each of the vertices in the sequence, but the last, there are two edges labeled in $\largea$.
	Nevertheless, the edge labeled by $\sigma_1$ points to a white circle vertex for which the invariant set is $M \cup \{\sigma_1\}$.
\end{proof}

\subsubsection{Selmer algorithms}

\begin{figure}[h!]
	\center
	  \begin{subfigure}{0.25\textwidth}
	     \begin{tikzpicture}[scale=.5]
	       \coordinate (a) at (0,0);
	       \coordinate (ac) at (2,0);
	       \coordinate (b) at (4,0);
	       \draw (a)node[below]{$(1,0,0)$} -- (ac) -- (b)node[below]{$(0,1,0)$} --++ (120:2)coordinate(bc) --++ (120:2)coordinate(c)node[above]{$(0,0,1)$} --++ (240:2)coordinate(cc) -- cycle;
	       \draw[name path=a--bc] (a) -- (bc);
	       \draw[name path=b--cc] (b) -- (cc);
	       \draw[name path=ac--c] (ac) -- (c);
	       \path [name intersections={of=a--bc and ac--c,by=e}];
	       \fill (a) circle (2pt);\fill (b) circle (2pt);\fill (c) circle (2pt);
	       \fill (ac) circle (2pt);\fill (bc) circle (2pt);\fill (cc) circle (2pt);\fill (cc) circle (2pt);
	       \path[pattern=north west lines,pattern color=black] (a)--(ac)--(e)--cycle;
	     \end{tikzpicture}
	  \end{subfigure}
	  {\large$\xrightarrow{F}$}%
	  \begin{subfigure}{0.2\textwidth}
	  \end{subfigure}
	  \begin{subfigure}{0.25\textwidth}
	     \begin{tikzpicture}[scale=.5]
	       \coordinate (a) at (0,0);
	       \coordinate (ac) at (2,0);
	       \coordinate (b) at (4,0);
	       \draw (a) -- (ac) -- (b) --++ (120:2)coordinate(bc) --++ (120:2)coordinate(c) --++ (240:2)coordinate(cc) -- cycle;
	       \draw[name path=a--bc] (a) -- (bc);
	       \draw[name path=b--cc] (b) -- (cc);
	       \draw[name path=ac--c] (ac) -- (c);
	       \draw (ac) -- (bc);
	       \path [name intersections={of=a--bc and ac--c,by=e}];
	       \fill (a) circle (2pt);\fill (b) circle (2pt);\fill (c) circle (2pt);
	       \fill (ac) circle (2pt);\fill (bc) circle (2pt);\fill (cc) circle (2pt);\fill (cc) circle (2pt);
	       \path[pattern=north east lines,pattern color=black] (a)--(ac)--(bc)--cycle;
	     \end{tikzpicture}
	  \end{subfigure}
	  \caption{Action on simplicial domains.}
	  \label{action:selmer}
\end{figure}
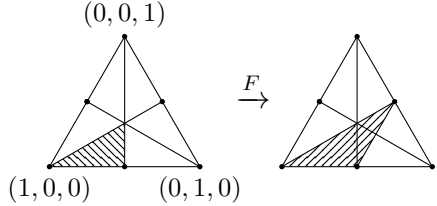

Introduced by Selmer in 1961 \cite{Selmer61}, the Selmer algorithm in dimension $3$ is defined by
$$F: (x_1, x_2, x_3)\in\Rp^3 \to (x'_1, x'_2, x'_3),$$
where if $\{i,j,k\} = \{1,2,3\}$ and $x_i > x_j > x_k$,
\begin{equation*}
	x'_i = x_i - x_k,\
	x'_j = x_j,\
	x'_k = x_k.
\end{equation*}
Figure~\ref{action:selmer} describes the action of Selmer algorithm on its simplicial domains.
Notice that unlike Brun algorithm, the image domains are not covering the simplicial domains defining the map.
This is related to the fact that the subsimplex $D$ defined by $x_i < x_j + x_k$ for all $\{i,j,k\} = \{1,2,3\}$ is an invariant attractive subset of this algorithm.\\

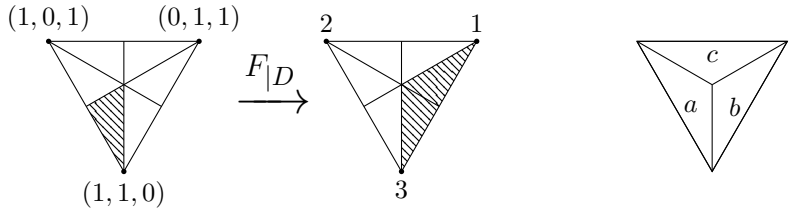
\begin{figure}[h!]
  \begin{subfigure}{0.25\textwidth}
	  \begin{tikzpicture}[scale=.5]
	       \coordinate (A) at (0,0);
	       \coordinate (Ac) at (2,0);
	       \coordinate (B) at (4,0);
	       \draw (A)node[above]{$(1,0,1)$} -- (Ac) -- (B)node[above]{$(0,1,1)$} --++ (-120:2)coordinate(Bc) --++ (-120:2)coordinate(C)node[below]{$(1,1,0)$} --++ (-240:2)coordinate(Cc) -- cycle;
	       \draw[name path=A--Bc] (A) -- (Bc);
	       \draw[name path=B--Cc] (B) -- (Cc);
	       \draw[name path=Ac--C] (Ac) -- (C);
	       \path [name intersections={of=A--Bc and Ac--C,by=E}];
	       \fill (A) circle (2pt);\fill (B) circle (2pt);\fill (C) circle (2pt);
	       \path[pattern=north west lines,pattern color=black] (C)--(E)--(Cc)--cycle;
	     \end{tikzpicture}
	  \end{subfigure}
	  {\LARGE$\xrightarrow{F_{|D}}$}%
	  \begin{subfigure}{0.25\textwidth}
		  \begin{tikzpicture}[scale=.5]
	       \coordinate (A) at (0,0);
	       \coordinate (Ac) at (2,0);
	       \coordinate (B) at (4,0);
	       \draw (A)node[above]{$2$} -- (Ac) -- (B)node[above]{$1$} --++ (-120:2)coordinate(Bc) --++ (-120:2)coordinate(C)node[below]{$3$} --++ (-240:2)coordinate(Cc) -- cycle;
	       \draw[name path=A--Bc] (A) -- (Bc);
	       \draw[name path=B--Cc] (B) -- (Cc);
	       \draw[name path=Ac--C] (Ac) -- (C);
	       \path [name intersections={of=A--Bc and Ac--C,by=E}];
	       \fill (A) circle (2pt);\fill (B) circle (2pt);\fill (C) circle (2pt);
	       \path[pattern=north west lines,pattern color=black] (C)--(B)--(E)--cycle;
	     \end{tikzpicture}
	\end{subfigure}
  \hfill \ \hfill
  \begin{subfigure}{0.3\textwidth}
	  \begin{tikzpicture}[scale=.5]
	       \coordinate (A) at (0,0);
	       \coordinate (Ac) at (2,0);
	       \coordinate (B) at (4,0);
	       \draw (A) -- (Ac) -- (B) --++ (-120:2)coordinate(Bc) --++ (-120:2)coordinate(C) --++ (-240:2)coordinate(Cc) -- cycle;
	       \draw[name path=A--Bc, opacity=0] (A) -- (Bc);
	       \draw[name path=B--Cc, opacity=0] (B) -- (Cc);
	       \draw[name path=Ac--C, opacity=0] (Ac) -- (C);
	       \draw (A) -- (E);
	       \draw (B) -- (E);
	       \draw (E) -- (C);
	       \path [name intersections={of=A--Bc and Ac--C,by=E}];
	       \draw (A) -- (B) node[pos=.5, below]{$c$};
	       \draw (B) -- (C) node[pos=.5, left]{$b$};
	       \draw (A) -- (C) node[pos=.5, right]{$a$};
	     \end{tikzpicture}
	  \end{subfigure}
  \caption{Action on the restriction and image domains.}
  \label{action:restriction}
\end{figure}

Restricted to $D$, the algorithm admits a simple description which was studied earlier by Mönkemeyer (see \cite{Monkemeyer54} and \cite{Panti08}).
In Figure~\ref{action:restriction}, we represent the action of the restriction of the Selmer algorithm on $D$, define new labels for a basis of the simplex $D$ and for image domains. On this domain, the vertex 2 is fixed and the central one is sent to 1.\\

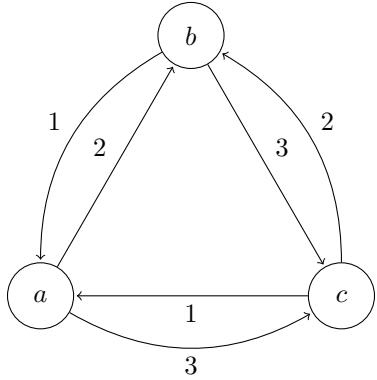
\begin{figure}[!htb]
	\center
	\begin{tikzpicture}[shorten >=1pt,node distance=4cm,on grid,auto]

	   \coordinate (a) at (0,0);
           \coordinate (b) at (2.0,3.464);
	   \coordinate (c) at (4,0);

           \node[state] (1) at (a) {$a$};
	   \node[state] (2) at (b) {$b$};
	   \node[state] (3) at (c) {$c$};

	   \path[->] (1) edge  node {$2$} (2);
	   \path[->] (2) edge  node {$3$} (3);
	   \path[->] (3) edge  node {$1$} (1);
	   \path[->, bend right] (2) edge  node[above left] {$1$} (1);
	   \path[->, bend right] (3) edge  node[above right]  {$2$} (2);
	   \path[->, bend right] (1) edge  node[below] {$3$} (3);

	\end{tikzpicture}
	\caption{Cassaigne algorithm as a simplicial system.}
	\label{cassaigne}
\end{figure}

This restriction of the algorithm is described by the graph in Figure~\ref{cassaigne}.
This graph is clearly \rauzytype \ and corresponds to Cassaigne algorithm given by the map
$$F_{|D}: (x_1, x_2, x_3)\in\Rp^3 \to
\left\{
\begin{array}{ll}
	(x_1-x_3, x_3, x_2) &\text{if } x_1 > x_3\\
	(x_2, x_1, x_3-x_1) &\text{if } x_3 > x_1
\end{array}
\right..
$$

The domains $a,b,c$ correspond to marking the permutation action of the Cassaigne algorithm on the coordinates of the vector.

\begin{proposition}
	Let $T$ be the win-lose induction induced by the simplicial system on Figure~\ref{cassaigne}, then we have $\pi \circ \rauzyacc = F_{|D} \circ \pi$.
\end{proposition}
As a straightforward consequence we obtain the following proposition.
\begin{proposition}
	Selmer algorithm in dimension $3$ restricted to $D$ is simply connected and \rauzytype.\\
\end{proposition}

Let us consider the generalization of this algorithm for $n \ge 3$,
$$F: (x_1, \dots, x_n)\in\conen \to (x'_1, \dots, x'_n),$$
where for $\sigma \in \mathfrak S_n$ defined such that $x_{\sigma_1} > \dots > x_{\sigma_n}$,
\begin{align*}
	x'_{\sigma_1} &= x_{\sigma_1} - x_{\sigma_n} \\
	x'_{\sigma_i} &= x_{\sigma_i} \text{ for all } i \geq 2.
\end{align*}
As for Brun algorithms, the domains of definition are labeled by $\mathfrak S_n$ and will be denoted by $D_\sigma$ for any $\sigma \in \mathfrak S_n$.
Similarly to dimension $3$, there is a stable subsimplex $D$ defined by the equations $x_{\sigma_1} < x_{\sigma_{n-1}} + x_{\sigma_n}$.
We consider the map $F_{|D}$ where $D_\sigma$ denotes the intersection of this set with $D$.\\

For any $\sigma \in \mathfrak S_n$ the domain $D_\sigma$ in $D$ is sent bijectively by $F_{|D}$ to the subsimplex defined by the equations $x'_{\sigma_2} > \dots > x'_{\sigma_n}$ and $x'_{\sigma_1} < x'_{\sigma_{n-1}}$ which will be denoted by $I_\sigma$.
In $I_\sigma$ the coordinate $x'_{\sigma_1}$ can either be in position $n-1$ or $n$, in other words, $I_\sigma = D_{(1 \dots (n-1)) \sigma} \cup D_{(1 \dots n) \sigma}$.
Thus the corresponding combinatoric graph has vertices pointing from $I_\sigma$ to $I_{(1 \dots (n-1)) \sigma}$ and from $I_\sigma$ to $I_{(1 \dots n) \sigma}$.\\

We first define a labeling for the basis which will help us keep track of the permutation of the extremal points of the simplex.
This is a generalization of what we did previously on Selmer algorithm in dimension 3.

The point for which all coordinates but one are equal to 1 and the other is equal to 0 is an extremal point of $D$ and is fixed by the algorithm.
We label each of these points by the label corresponding to its zero coordinate:
\[v_\alpha = 11 \dots 1 \underset{\alpha} 01 \dots 1.\]
This is what we did before in Figure~\ref{action:restriction}.
Now observe that $I_\sigma$ is the convex hull of $v_{\sigma_n}, v_{\sigma_{1}}, c$ and $w_k$ for $2 \le k \le n-2$, where $c$ is the point for which all coordinates are equal to $1$ and
$$w_k(i) =
\left\{
\begin{array}{ll}
	1 	      &\text{if } \ i = \sigma_2, \dots, \sigma_k \\
	\frac {1} {2} &\text{otherwise}
\end{array}
\right..
$$

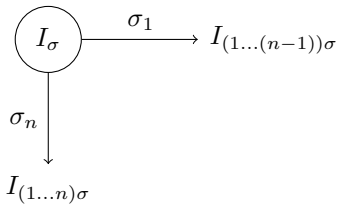
\begin{figure}[!h]
	\center
	\begin{tikzpicture}[shorten >=1pt,node distance=2cm,on grid,auto]
	   \node[state] (sigma) {$I_\sigma$};
	   \node (nsigma) [right=3cm of sigma] {$I_{\left(1\dots(n-1)\right)\sigma}$};

	   \node (1sigma) [below=of sigma] {$I_{(1\dots n) \sigma} $};

	   \path[->] (sigma) edge node[left] {$\sigma_n$} (1sigma);
	   \path[->] (sigma) edge node[above] {$\sigma_{1}$} (nsigma);
	\end{tikzpicture}
	\caption{Selmer algorithm as a simplicial system in dimension $n$.}
	\label{graph:selmer_multi}
\end{figure}

On each of these subsimplices the algorithm only compares coordinates in $v_{\sigma_1}$ and $v_{\sigma_n}$, thus these two labels are the only ones that matter.
In this labeling, $\sigma_1$ loses when $x_{\sigma_1} > x_{\sigma_{n}}$ and vice-versa, which may be counter-intuitive.
The graph for Selmer algorithm is thus described by Figure~\ref{graph:selmer_multi}.

\begin{proposition}
	Let $T$ be the win-lose induction induced by the simplicial system on Figure~\ref{graph:selmer_multi}, then we have $\pi \circ \rauzyacc = F_{|D} \circ \pi$.
\end{proposition}
\begin{proposition}
	Selmer algorithm restricted to $D$ in all dimensions is simply connected and \rauzytype.
\end{proposition}

In particular Theorem~\ref{ergodicity} and Theorem~\ref{unstable} give an alternative proof of the following result as well as unicity of the ergodic measure and the fact that it induces the measure of maximal entropy on the canonical suspension.
\begin{theorem*}[\cite{Monkemeyer54}, \cite{Schweiger00}]
	Sorted Selmer algorithms in all dimensions admit an invariant ergodic measure equivalent to Lebesgue measure.
\end{theorem*}
Recall that the sorted algorithms denote the same maps on vectors with ordered coordinates composed with a permutation that orders the coordinates afterwards.
Our proof actually shows something stronger.
\begin{theorem}
	Selmer algorithms in all dimensions admit a unique invariant ergodic measure equivalent to Lebesgue measure which induces the unique measure of maximal entropy on its canonical suspension.
\end{theorem}

\begin{proof}
	The graph is strongly connected since the permutation group is generated by the two cycles $(1 \dots n)$ and $(1 \dots (n-1))$.
	Moreover, all labels in $\{1, \dots, n\}$ appear at least once.\\

	Let us denote by $G$ the graph in Figure~\ref{graph:selmer_multi} and let $\largea$ be some non trivial subset of $\alphabet$.
	The property that $\sigma_n$ is in $\alphabet \setminus \largea$ is invariant in the subgraph $\graph_\largea$, so in a strongly connected component $\sigma_n$ is either always or never in $\largea$.

	If $\sigma_n$ is in $\alphabet \setminus \largea$ then for all vertices in the strongly connected component one of the two edges in not labeled in $\largea$.

	If $\sigma_n$ is in $\largea$, it remains so in the next step unless $\sigma_1$ is in $\alphabet \setminus \largea$.
	But at each step the numbers $\sigma_1, \dots, \sigma_{n-1}$ are shifted to the left in the permutation.
	Hence in less than $n$ steps, the permutation is such that $\sigma_1$ is in $\alphabet \setminus \largea$.
\end{proof}

\subsection{Rauzy Gaskets and Arnoux-Rauzy-Poincaré}
\label{gasket}

Following \cite{ArnouxStarosta13}, we define the Rauzy gasket in arbitrary dimension $n \ge 2$.
Let $C = \{ (x_1, \dots, x_n) \in \conen \mid x_j \le \sum_{i \ne j} x_i, \ \forall j \}$ and the Arnoux-Rauzy map,
$$F: (x_1, \dots, x_n)\in\conen\setminus C \to (x'_1, \dots, x'_n),$$
where for $\sigma \in \mathfrak S_n$ defined such that $x_{\sigma_1} > \dots > x_{\sigma_n}$,
\begin{align*}
	x'_{\sigma_1} &= x_{\sigma_1} - \sum_{i=2}^{n} x_{\sigma_i} \\
	x'_{\sigma_i} &= x_{\sigma_i} \text{ for all } i \geq 2.
\end{align*}

Consider the limit set,
$$ \gasket = \bigcap_{n\ge 0} F^{-n}(\conen \setminus C).$$
The Rauzy gasket is the intersection $\gasket \cap \simplex$, where
$$\simplex := \{(x_1, \dots, x_n) \mid \sum x_i = 1 \}.$$

Observe that in a simplicial system point of view, the map $F$ first splits the simplex depending on the order of the coordinates then for each ordering $\sigma \in \mathcal S_n$ sends the subsimplex defined by $x_{\sigma_1} > \sum_{i=2}^n x_{\sigma_i}$ to the whole simplex and is not defined on the other parts.
The graph will thus have two main parts: one connecting $I_\sigma$ states to $\tilde I_\sigma$ which will be the same as for Brun algorithm and another one which connects states $\tilde I_\sigma$ to $I_\sigma$ cutting out the parts on which the algorithm is not defined.

\begin{figure}[!h]
	\center
	\begin{tikzpicture}[shorten >=1pt,node distance=2cm,on grid,auto]
	   \node[state] (sigma) {$I_\sigma$};
	   \node (2) [right=of sigma] {$\bullet$};
	   \node (3) [right=of 2] {$\dots$};
	   \node (4) [right=of 3] {$\bullet$};
	   \node (5) [right=of 4] {$\bullet$};
	   \node[state] (nsigma) [right=of 5] {$\tilde I_\sigma$};

	   \node (1sigma) [below=of sigma] {$\tilde I_{(1\dots n)\sigma} $};
	   \node (2sigma) [below=of 2] {$\tilde I_{\left(1 \dots (n-1)\right)\sigma}$};
	   \node (4sigma) [below=of 4] {$\tilde I_{\left(123\right)\sigma}$};
	   \node (5sigma) [below=of 5] {$\tilde I_{\left(12\right)\sigma}$};

	   \path[->] (sigma) edge node[left] {$\sigma_1$} (1sigma);
	   \path[->] (2) edge node[left] {$\sigma_1$} (2sigma);
	   \path[->] (4) edge node[left] {$\sigma_{1}$} (4sigma);
	   \path[->] (5) edge node[left] {$\sigma_{1}$} (5sigma);

	   \path[->] (sigma) edge node[above] {$\sigma_{n}$} (2);
	   \path[->] (2) edge node[above] {$\sigma_{n-1}$} (3);
	   \path[->] (3) edge node[above] {$\sigma_4$} (4);
	   \path[->] (4) edge node[above] {$\sigma_3$} (5);
	   \path[->] (5) edge node[above] {$\sigma_2$} (nsigma);
	\end{tikzpicture}
	\caption{Part of the graph for Rauzy gasket connecting $I_\sigma$ to $\tilde I_\sigma$.}
	\label{graph:gasket}
\end{figure}
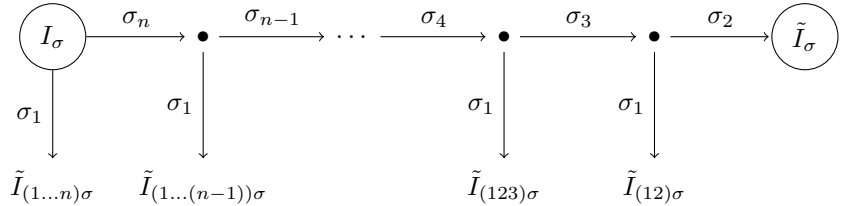

Now consider the compatible basis introduced for Brun algorithm
$$y_{\sigma_n} = x_{\sigma_n}, \ y_{\sigma_{n-1}} = x_{\sigma_{n-1}} - x_{\sigma_n}, \ \dots \ , \ y_{\sigma_2} = x_{\sigma_2} - x_{\sigma_3} \ \text{ and } y_{\sigma_1} = x_{\sigma_1}.$$
In this basis, the condition $x_{\sigma_1} < \sum_{i=2}^n x_{\sigma_i}$ is given by
$$y_{\sigma_1} < y_{\sigma_3} + 2 y_{\sigma_4} + \dots + (n-2) y_{\sigma_n}.$$
This is given by a graph with a sequence of edges from $\tilde I_{\sigma}$ to $I_\sigma$ labeled in the following order: $\sigma_3$, twice $\sigma_4$, \dots, $n-2$ times $\sigma_n$.
For each vertex in this sequence, starting with $\tilde I_\sigma$, there is an edge labeled by $\sigma_1$ and going out.\\

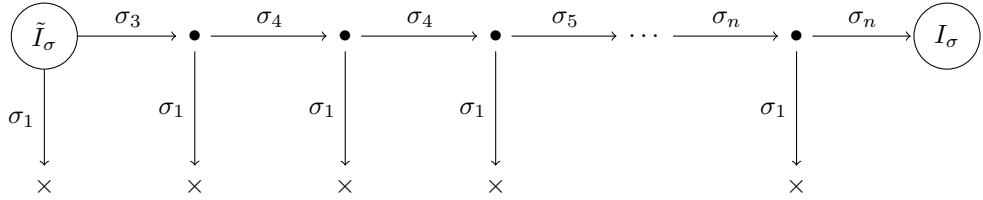
\begin{figure}[!h]
	\center
	\hspace*{-.5cm}
	\begin{tikzpicture}[shorten >=1pt,node distance=2cm,on grid,auto]
	   \node[state] (sigma) {$\tilde I_\sigma$};
	   \node (2) [right=of sigma] {$\bullet$};
	   \node (3) [right=of 2] {$\bullet$};
	   \node (4) [right=of 3] {$\bullet$};
	   \node (5) [right=of 4] {$\dots$};
	   \node (6) [right=of 5] {$\bullet$};
	   \node[state] (nsigma) [right=of 6] {$I_\sigma$};

	   \node (1sigma) [below=of sigma] {$\times$};
	   \node (2sigma) [below=of 2] {$\times$};
	   \node (3sigma) [below=of 3] {$\times$};
	   \node (4sigma) [below=of 4] {$\times$};
	   \node (6sigma) [below=of 6] {$\times$};

	   \path[->] (sigma) edge node[left] {$\sigma_1$} (1sigma);
	   \path[->] (2) edge node[left] {$\sigma_1$} (2sigma);
	   \path[->] (3) edge node[left] {$\sigma_1$} (3sigma);
	   \path[->] (4) edge node[left] {$\sigma_{1}$} (4sigma);
	   \path[->] (6) edge node[left] {$\sigma_{1}$} (6sigma);

	   \path[->] (sigma) edge node[above] {$\sigma_{3}$} (2);
	   \path[->] (2) edge node[above] {$\sigma_{4}$} (3);
	   \path[->] (3) edge node[above] {$\sigma_4$} (4);
	   \path[->] (4) edge node[above] {$\sigma_5$} (5);
	   \path[->] (5) edge node[above] {$\sigma_{n}$} (6);
	   \path[->] (6) edge node[above] {$\sigma_n$} (nsigma);
	\end{tikzpicture}
	\caption{Part of the graph for Rauzy gasket connecting $\tilde I_\sigma$ to $I_\sigma$.}
	\label{graph:gasket2}
\end{figure}

Let $\subgraph$ be the subgraph of the graph defined in Figure~\ref{graph:gasket} and Figure~\ref{graph:gasket2}.
From the construction it is clear that we have,
\begin{proposition}
	$\pi \left(\subparameter\right) = \gasket$.
\end{proposition}

Moreover, the subgraph $\subgraph$ is dynamically equivalent to the graph defining Brun algorithm (it can be accelerated to the Brun algorithm) since the only added edges are of degree one.
Thus we have,
\begin{theorem}
	The Rauzy gasket in any dimension $n \ge 3$ has Hausdorff dimension strictly smaller than $n-1$ and its canonical suspension flow has a unique measure of maximal entropy.\\
\end{theorem}

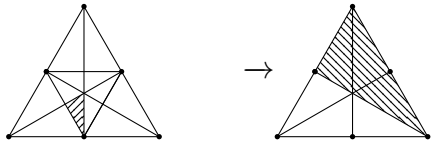
\begin{figure}[h!]
	\center
	  \begin{subfigure}{0.2\textwidth}
	  \end{subfigure}
	  \begin{subfigure}{0.25\textwidth}
	     \begin{tikzpicture}[scale=.5]
	       \coordinate (a) at (0,0);
	       \coordinate (ac) at (2,0);
	       \coordinate (b) at (4,0);
	       \draw (a) -- (ac) -- (b) --++ (120:2)coordinate(bc) --++ (120:2)coordinate(c) --++ (240:2)coordinate(cc) -- cycle;
	       \draw[name path=a--bc] (a) -- (bc);
	       \draw[name path=b--cc] (b) -- (cc);
	       \draw[name path=ac--c] (ac) -- (c);
	       \draw[name path=ac--bc] (ac) -- (bc);
	       \draw[name path=bc--cc] (bc) -- (cc);
	       \draw[name path=cc--ac] (cc) -- (ac);
	       \draw (ac) -- (bc);
	       \path [name intersections={of=a--bc and ac--c,by=e}];
	       \path [name intersections={of=cc--ac and a--bc,by=ec}];
	       \fill (a) circle (2pt);\fill (b) circle (2pt);\fill (c) circle (2pt);
	       \fill (ac) circle (2pt);\fill (bc) circle (2pt);\fill (cc) circle (2pt);\fill (cc) circle (2pt);
	       \path[pattern=north east lines,pattern color=black] (e)--(ec)--(ac)--cycle;
	     \end{tikzpicture}
	  \end{subfigure}
	  {\large$\xrightarrow{}$}%
	  \begin{subfigure}{0.25\textwidth}
	     \begin{tikzpicture}[scale=.5]
	       \coordinate (a) at (0,0);
	       \coordinate (ac) at (2,0);
	       \coordinate (b) at (4,0);
	       \draw (a) -- (ac) -- (b) --++ (120:2)coordinate(bc) --++ (120:2)coordinate(c) --++ (240:2)coordinate(cc) -- cycle;
	       \draw[name path=a--bc] (a) -- (bc);
	       \draw[name path=b--cc] (b) -- (cc);
	       \draw[name path=ac--c] (ac) -- (c);
	       \path [name intersections={of=a--bc and ac--c,by=e}];
	       \fill (a) circle (2pt);\fill (b) circle (2pt);\fill (c) circle (2pt);
	       \fill (ac) circle (2pt);\fill (bc) circle (2pt);\fill (cc) circle (2pt);\fill (cc) circle (2pt);
	       \path[pattern=north west lines,pattern color=black] (b)--(c)--(cc)--cycle;
	     \end{tikzpicture}
	  \end{subfigure}
	  \caption{Action of Poincaré algorithm on a subdomain of $C$.}
	  \label{action:poinca}
\end{figure}

Finally we remark that in dimension 3, the Poincaré algorithm acts on $C$ as described on Figure~\ref{action:poinca}.

\begin{figure}[!h]
	\center
	\begin{tikzpicture}[shorten >=1pt,node distance=2cm,on grid,auto]
	   \node[state] (sigma) {$\tilde I_\sigma$};
	   \node[state] (nsigma) [right=3cm of sigma] {$I_\sigma$};
	   \node (1sigma) [below=of sigma] {$I_{(123) \sigma}$};

	   \path[->] (sigma) edge node[left] {$\sigma_1$} (1sigma);
	   \path[->] (sigma) edge node[above] {$\sigma_{3}$} (nsigma);
	\end{tikzpicture}
	\caption{Connection for Arnoux-Rauzy-Poincaré in dimension 3.}
	\label{graph:arp}
\end{figure}
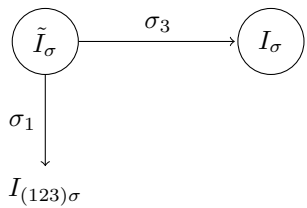

This gives us a natural way to describe Arnoux-Rauzy-Poincaré algorithm in dimension 3, consisting in applying the Arnoux-Rauzy map on $\conen \setminus C$ and the restriction of Poincaré map on $C$ (see \cite{BertheLabbe15}).
We only need to make the edges pointing to the hole vertex $\times$ from $\tilde I_\sigma$ point to $I_{(123) \sigma}$ as represented on Figure~\ref{graph:arp}.\\

As for Brun algorithm in dimension $3$, we only need to check the \unstability\ for two letter subgraphs, say $G_{1,2}$.
Here again the strongly connected components will be two loops around $1>2$ and $2>1$ formed by $3$ edges.
Which implies,
\begin{proposition}
	The Arnoux-Rauzy-Poincaré algorithm in dimension 3 is simply connected and \rauzytype.\\
\end{proposition}

Observe that the generalization of this algorithm to higher dimension will have more complicated combinatorics, since the images induced by the edges going out of the graph of Arnoux-Rauzy will for new sets of images.
Perhaps another more natural way to generalize this algorithm in the simplicial system point of view would be to connect all these edges to $I_{(1\dots n) \sigma}$.
This will again be a \unstable\ simplicial system.

\paragraph{Link with Baragar constants}

In this paragraph we show that constants computed by Baragar in \cite{Baragar98} correspond to $\kappa_\subgraph$ in Theorem~\ref{fractal:measure}.\\

By Proposition~\ref{jacobian:roof} the transfer operator $\mathcal L_s$ considered in \cite{GamburdMageeRonan19} (see formula (4.3) and Lemma 44 in their paper) is the Ruelle operator considered above
$$\mathcal L_s (f)(x) = \sum_{\rauzyacc(y)=x} e^{-s \cdot \roofacc(x)} f(y)$$
with an alternative acceleration of the algorithm.
Following \cite{GamburdMageeRonan19}, according to Ruelle--Perron--Frobenius theorem (Theorem 39) there exists a positive real number $\lambda_s$ which is the eigenvalue with the largest real part.
This eigenvalue satisfies $\lambda_s = e^{P_G(-s \cdot \roofacc)}$ in this setting of full shift with bounded variations (see Theorem~2 in \cite{Sarig03}).
As remarked after Proposition~43, the number $s$ such that $\lambda_s = 1$ corresponds to Baragar's constants.

Moreover these numbers coincide with the solution of the equation $P_G(-s \cdot \roofacc)$ which has an intrinsic definition as the entropy of the suspension flow on the fractal as described in Theorem~\ref{entropy_fractal}.\\

As a consequence of this remark, the computations in \cite{Baragar98} and Theorem~\ref{hausdorff} we have the following result.

\begin{theorem}
	If $\mathcal G^d$ denotes the Rauzy gasket in dimension $d$, we have the bounds
\begin{align*}
	\dim_H(\gasket^2) &< 1.825 \\
	\dim_H(\gasket^3) &< 2.7 \\
	\dim_H(\gasket^4) &< 3.612 \ .
\end{align*}
And for $d$ going to infinity,
\begin{equation*}
	\dim_H(\gasket^d) < d-1 + \frac {\log d} {\log 2 \cdot (d+1)} + o(d^{-1.58}).
\end{equation*}
\end{theorem}

This stengthens and generalizes the only previous known bound proved in \cite{AvilaHubertSkripchenko16a} to be $\dim_H(\gasket^2) < 2$.

\section*{Acknowledgments}
I am very grateful to Valérie Berthé, Vincent Delecroix, Sébastien Gouëzel, Pascal Hubert, Malo Jezequel and Sasha Skripchenko for multiple discussions related to this project as well as to Paul Mercat and Thierry Coulbois for their careful reading of early versions of this articles.
This work was done during my stay at MPIM in Bonn and MSRI in Berkeley, I am grateful for the great working environment they have provided.\\

This material is based upon work supported by the National Science Foundation under Grant No. 1440140, while the author was in residence at the Mathematical Sciences Research Institute in Berkeley, California, during the semester of fall 2019 \textit{Holomorphic differentials in Mathematics and Physics}.\\

This work was supported by the Agence Nationale de la Recherche through the project Codys (ANR- 18-CE40-0007).

\newpage
\bibliography{all.bib}
\bibliographystyle{alpha}
\end{document}